\documentclass[12pt,a4paper,reqno]{amsart}
\usepackage{amsmath, amssymb, verbatim,amsthm,amsbsy, a4,amsfonts,latexsym,
amscd,amsopn,bbm,stmaryrd,wasysym,
mathbbol, srcltx}
\usepackage[latin1]{inputenc}
\usepackage{color}              
\usepackage{epsfig}
\newtheorem{theorem}{Theorem}[section]
\newtheorem{lemma}[theorem]{Lemma}
\newtheorem{corollary}[theorem]{Corollary}
\newtheorem{proposition}[theorem]{Proposition}

\newtheorem{definition}[theorem]{Definition}

\newenvironment{Proof}{\removelastskip\par\medskip
\noindent{\em Proof.} \rm}{\penalty-20\null\hfill$\square$\par\medbreak}

\numberwithin{equation}{section}

\begin{document}
\hoffset  0.25  true cm
\voffset  1.5 true cm
\hsize=15.5 true cm
\vsize=21 true cm


\renewcommand\Re{\operatorname{Re}}
\renewcommand\Im{\operatorname{Im}}

\newcommand\Fr{\operatorname{Fr}}
\newcommand{\Cal}{\mathcal}
\newcommand{\ffrac}{\frac}

\newcommand{\llapp}[1]{\hbox{\llap{$\dsize{#1}$}}}
\newcommand{\rlapp}[1]{\hbox{\rlap{$\dsize{#1}$}}}
\renewcommand{\(}{ \left(}
\renewcommand{\)}{\right)}
\newcommand{\E}{{\,{\bold E}\,}}
\renewcommand{\P}{\,{\bold P}\,}
\newcommand{\pr}{\,{\bold P}\,}
\renewcommand{\a}{\alpha}
\newcommand{\be}{\beta}
\newcommand{\g}{\gamma}
\newcommand{\GA}{\Gamma}
\newcommand{\la}{\lambda}
\newcommand{\LA}{\Lambda}
\newcommand{\de}{\delta}
\newcommand{\om}{\omega}
\newcommand{\DE}{\Delta}
\newcommand{\cov}{\operatorname{cov}}
\newcommand{\p}{\varphi}
\newcommand{\<}{\bigl <}\newcommand{\8}{\left<}\newcommand{\9}{\right>}
\renewcommand{\>}{\bigr>}
\newcommand{\ovln}[1]{\,{\overline{\!#1}}}
\newcommand{\ov}{\overline}
\renewcommand{\L}{\Cal L}
\newcommand\N{\Cal N}
\newcommand\e{\varepsilon}
\newcommand\vk{\varkappa}
\renewcommand\t{\tau}
\renewcommand\th{\theta}
\newcommand\PM{^{\2\prime}}
\newcommand{\BT}[1]{\Cal B_d(#1)}
\newcommand{\BBBT}[1]{\Cal B_d^*(#1)}

\newcommand{\BD}[2]{\Cal B_{#2}(#1)}
\newcommand{\define}[2]{\newcommand{#1}{#2}}
\define\BB{\Cal B_d(\tau)}
\define\BBB{\Cal B_d^*(\tau)}
\newcommand{\AT}[1]{\Cal A_d(#1)}

\newcommand{\AD}[2]{\Cal A_{#2}(#1)}
\define\A{\Cal A_d(\tau)}
\newcommand{\ATq}[1]{\Cal A_d'(#1,\de,\rho)}

\newcommand{\ADq}[2]{\Cal A_{#2}'(#1,\de,\rho)}
\define\Aq{\Cal A_d'(\tau,\de,\rho)}

\newcommand{\Atdr}[3]{\Cal A_{d}'(#1,#2,#3)}
\newcommand{\ATQ}[1]{\Cal A_d''(#1,\de,\rho)}
\newcommand{\ADQ}[2]{\Cal A_{#2}''(#1,\de,\rho)}
\define\AQ{\Cal A_d''(\tau,\de,\rho)}

\newcommand{\ATDR}[3]{\Cal A_{d}''(#1,#2,#3)}
\newcommand{\ATQZ}[1]{\Cal A_d^*(#1,\de,\rho)}

\newcommand{\ADQZ}[2]{\Cal A_{#2}^*(#1,\de,\rho)}
\define\AQZ{\Cal A_d^*(\tau,\de,\rho)}

\newcommand{\ATDRZ}[3]{\Cal A_{d}^*(#1,#2,#3)}
\newcommand{\ADRZ}[4]{\Cal A_{#4}^*(#1,#2,#3)}
\newcommand{\Duu}[3]{{\<D_{#1}{#2},{#3}\>}}
\newcommand{\Dpuu}[3]{{\<D\PM_{#1}{#2},{#3}\>}}
\newcommand{\Dhuu}[3]{{\<D_{#1}(h)\,{#2},\,{#3}\>}}
\newcommand{\Dhu}[2]{{\<D_{#1}(h)\,{#2},\,{#2}\>}}
\newcommand{\Dwhu}[3]{{\<D_{#1}\bigl(#3\bigr)\4{#2},\,{#2}\4\>}}
\newcommand{\D }[1]{{\<D{#1},{#1}\>}}
\newcommand{\Dhuv}[2]{{\<D(h)\4{#1},{#2}\>}}
\define\si{\sigma}
\newcommand{\Duv}[2]{{\<D{#1},{#2}\>}}
\newcommand{\Du}[2]{{\<D_{#1}{#2},{#2}\>}}
\newcommand{\Dpu}[2]{{\<D\PM_{#1}{#2},{#2}\>}}
\newcommand{\norm}[1]{\left\|#1\right\|}
\newcommand{\nnorm}[1]{\bgl\|#1\bgr\|}
\newcommand{\nm}[1]{\bgl|#1\bgr|}
\newcommand{\nrm}[1]{|\2#1\2|}

\newcommand{\nnnorm}[1]{\left\|\2\smash{#1}\2\right\|}
\newcommand{\mnorm}[1]{\left|\2\smash{#1}\2\right|}
\newcommand{\snorm}[1]{\bigl\|\2#1\2\bigr\|}
\define\ssqrt{^{1/2}}
\define\msqrt{^{-1/2}}
\define\me{^{-1}}
\define\htau{\norm h\tau}
\define\ttau{\norm t\tau}
\newcommand{\bed}[1]{\bold1_{{#1}}}             
\newcommand{\ind}[1]{\bold I\bgl\{{#1}\bgr\}}             
\define\eht{\bigl(\2\e+\htau\2\bigr)}
\newcommand{\weht}[1]{\bigl(\2\e+\nnnorm{#1}\t\2\bigr)}
 \newcommand{\leei}[2]{\log\E e^{i\8\right.\!{#1},{#2}\!\left.\9}}
 \newcommand{\eei}[2]{\E e^{i\8\right.\!{#1},{#2}\!\left.\9}}
 \newcommand{\lee}[2]{\log\E e^{\8\right.\!{#1},{#2}\!\left.\9}}
 \newcommand{\ee}[2]{\E e^{\8\right.\!{#1},{#2}\!\left.\9}}
\newcommand{\xxxi}[1]{{\overline\xi}_{#1}(h)}
\define\xh{{\overline\xi}(h)}
\newcommand{\oxh}[2]{{\overline\xi_{#1}}(h_{#2})}
\newcommand{\ooxh}[2]{{\overline\xi{}_{#1}\PM}(h'_{#2})}
\newcommand{\xxi}[1]{\xi_#1}
\define\eq{{}}
\define\wth{\widetilde h}
\define\wthk{\widetilde h_k}
\newcommand{\wthj}[1]{\widetilde h_{#1}}
\define\tht{{\th\4\htau}}
\define\twhtk{{\th\4\nnnorm{\wthk}\tau}}
\newcommand{\twhtj}[1]{{\th\4\nnnorm{\wthj#1}\tau}}
\define\wtht{{\th\4\wthk\tau}}
\define\ttt{{\th\4\ttau}}
\define\Rd{\Bbb R^d}
\define\Rk{\bold R^k}
\define\Rdm{\Bbb R^{d-1}}
\define\Cd{\Bbb C^d}
\define\Bd{\frak B_d}
\define\Gd{\frak G_d}
\define\Fd{\frak F_d}
\define\xwth{{\overline\xi}(\wth)}
\newcommand{\xwthj}[1]{{\overline\xi_{{#1}}}(\wthj{#1})}
\newcommand{\xhpj}[2]{{\overline\xi{}\PM_{{#1}}}(h'_{#2})}
\define\tdsi{{\t\4d^{3/2}\!\big/\si}}
\define\ftdsi{\ffrac{\t\4d^{3/2}}{\si}}
\define\detD{(2\4\pi)^{d/2}\,(\det D)\ssqrt}
\newcommand{\detDj}[1]{(2\4\pi)^{d/2}\,(\det D_{#1})\ssqrt}
\define\wt{\widetilde}
\define\wh{\widehat}
\define\RE{\operatorname{Re}}
\define\IM{\operatorname{Im}}
\define\Var{\operatorname{Var}}
\define\jd{{j=1,\dots,d}}
\newcommand{\Y}[3]{{#1}_{#2}^{[#3]}}
\newcommand{\YY}[3]{{#1}_{#2}^{(#3)}}
\newcommand{\UU}[3]{\wt{#1}_{#2}^{(#3)}}
\newcommand{\YYY}[3]{{#1}_{#2}^{\{#3\}}}
\define\4{\kern1pt}
\define\gm{\4|\4}
\define\bgm{\,\big|\4}
\define\2{\kern.5pt}
\define\6{\phantom0}
\newcommand{\bgr}[1]{\4\bigr#1}
\newcommand{\bgl}[1]{\bigl#1\4}
\define\3{\kern-10pt}
\newcommand{\7}[1]{_{(#1)}}
\define\cdt{\4\cdot\4}
\newcommand{\R}[2]{\Cal R_d\bgl(#1,#2\bgr)}
\newcommand{\RB}[2]{\Cal R_d\Bgl(#1,#2\Bgr)}
\newcommand{\RBM}[2]{\Cal R_{d-1}\Bgl(#1,#2\Bgr)}
\define\half{{}^1\kern-.5pt\!\big/\!\2{}_2}
\newcommand{\sfrac}[2]{{}^{#1}\kern-.5pt\!\big/\!\kern.5pt{}_{#2}}
\newcommand{\pfrac}[2]{{#1}\big/{#2}}
\newcommand{\zfrac}[2]{\raise.5pt\hbox{\eightpoint$\dfrac{\,#1\,}{\,#2\,}$}}
\newcommand{\Bgr}[1]{{\eightpoint\4\raise.5pt\hbox{{$\Bigr#1$}}}}
\newcommand{\Bgl}[1]{{\eightpoint\raise.5pt\hbox{{$\Bigl#1$}}\4}}
\newcommand{\step}[1]{\raise9pt\hbox{\eightpoint$\tsize{#1}$}}
\renewcommand{\le}{\leqslant}
\renewcommand{\ge}{\geqslant}
\define\ph{\phantom}
\newcommand{\LL}[1]{\noindent\llapp{{[#1]{\kern10pt }}}\kern-2pt}
\define\ONE{\text{1\kern-2.5pt l}}
\define\bR{\bold R}
\define\BR{\bold R}
\define \CG{\Cal G}
\newcommand{\cit}[1]{[\2#1\2]}
\renewcommand{\=}{\overset{ \text{def} }\to =}
\define\na{\,\, {\raise.4pt\hbox{$\shortmid$}}{\hskip-2.0pt\to}\, \, }


\title[Free Probability Theory]{ The~Arithmetic of Distributions in Free Probability Theory} 

\author{ G. P. Chistyakov$^{1,2,3}$ and F. G\"otze$^{1,3}$}
\thanks{$^3$Research supported by SFB 701.
Partly supported by INTAS grant N 03-51-5018}


\date{September 2010}
\keywords{  Free random variables, Cauchy transforms, free convolutions,
Delphic semigroups, Khintchine's theorems}

\subjclass{
Primary 46L50}  
\maketitle
\begin{center}

$^1$University of Bielefeld,\\
$\phantom 0$\\
$^2$Institute for Low Temperature Physics and Engineering,\\
Kharkov\\
\end{center}
\markboth{ Free convolutions }{G. P. Chistyakov and  F. G\"otze }
\leftmark
\rightmark
\begin{abstract}
We give an~analytical approach to the~definition of additive and
multiplicative free convolutions which is based on 
the~theory of Nevanlinna and of Schur functions. We consider the set~of
probability distributions as a~semigroup $\bold M$ equipped with 
the~operation of free convolution and prove a~Khintchine type theorem 
for the  factorization of elements of this semigroup.
An~element of $\bold M$ contains either indecomposable (``prime'') 
factors or it belongs to a~class, say $I_0$, of distributions without
indecomposable factors. In contrast to the~classical convolution semigroup
in  the~free additive and multiplicative convolution semigroups 
the~class $I_0$ consists of units (i.e. Dirac measures) only. Furthermore
we show that the~set of indecomposable elements is dense in $\bold M$.

\end{abstract}

\maketitle

\section{Introduction} 

In recent years a~larger number of papers has been devoted 
to applications and extensions of the~definition
of free convolution of measures introduced by D. Voiculescu.
The~key concept of this definition is the~notion of freeness,
which can be interpreted as a~kind of independence for non-commutative
random variables. As in the~classical probability where the~concept of
independence gives rise to the~classical convolution, the~concept
of freeness leads to a~binary operation on the~probability measures
on the~real line, called free convolution. As one might expect 
there are many classical results for sums of independent random variables 
having a~counterpart in this theory, such as the~law of large numbers,
the~central limit theorem, the~L\'evy-Khintchine formula and others.
We refer to Voiculescu, Dykema and Nica~\cite{VoDyNi:1992} (1992) 
and to Hiai and Petz~\cite{HiPe:2000} (2000) for 
an~introduction to these topics. One of the~main problems in dealing with free
convolution is that its definition is rather indirect. In the~first
part of this paper we propose an~analytical approach to the~definition of
additive and multiplicative free convolution. This approach which develops
an analytic method due to Maassen~\cite{Maa:1992} 
is based on the~classical theory of Nevanlinna and Schur functions 
and allows us to give a~direct definition of free convolutions
(see Theorem~2.1, Theorem~2.4, Theorem~2.7) by purely analytic methods.
This approach shows the~equivalence of the  ``characteristic function'' approach
and the probabilistic approach for  additive and multiplicative free convolutions.
Note that Bercovici and Belinschi~\cite{BelBe:2007a} gave a~related approach using other analytic methods.

In the~second part of the~paper we study the~arithmetic structure of
the~Voi\-culescu semigroups  $(\mathcal M,\boxplus )$, 
$(\mathcal M_+,\boxtimes)$, and $(\mathcal M_*,\boxtimes )$ of probability 
measures on $\mathbb R$ with additive free convolution ($\boxplus$), 
on  $\mathbb R_+$ and on $\mathbb T$ with multiplicative free convolution 
($\boxtimes$), respectively. This subject had its origin
in the~work of Khintchine on the~convolution semigroup $(\mathcal M,*)$
of probability measures on the~real line. He derived for this semigroup
the~three basic theorems listed in Section~2. Kendall~\cite{Ken:1967} (1967),
\cite{Ken:1968} (1968)  
introduced the~so called Delphic semigroups, which are commutative 
topological semigroups satisfying the~central limit theorem for 
triangular arrays. Their~arithmetic is similar to 
the~convolution arithmetic of probability measures on $\mathbb R$.
A~characteristic feature of all these semigroups is the~presence of
{\it infinitely divisible} (i.d.) elements, which for every positive
integer $n$ may be represented as $n$-th power of some element of
the~semigroup.

In any Delphic semigroup there are three classes of its elements:
\begin{itemize}
\item
the~indecomposable or {\it simple} elements, which have no factors
besides themselves and the~identity, (a~set we shall denote by ``$S$'');

\item
the~elements which are decomposable and have an~indecomposable factor,
(a~set we shall denote by ``$D$'');

\item
the~infinitely divisible elements which have no indecomposable 
factors, (a~set we shall denote by ``$I_0$'').
\end{itemize}

It turns out that  a~lot of important semigroups, in particular $(\mathcal M,*)$,
are Delphic or almost Delphic.

It is convenient to formulate the~Delphic hypothesis rather restrictively,
and then to show that semigroups like $(\mathcal M,*)$ are 'almost' Delphic,
that means they satisfy these hypothesis with nonessential modifications.
In this context Davidson \cite{Dav:1968(I)}--\cite{Dav:1969} (1968), (1969)
introduced the~concept of an~hereditary sub-semigroup 
to verify that semigroups are almost (or properly) Delphic. 
Using a~multivariate analytic description
of free convolutions we show that the~Voiculescu semigroups  
$(\mathcal M,\boxplus)$, $(\mathcal M_+,\boxtimes)$, and $(\mathcal M_*,\boxtimes)$  
are close by the~structure to almost Delphic semigroups (but are not
almost Delphic semigroups).
Using some ideas of Khintchine, Kendall and Davidson, we deduce the~three
basic Delphic theorems for these semigroups. One of them states that each
element of a~Delphic semigroup may be written as a~product of a~countable
number of indecomposable elements and an~element of $I_0$, so a~knowledge
of $I_0$ and the~set of indecomposable elements is essential for 
the~arithmetic of the~semigroup. As a~consequence we show that in 
the~Voiculescu semigroups the~class $I_0$ consists of Dirac measures
$\delta_a$ only and for each semigroup there is a~dense set of 
indecomposable elements.

As another consequence of this approach we obtain an~analogue of 
Khintchine's limit theorem
in Voiculescu's semigroups.

The~paper is organized as follows. In Section~2 we discuss the~results
of the~paper. In Section~3 we collect auxiliary results from complex
analysis, free probability.
In Sections~4 and 5 
we prove the~necessary analytical results for our 
approach to free convolutions. In Section~6 we prove 
basic Delphic theorems for the~Voiculescu semigroups, in Section~7
we describe the~class $I_0$, and in Section~8 we describe some dense 
classes of indecomposable elements in these semigroups. 

Note that the~contents of this paper is close to the~contents of \cite{ChG:2005a}
with some changes. We changed the~text of the Introduction, added the~assertion (3)
in Proposition~\ref{3.11pro}, added some details in the proof of Lemma~\ref{4.2l}. We corrected
the text of Section 6. Here  we omitted some  propositions about the reduction to theory of 
Delphic semigroups which were not correct. In additional arguments  in Section~6 still based on the ideas of 
\cite{ChG:2005} we now only use  the hereditary property and the well-known methods
of Delphic semigroups for the proof of Theorems 2.10--2.12. In Section~8 we give some new 
free indecomposability conditions for probability measures using the~arguments of \cite{ChG:2005a}. 

{\bf Acknowledgment}.
We would like to thank Leonid Pastur who has drawn our attention
 to the~problem of the~arithmetic of Voiculescu  semigroups. 

\section{Results} 

Denote by $\mathcal M$ the~family of all Borel probability measures
(p-measures for short)
defined on the~real line $\mathbb R$. On the~set $\mathcal M$ there are
defined two associative composition laws  denoted $*$ and $\boxplus$.
Let $\mu_1,\mu_2\in\mathcal M$. The~measure $\mu_1*\mu_2$ is the~classical 
convolution of $\mu_1$ and $\mu_2$. In probabilistic terms, $\mu_1*\mu_2$
is the~probability distribution of $X+Y$, where $X$ and $Y$ are
(commuting) independent random variables with distributions $\mu_1$ 
and $\mu_2$ respectively. The~measure $\mu_1\boxplus\mu_2$ denotes the~free
(additive) convolution of $\mu_1$ and $\mu_2$ introduced by 
Voiculescu~\cite{Vo:1986} (1986)
for compactly supported measures. The~notion of free convolution 
was extended by Maassen~\cite{Maa:1992} (1992) to measures with 
finite variance 
and by Bercovici and Voiculescu~\cite{BercVo:1993} (1993) 
to all measures in $\mathcal M$.
Here, $\mu_1\boxplus\mu_2$ may be considered as the~probability 
distribution of $X+Y$, where $X$ and $Y$ are free random variables 
with distributions $\mu_1$ and $\mu_2$, respectively. 
For positive random variables and for random variables with values
on $\mathbb T$ we consider multiplicative convolutions as well and their
free analogues of multiplicative convolutions which were introduced 
by Voiculescu~\cite{Vo:1987} (1987). 

In this section we give an~analytical approach to the~definition of
$\mu_1\boxplus\mu_2$ which extends Maassen's definition.
Furthermore, we shall present an~analytical approach to the~definition of 
multiplicative free convolution $\boxtimes$.

Let $\mathbb C^+\,(\mathbb C^-)$ denote the~open upper (lower) half of
the~complex plane. If $\mu\in\mathcal M$,
then its Cauchy transform
\begin{equation}\label{2.1}
G_{\mu}(z)=\int\limits_{-\infty}^{\infty}\frac {\mu(dt)}{z-t},
\qquad z\in\Bbb C^+.
\end{equation}

Following Maassen~\cite{Maa:1992} and Bercovici and 
Voiculescu~\cite{BercVo:1993},
in the~sequel we will consider the~{\it reciprocal Cauchy transform}
\begin{equation}\label{2.2}
F_{\mu}(z)=\frac 1{G_{\mu}(z)}.
\end{equation}
Let $\mathcal F$ denote the~corresponding class of reciprocal Cauchy
transforms of all $\mu\in\mathcal M$.
This class admits the~following simple description. 


The~class $\mathcal F\subset \mathcal N$ of reciprocal Cauchy transforms of
p-measures introduced above coincides with the~subclass of 
Nevanlinna functions $F\in\mathcal N$ such that 
$F(z)/z\to 1$ as $z\to \infty$ nontangentially to 
$\mathbb R$ (i.e., such that $|\RE z|/\IM z$ stays bounded). See Section~3.
Recall that the~class of all analytic functions 
$F:\mathbb C^+\to\mathbb C^+\cup\mathbb R$, say
$\mathcal N$, is called {\it Nevanlinna class}.

This implies that $F_{\mu}$ has certain invertibility properties.
To be precise, for two numbers $\alpha>0,\beta>0$ we set
$$
\Gamma_{\alpha}=\{z=x+iy\in\Bbb C^+:|x|<\alpha y\}\quad\text{and}\quad
\Gamma_{\alpha,\beta}=\{z=x+iy\in\Gamma_{\alpha}:y>\beta\}.
$$
Then, by relation (\ref{3.4}) of Section~3, for every $\alpha>0$ there 
exists $\beta=\beta(\mu,\alpha)$ such that $F_{\mu}$ has 
a~right inverse $F_{\mu}^{(-1)}$ defined on $\Gamma_{\alpha,\beta}$.
The~function 
$$
\varphi_{\mu}(z)=F_{\mu}^{(-1)}(z)-z
$$
is called the~Voiculescu transform of $\mu$. It is not hard to show that
$\IM \varphi_{\mu}(z)\le 0$ for $z\in \Gamma_{\alpha,\beta}$ where 
$\varphi_{\mu}$ is defined. We also have $\varphi_{\mu}(z)=o(z)$
as $|z|\to\infty$, $z\in\Gamma_{\alpha}$.

If $\mu$ is the~point measure $\delta_a$ at $a$, 
then $F(z)=z-a$ whereas $F(z)=z+ib,\,b\in\mathbb R$, corresponds to 
the~Cauchy distribution with density $x\mapsto b/(\pi(x^2+b^2))$ 
which has infinite variance.
 
\section*{Additive free convolution}

Let $\mu_1$ and $\mu_2$ be p-measures in $\mathcal M$ and let
$F_{\mu_1}(z)$ and $F_{\mu_2}(z)$ denote their reciprocal Cauchy 
transforms respectively. 
We shall define the~free convolution $\mu_1\boxplus\mu_2$,
based on $F_{\mu_1}(z)$ and $F_{\mu_2}(z)$ using the~following result.

\begin{theorem}
There exist unique functions $Z_1(z)$ and $Z_2(z)$ in the~class $\mathcal F$
such that, for $z\in\mathbb C^+$, 
\begin{equation}\label{2.3}
z=Z_1(z)+Z_2(z)-F_{\mu_1}(Z_1(z))\quad\text{and}\quad
F_{\mu_1}(Z_1(z))=F_{\mu_2}(Z_2(z)). 
\end{equation}
\end{theorem}

The~function $F_{\mu_1}(Z_1(z))$ is in $\mathcal F$ again, hence 
there exists some p-measure $\mu$ such that $F_{\mu_1}(Z_1(z))
=F_{\mu}(z)$, where $F_{\mu}(z)=1/G_{\mu}(z)$ and 
$G_{\mu}(z)$ denotes the~Cauchy transform (\ref{2.1}) of $\mu$. 

Since the~p-measure $\mu$ depends on $\mu_1$ and $\mu_2$ only, 
we define $\mu_1\boxplus\mu_2:=\mu$. 

Thus, we defined the~free additive convolution by purely complex analytic
me\-thods (see Chistyakov and G\"otze (2005) \cite{ChG:2005}, 
\cite{ChG:2005a} for a~previous version of this paper). 
A~related approach has been suggested later by Belinschi (2006)
in~\cite{Bel:2006a} and Belinschi and Bercovici (2007) in~\cite{BelBe:2007a} .

The~symmetry of relation (\ref{2.3})
obviously implies that this operation is commutative.
Furthermore, choosing $\mu_2=\delta_a$ in (\ref{2.3}),
where $\delta_a$ denotes a~Dirac measure concentrated at the~point $a$,
we get $\mu_1\boxplus\delta_a=\mu_1*\delta_a$. Definition (\ref{2.3})
does not restrict the~class of p-measures and allows an~obvious extension
to the~case of multiplicative convolutions, described below. 

Moreover,
on any set $\Gamma_{\alpha,\beta}$, where the~functions $\varphi_{\mu_1}(z)$,
$\varphi_{\mu_2}(z)$ and $\varphi_{\mu_1\boxplus\mu_2}(z)$ are defined, 
we obtain immediately from (\ref{2.3}) that
\begin{equation}\label{2.4}
\varphi_{\mu_1\boxplus\mu_2}(z)=\varphi_{\mu_1}(z)+\varphi_{\mu_2}(z).
\end{equation}

The~relation~(\ref{2.4}) implies at once that the~operation $\boxplus$
is associative. 

The~equation (\ref{2.4}) for the~distribution  $\mu_1\boxplus\mu_2$ of $X+Y$, 
where $X$ and $Y$ are free random variables is due to 
Voiculescu~\cite{Vo:1986}.
He considered p-measures $\mu$ with compact support. The~result was 
extended by Maassen~\cite{Maa:1992} to p-measures with finite 
variance; the~general case was proved by Bercovici and 
Voiculescu~\cite{BercVo:1993}.
Note here that Voiculescu's and Bercovici's approach to the~definition
of $\mu_1\boxplus\mu_2$ based on the~operator algebras.
Maassen's analytic approach to the~definition is closer to the~one
presented here.

We see from (\ref{2.4}) that our definition of $\mu_1\boxplus\mu_2$
coincides with the~Voiculescu, Bercovici, Maassen definition.

Since one can investigate free convolutions of several {\it different}
p-measures $\mu_1$, $\dots$, $\mu_n$ using the~multivariate 
approach (\ref{2.3}) with the~help of Nevanlinna functions
which are defined on the~whole half-plane 
$\Bbb C^+$ (see Corollary~2.2 below), 
limit laws can be sucessfully described
as in classical probability theory
(see~\cite{ChG:2006}, \cite{ChG:2006a}).

Voiculescu~\cite{Vo:1993} showed for compactly supported p-measures
that there exist unique functions $F_1,F_2\in\Cal F$ such that 
$G_{\mu_1\boxplus
\mu_2}(z)=G_{\mu_1}(F_1(z))=G_{\mu_2}(F_2(z))$ for all $z\in\Bbb C^+$.
Using Speicher's combinatorial approach~\cite{Spe:1998} (1998) to freeness, 
Biane~\cite{Bia:1998} (1998) proved this result in the~general case.
It follows from Theorem~2.1 that $F_1(z)$ and $F_2(z)$ are $Z_1(z)$ 
and $Z_2(z)$ in (\ref{2.3}), respectively. About the~existence and 
uniqueness of the~subordinating functions $Z_1$ and $Z_2$ see also 
Voiculescu \cite{Vo:2000}, \cite{Vo:2002}.

Pastur and Vasilchuk~\cite{PasVas:2000} (2000) 
studied the~normalized eigenvalue counting measure of the~sum of two
$n\times n$ unitary matrices rotated independently by random unitary 
Haar distributed measures. They established 
the~convergence in probability as $n\to\infty$ to a~limiting
nonrandom measure. They derive functional equations for the~Cauchy 
transforms of limiting distributions 
assuming the~existence of the~mean of the~limiting measures. 
It follows from Theorem~2.1 that
their equations are equivalent to (\ref{2.3}).

\begin{corollary}
Let $\mu_1,\dots,\mu_n\in\mathcal M$.
There exist unique functions $Z_1(z),\dots,Z_n(z)$ in the~class $\mathcal F$
such that, for $z\in\mathbb C^+$,
$$
z=Z_1(z)+\dots+Z_n(z)-(n-1)F_{\mu_1}(Z_1(z)),\quad\text{and}
\quad F_{\mu_1}(Z_1(z))=\dots=F_{\mu_n}(Z_n(z)).
$$
Moreover, $F_{\mu_1\boxplus\dots\boxplus\mu_n}(z)=F_{\mu_1}(Z_1(z))$
for all $z\in\mathbb C^+$.
\end{corollary}

Specializing to $\mu_1=\mu_2=\dots=\mu_n=\mu$ write 
${\mu_1\boxplus\dots\boxplus\mu_n}=\mu^{n\boxplus}$. Then we get
\begin{corollary}
Let $\mu\in\mathcal M$. There exists an~unique function $Z\in\mathcal F$ 
such that
\begin{equation}\label{2.5}
z=nZ(z)-(n-1)F_{\mu}(Z(z)),\quad z\in\mathbb C^+,
\end{equation}
and $F_{\mu^{n\boxplus}}(z)=F_{\mu}(Z(z)),\,z\in\mathbb C^+$.
\end{corollary}

By (\ref{2.5}), we see that $(Z^{(-1)}(z)-z)/(n-1)=z-F_{\mu}(z)$ for $z$
from some domain $\Gamma_{\alpha,\beta}$. It follows from this that
$Z^{(-1)}(z)-z$ has an~analytic continuation to $\mathbb C^+$ with
values in $\mathbb C^-\cup\mathbb R$. Since $Z\in\mathcal F$, it is easy to see,
that $(Z^{(-1)}(iy)-iy)/y\to 0$ as $y\to+\infty$. Hence $Z(z)=F_{\nu}(z),\,
z\in\mathbb C^+$, where $\nu\in\mathcal M$ is infinitely divisible relative
to the~free additive convolution (the~definition and the~characterization
of the~$\boxplus$-infinitely divisibility see in this section below). 
Note that relation (\ref{2.5}) holds if the~integers $n$ are replaced 
by real numbers $t\ge 1$. This shows that there is a~semigroup
$\nu_t\in\mathcal M,\,t\ge 1$, such that $t\varphi_{\nu}=\varphi_{\nu_t}$. 
Thus we give an~analytical approach to the~existence of the~semigroup
$\nu_t\in\mathcal M,\,t\ge 1$, with the~described property.

This fact was shown in  \cite{BercVo:1995}, \cite{NiSp:1996},
\cite{Bel:2006} by other methods.

Having defined the~Voiculescu semigroup $(\mathcal M,\boxplus)$, based
on the~properties of the~functions of the~subclass $\Cal F$  
of the~class $\mathcal N$ of Nevanlinna functions, we shall proceed
by studying multiplicative free convolutions.

\section*{Multiplicative free convolution on $\mathbb R_+$}

Let $\mathcal M_+$ be the~set of p-measures  $\mu$ on 
$\mathbb R_+=[0,+\infty)$ such that $\mu(\{0\})<1$. Define, following 
Voiculescu~\cite{Vo:1987}, the~$\psi_{\mu}$-function of 
a~p-measure $\mu\in\mathcal M_+$, by
\begin{equation}\label{2.6}
\psi_{\mu}(z)=\int\limits_{\mathbb R_+}\frac{z\xi}{1-z\xi}\,\mu(d\xi)
\end{equation}
for $z\in\mathbb C\setminus\mathbb R_+$. The~measure $\mu$ is completely
determined by $\psi_{\mu}$ because $z(\psi_{\mu}(z)+1)=G_{\mu}(1/z)$.
Note that $\psi_{\mu}:\mathbb C\setminus\mathbb R_+\to \mathbb C$ is 
an~analytic function such that $\psi_{\mu}(\bar z)
=\overline{\psi_{\mu}(z)}$, and $z(\psi_{\mu}(z)+1)\in\mathbb C^+$
for $z\in\mathbb C^+$. Consider the~function 
\begin{equation}\label{2.7}
K_{\mu}(z):=\psi_{\mu}(z)/(1+\psi_{\mu}(z)),
\quad z\in\mathbb C\setminus\mathbb R_+.
\end{equation}  
It is easy to see that $K_{\mu}(z)\in\mathcal N$ and $K_{\mu}(z)$
is analytic and nonpositive on the~negative real axis $(-\infty,0)$.
In addition, for $x>0$, $K_{\mu}(-x)\to 0$ as $x\to 0$.

Denote by $\mathcal K$ the~subclass of $\mathcal N$ of functions $f$ such that 
$f(z)$ is analytic and nonpositive on the~negative real axis, and, 
for $x>0$, $f(-x)\to 0$ as $x\to 0$. 

By Krein's results (see Section~3), 
the~function $K_{\mu}(z)$, being analytic and nonpositive on $(-\infty,0)$,
and $\lim_{x\to 0,x>0}K_{\mu}(-x)=0$, belongs to $\mathcal K$ and  
admits by Corollary~3.3 (2) (like {\it all} functions in $\mathcal K$)
the~following representation 
\begin{equation}\label{2.8}
\frac {K_{\mu}(z)}z=a+\int\limits_{(0,\infty)}\frac{\tau(dt)}{t-z},\quad
0<\arg z<2\pi,
\end{equation} 
where $a\ge 0$ and $\tau$ is a~nonnegative measure such that
\begin{equation}\label{2.9}
\int_{(0,\infty)}\frac{\tau(dt)}{1+t}<\infty.
\end{equation}

In view of Proposition~3.4, see Section~3, the~function $K_{\mu}$ 
has the~inverse function $\tilde{\chi}_{\mu}$ on the~image 
$K_{\mu}(i\mathbb C^+)$.
We define the~$\Sigma$-transform of $\mu$ as the~function 
$$
\Sigma_{\mu}(z):=\tilde{\chi}_{\mu}(z)/z,\quad z\in K_{\mu}(i\mathbb C^+). 
$$ 
Note that $K_{\mu}(i\mathbb C^+)\supseteq
\Gamma_{\alpha,\beta,\Delta}^+:=\{z\in\mathbb C:\beta<|z|<\Delta,
\alpha<\arg z<2\pi-\alpha\}$ for some $0<\beta<\Delta$ and $\alpha\in(0,\pi)$.
In addition we conclude
from (\ref{2.8}) that $\arg K_{\mu}(z)\ge\arg z$ for $z\in\Bbb C^+$ and
$\arg K_{\mu}(z)=\pi$ for $z\in(-\infty,0)$.
Therefore $\arg\Sigma_{\mu}(z)\le 0,\,\Im z\ge 0$, and 
$\arg\Sigma_{\mu}(z)=0,\, z\in(-\infty,0)$, where 
$\Sigma_{\mu}(z)$ is defined.  
 
Let $\mu_1$ and $\mu_2$ denote $p$-measures in  $\mathcal M_+$ with
corresponding transforms $K_{\mu_1}$
and $K_{\mu_2}$ defined in (\ref{2.7}), which are in the~klass $\mathcal K$.

We shall define the~free multiplicative convolution using the~transforms 
$K_{\mu_1}$ and $K_{\mu_2}$ by means of the~following characterization
which (after exchanging addition with multiplication) is identical
to characterization (\ref{2.3}) for the~additive convolution. 

\begin{theorem}
There exist two uniquely determined functions $Z_1(z)$ and $Z_2(z)$ 
in the~Krein class $\mathcal K$ such that 
\begin{equation}\label{2.10}
Z_1(z)Z_2(z)=zK_{\mu_1}(Z_1(z))\quad\text{and}\quad 
K_{\mu_1}(Z_1(z))=K_{\mu_2}(Z_2(z)),
\quad z\in\mathbb C^+.
\end{equation}

\end{theorem}

Introduce
$$
K(z):=K_{\mu_1}(Z_1(z))\quad\text{and}\quad\psi(z):=K(z)/(1-K(z)),
\quad z\in\mathbb C^+.
$$
Then, by (\ref{2.8}) for $K_{\mu_1}$ and $Z_1$, we note that 
$K(z)$ and $K(z)/z$ are functions in the~class $\mathcal N$
and $K(-x)\to 0$ as $x\to 0$ for $x>0$. Hence, by Corollary~3.3 (1), 
$K(z)\in\mathcal K$. Using this assertion we easily see 
that $\psi(z)\in\mathcal N$, $\psi(z)/z\in\mathcal N$, and $\lim_{x\to 0,x>0}
\psi(-x)=0$. Moreover, by representation (\ref{2.8}) for $K(z)$, we have 
$\lim_{x\to-\infty}\psi(x)/x=0$. 
Hence the~function $\psi(z)/z$ admits the~representation (\ref{2.8})
with $a=0$, i.e.,
\begin{equation}\label{2.11}
\frac{\psi(z)}z=\int\limits_{(0,\infty)}\frac{\tau_{\psi}(dt)}{t-z}
=\int\limits_{(0,\infty)}\frac{u}{1-uz}\mu_{\psi}(du),\quad z\in\mathbb C^+,
\end{equation}
where $\tau_{\psi}$ is a~nonnegative measure satisfying condition~(\ref{2.9})
and $\mu_{\psi}$ is a~nonnegative finite measure.

Note that $\lim_{x\to-\infty}\psi(x)=-1$ if and only if
in representation (\ref{2.8})
for $K(z)$ either $a>0$ or $\tau((0,\infty))=\infty$. In this case,
by (\ref{2.11}), we may represent $\psi(z)\in\mathcal K$ as $\psi=\psi_{\mu}$,
see (\ref{2.6}), with some p-measure $\mu\in\Cal M_+$ such that
$\mu(\{0\})=0$. In addition $K\in\mathcal K$ may be represented as
$K(z)=K_{\mu}(z),\, z\in\mathbb C^+$. Therefore
$\psi_{\mu_1}(Z_1(z))=\psi_{\mu}(z)$. 

Let in (\ref{2.8}) for $K(z)$ $a=0$ and $\tau((0,\infty))<\infty$. 
Then, as it is easy to see, 
$\lim_{x\to-\infty}\psi(x)=-p=-\tau((0,\infty))/(1+\tau((0,\infty)))$. 
Again, by (\ref{2.11}), we get for $\psi(z)$ representation (\ref{2.6}) with 
some p-measure $\mu\in\mathcal M_+$ and $\mu(\{0\})=1-p$. 
Thus, $\psi(z)=\psi_{\mu}(z)$ and $K(z)=K_{\mu}(z),\,z\in\mathbb C^+$.
Hence $\psi_{\mu_1}(Z_1(z))=\psi_{\mu}(z)$.

The~p-measure $\mu$ is determined 
uniquely by the~p-measures $\mu_1$ and $\mu_2$. 

We define $\mu:=\mu_1\boxtimes\mu_2$. 

This defines the~free multiplicative convolution on the~non-negative
half-line by purely complex analytic me\-thods as above  on p. 4. 

Since
$K_{\mu_1}(Z_1(z))=K_{\mu_2}(Z_2(z))$ for $z\in\mathbb C^+$, 
we have $\mu_1\boxtimes\mu_2=\mu_2\boxtimes\mu_1$ and it is easily 
verified that this convolution is associative as well. 

From Theorem~2.4 we
conclude that the~relation 
\begin{equation}\label{2.11a}
\Sigma_{\mu_1}(z)\Sigma_{\mu_2}(z)=\Sigma_{\mu}(z)
\end{equation} 
holds for $z$, where $\Sigma_{\mu_1}(z),\,\Sigma_{\mu_2}(z)$
and $\Sigma_{\mu}(z)$ are defined.  
This relation is due to Voiculescu~\cite{Vo:1987} and Bercovici
and Voiculescu~\cite{BercVo:1993}. 
 
Using Speicher's combinatorial approach~\cite{Spe:1998} (1998) to freeness, 
Biane~\cite{Bia:1998} (1998) showed that there exist unique functions
$R_1,R_2\in\mathcal K$ such that $\psi_{\mu_1\boxtimes\mu_2}(z)=
\psi_{\mu_1}(R_1(z))=\psi_{\mu_2}(R_2(z))$ for all $z\in\mathbb C^+$.
It follows from Theorem~2.4 that $R_1(z)$ and $R_2(z)$ are $Z_1(z)$ 
and $Z_2(z)$ in (\ref{2.10}), respectively.

\begin{corollary}
Let $\mu_1,\dots,\mu_n\in\mathcal M_+$.
There exist unique functions $Z_1(z),\dots,Z_n(z)$ in the~class 
$\mathcal K$ such that, for $z\in\mathbb C^+$,
$$
Z_1(z)\dots Z_n(z)=z\big(K_{\mu_1}(Z_1(z))\big)^{n-1},\quad\text{and}
\quad K_{\mu_1}(Z_1(z))=\dots=K_{\mu_n}(Z_n(z)).
$$
Moreover, $K_{\mu_1\boxtimes\dots\boxplus\mu_n}(z)=K_{\mu_1}(Z_1(z))$
for all $z\in\mathbb C^+$.
\end{corollary} 

Let $\mu_1=\mu_1=\dots=\mu_n=\mu$. Denote $\mu_1\boxtimes\dots\boxtimes
\mu_n=\mu^{n\boxtimes}$.
\begin{corollary}
Let $\mu\in\mathcal M_+$. There exists an~unique function $Z\in\mathcal K$ 
such that
\begin{equation}\label{2.12}
(Z(z))^n=z\big(K_{\mu}(Z(z))\big)^{n-1},\quad z\in\mathbb C^+,
\end{equation}
and $K_{\mu^{n\boxtimes}}(z)=K_{\mu}(Z(z)),\,z\in\mathbb C^+$.
\end{corollary}

Again relation (\ref{2.12}) holds if we replace integers
$n$ by real numbers $t\ge 1$. This shows that there is a~semigroup
$\nu_t\in\mathcal M_+,\,t\ge 1$, such that $(\Sigma_{\nu}(z))^t
=\Sigma_{\nu_t}(z)$.

Thus we have defined the~Voiculescu semigroup $(\mathcal M_+,\boxtimes)$,
based on properties of functions of the~Krein subclass $\mathcal K$  
of the~class $\mathcal N$ of Nevanlinna functions. In the~following
we consider the~case of spectral p-measures on the~unit circle $\mathbb T$.

\section*{Multiplicative free convolution on the~unit circle}

Let $\mu$ be a~p-measure on $\mathbb T$. Following 
Voiculescu~\cite{Vo:1987}, we define a~transform of the~p-measure $\mu$ 
on $\mathbb T$, by
$$
\psi_{\mu}(z)=\int\limits_{\mathbb T}\frac{z\xi}{1-z\xi}\,\mu(d\xi).
$$
The~function $\psi$ has a~convergent power series representation 
in $\mathbb D:=\{z\in\mathbb C:\,|z|<1\}$,
the~open unit disk of $\Bbb C$, such that $\psi_{\mu}(0)=0$.

Let $\mathcal M_*$ denote the~set of p-measures on $\mathbb T$
such that $\int_{\mathbb T}\xi\,\mu(d\xi)\ne 0$. 

If $\mu\in\mathcal M_*$, it follows that
the~function 
$$
Q_{\mu}:=\psi_{\mu}/(1+\psi_{\mu})
$$ 
has a~right inverse
$Q_{\mu}^{(-1)}$, defined in a~neighborhood of $0$ denoted by
$\mathbb D_{\alpha}:=\{z\in\mathbb C:\,|z|<\alpha\}$ with some $0<\alpha\le 1$,
such that $Q_{\mu}^{(-1)}(0)=0$. Let 
$$
\Sigma_{\mu}(z)=Q_{\mu}^{(-1)}(z)/z
$$ 
denote the~so called $\Sigma$-transform of $\mu$.

Denote by $\mathcal C$ the~class of analytic functions $H(z)$ on
$\mathbb D\to-i(\mathbb C^+\cup\mathbb R)$ introduced by Carath\'eodory.

Note that
\begin{equation}\label{2.13}
Q_{\mu}(z)=\frac{\psi_{\mu}(z)}{1+\psi_{\mu}(z)}=\frac{H(z)-1}{H(z)+1}
\end{equation}
where $H(z):=1+2\psi_{\mu}(z)$ is a~function of Carath\'eodory's 
class $\mathcal C$. Such functions $H(z)$, by (\ref{3.1}) and $H(0)=1$, 
see Section~3, have the~form
$$
H(z)=\int\limits_{\mathbb T}\frac{\xi+z}{\xi-z}
\,\sigma(d\xi),
$$
where $\sigma$ is a~p-measure.

Define $\mathcal S$ to be so called {\it Schur class} of analytic functions
$\mathbb D\to\overline{\mathbb D}$ (see Section~3), where $\overline{\mathbb D}$
is the~closure of $\mathbb D$.

We see from (2.13), that $Q_{\mu}\in \mathcal S$ and since $\psi_{\mu}(0)=0$ 
and $\mu\in\mathcal M_*$, $Q_{\mu}(0)=0,\,
Q_{\mu}'(0)\ne 0$. 

In the~sequel we denote
by $\mathcal S_*$
the~subclass of $\mathcal S$ which consists of Schur functions
$Q$ with properties $Q(0)=0$ and $Q'(0)\ne 0$.  

Since $Q_{\mu} \in\mathcal S_*$, by the~Schwarz lemma, we have 
$|Q_{\mu}(z)/z|\le 1$ for $z\in\mathbb D$. Hence $|Q_{\mu}^{(-1)}(z)/z|\ge 1$ 
in a~neighborhood of $0$. Moreover, both (\ref{3.2}) induces a~one-to-one
correspondence between the~classes $\mathcal C$ and $\mathcal S$, and 
(\ref{2.13}) induces
a~one-to-one correspondence between functions $H\in\mathcal C$ such that
$H(0)=1$ and $H'(0)\ne 0$, and
functions $Q_{\mu}$ of the~class $\mathcal S_*$.

Let $\mu_1$ and $\mu_2$ denote p-measures in $\mathcal M_*$ and let $Q_{\mu_1}$ 
and $Q_{\mu_2}$ be Schur functions which correspond to these 
measures, by (\ref{2.13}). We now define the free multiplicative free 
convolution $\mu_1\boxtimes\mu_2$ based on $Q_{\mu_1}$ and
$Q_{\mu_2}$ using the~following characterization.
   
\begin{theorem}
There exist two functions $Z_1(z)$ and $Z_2(z)$ in the~set $\mathcal S_*$ 
such that 
\begin{equation}\label{2.14}
Z_1(z)Z_2(z)=zQ_{\mu_1}(Z_1(z))\quad\text{and}\quad 
Q_{\mu_1}(Z_1(z))=Q_{\mu_2}(Z_2(z)),\quad z\in\mathbb D.
\end{equation}
The~functions $Z_1(z)$ and $Z_2(z)$ are unique solutions of $(\ref{2.14})$
in the~class $\mathcal S_*$.
\end{theorem}

Consider the~function $Q_{\mu_1}(Z_1(z))$. It is easy
to see that it belongs to Schur's class $\mathcal S$ and 
$Q_{\mu_1}(Z_1(0))=0,\,Q_{\mu_1}'(0)Z_1'(0)\ne 0$. Therefore 
$Q_{\mu_1}(Z_1)\in\Cal S_*$ and $Q_{\mu_1}(Z_1(z))=Q_{\mu}(z)$  
for $z\in\Bbb D$, where $Q_{\mu}(z)$ has form (\ref{2.13}) 
for some p-measure $\mu\in\mathcal M_*$.
This measure is determined uniquely by the~p-measures $\mu_1$
and $\mu_2$. Define $\mu:=\mu_1\boxtimes\mu_2$. 

This defines the~free multiplicative convolution on the~unit
circle by purely complex analytic me\-thods as on p. 4. 

Since
$Q_{\mu_1}(Z_1(z))=Q_{\mu_2}(Z_2(z))$ for $z\in\mathbb D$, we get
$\mu_1\boxtimes\mu_2=\mu_2\boxtimes\mu_1$. It is easy to verify that 
the~ope\-ration $\boxtimes$ is associative.
Using relation (\ref{2.13}) between
the~function $Q_{\mu}(z)\in \mathcal S_*$ and the~function $\psi_{\mu}(z)$, we
conclude that $\psi_{\mu}(z)=\psi_{\mu_1}(Z_1(z))$ for $z\in\mathbb D$.
In addition we have in some neighborhood of $0$
\begin{equation}\label{2.15}
\frac{Q_{\mu_1}^{(-1)}(z)}z\frac{Q_{\mu_2}^{(-1)}(z)}z
=\frac{Q_{\mu}^{(-1)}(z)}z\quad\text{or}\quad
\Sigma_{\mu_1}(z)\Sigma_{\mu_2}(z)=\Sigma_{\mu}(z).
\end{equation} 
This formula is due to Voiculescu~\cite{Vo:1987}, \cite{BercVo:1992}
(1992).

Biane~\cite{Bia:1998} (1998) showed that there exist unique functions
$Q_1,Q_2\in\Cal S_*$ such that $\psi_{\mu_1\boxtimes\mu_2}(z)=
\psi_{\mu_1}(Q_1(z))=\psi_{\mu_2}(Q_2(z))$ for all $z\in\mathbb D$.
It follows from Theorem~2.7 that $Q_1(z)$ and $Q_2(z)$ are $Z_1(z)$ 
and $Z_2(z)$ in (\ref{2.14}), respectively.

Note that Vasilchuk~\cite{Vas:2001} (2001) studied the~normalized
eigenvalue counting measure of the~product of two $n\times n$ unitary
matrices and the~measure of product of three $n\times n$ Hermitian positive
matrices rotated independently by random unitary Haar distributed measures.
He established 
the~convergence in probability as $n\to\infty$ to a~limiting
nonrandom measure and derived functional equations for 
the~Herglotz and Cauchy transforms of limiting distributions under 
some restriction on counting measures. From Theorem~2.4 and Theorem~2.7 
it follows that his equations are equivalent to (\ref{2.10}) and (\ref{2.14}).

For the~multiplicative free convolution the~analogues of Corrolary~2.2
and Corollary~2.3 hold.
\begin{corollary}
Let $\mu_1,\dots,\mu_n\in\mathcal M_*$.
There exist uniquely determined functions $Z_1(z),\dots,Z_n(z)$ 
in the~class $\mathcal S_*$ such that, for $z\in\mathbb D$,
$$
Z_1(z)\dots Z_n(z)=z\big(Q_{\mu_1}(Z_1(z))\big)^{n-1},\quad\text{and}
\quad Q_{\mu_1}(Z_1(z))=\dots=Q_{\mu_n}(Z_n(z)).
$$
Moreover, $Q_{\mu_1\boxtimes\dots\boxplus\mu_n}(z)=Q_{\mu_1}(Z_1(z))$
for all $z\in\mathbb D$.
\end{corollary}

\begin{corollary}
Let $\mu\in\mathcal M_*$. There exists an~unique function $Z\in\mathcal S_*$ 
such that
\begin{equation}\label{2.16}
(Z(z))^n=z\big(Q_{\mu}(Z(z))\big)^{n-1},\quad z\in\mathbb D,
\end{equation}
and $Q_{\mu^{n\boxtimes}}(z)=Q_{\mu}(Z(z)),\,z\in\mathbb D$.
\end{corollary}

Rewrite (\ref{2.16}) in the~form
\begin{equation}\label{2.16a}
\tilde{Z}(z)=\big(\tilde{Q}_{\mu}(Z(z))\big)^{n-1},\quad z\in\mathbb D,
\end{equation}
where $\tilde{Z}(z):=Z(z)/z$ and $\tilde{Q}_{\mu}(z)
:=Q_{\mu}(z)/z,\,z\in\Bbb D$, are functions
of the~class $\mathcal S$ and $\tilde{Z}(0)\ne 0$ and $\tilde{Q}_{\mu}(0)\ne 0$. 
From (\ref{2.16a})
it follows that $\tilde{Z}(z)\ne 0$ and $\tilde{Q}_{\mu}(Z(z))\ne 0$ for 
$z\in\mathbb D$.

Note that relation (\ref{2.16a}) holds if we replace integers 
$n$ by real $t\ge 1$ for the~functions $\tilde{Q}_{\mu\boxtimes\mu}(z)$ and
for $\tilde{Q}_{\mu}(z)$ if $Q_{\mu}(z)\ne 0$ for $z\in\mathbb D\setminus \{0\}$.  
This shows that there is a~semigroup $\nu_t\in\mathcal M_*,\,t\ge 2$, 
such that $(\Sigma_{\nu}(z))^t=\Sigma_{\nu_t}(z)$ in the~general case
and there is a~semigroup $\nu_t\in\mathcal M_*,\,t\ge 1$,
such that $(\Sigma_{\nu}(z))^t=\Sigma_{\nu_t}(z)$ in the~case
$Q_{\nu}(z)\ne 0$ for $z\in\mathbb D\setminus \{0\}$. 

Thus we have defined the~Voiculescu semigroup $(\mathcal M_*,\boxtimes)$,
based on properties of functions of the~subclass $\mathcal S_*$  
of the~class $\Cal S$ of Schur functions.

The~relations~(2.3), (2.10), and (2.14) were used to good effect 
in the~papers of Belinshi~\cite{Bel:2003}, \cite{Bel:2006}, Belinshi and 
Bercovici~\cite{BelBe:2004}, \cite{BelBe:2005}, Bercovici and 
Voiculescu~\cite{BeVo:1998}, Biane~\cite{Bia:1997}. In our paper
we use the~relations~(2.3), (2.10), and (2.14) to study 
the~arithmetic of p-measures in Voiculescu's semigroups.

\section*{Arithmetic of p-measures in Voiculescu's semigroups}

Now we consider the~problem of the~decomposition of measures $\mu$
of the~commutative semigroups $(\mathcal M,\boxplus)$, $(\mathcal M_+,\boxtimes)$, 
and $(\mathcal M_*,\boxtimes)$. In the~sequel we shall denote these semigroups
by a~symbol $(\mathbf M,\circ)$, where $\mathbf M$ means $\mathcal M,\mathcal M_+,
\mathcal M_*$ and $\circ$ means the~operations $\boxplus,\boxtimes$. 
The~following notions are analogues of the~classical ones for $(\mathcal M,*)$.

We shall say that $\mu_1\in\mathbf M$ is a~{\it free factor} or just
{\it factor} of $\mu\in\mathbf M$ if there exists 
$\mu_2\in\mathbf M$ such that $\mu=\mu_1\circ\mu_2$.
Every $\mu\in\mathbf M$ has factors. Indeed, we have $\mu=\delta_a
\circ(\mu\circ\delta_b)$, where $b=-a, a\in\mathbb R$, in the~case
of the~semigroup $(\mathcal M,\boxplus)$, $b=1/a, a>0$, in the~case
of $(\mathcal M_+,\boxtimes)$, and $b=1/a, a\in\Bbb T$, in the~case
of $(\mathcal M_*,\boxtimes)$.
Hence $\delta_a$ and $\mu\circ\delta_b$ are factors of $\mu$. 
Such factors are called {\it improper}. A~p-measure $\mu$ 
which is not a~Dirac measure
is called {\it indecomposable} if it has improper factors only. 
Such p-measures may be regarded as simple elements
of this semigroup. If $\mu$ in $(\mathbf M,\circ)$ is not 
indecomposable it is called {\it decomposable}.
Two measures $\mu_1$ and $\mu_2$ are called {\it equivalent}, 
$\mu_1\sim\mu_2$, if $\mu_1=\mu_2\circ\delta_a$, where $a\in\mathbb R$ 
in the~case of $(\mathcal M,\boxplus)$, $a>0$ in the~case of
$(\mathcal M_+,\boxtimes)$, and $a\in\mathbb T$ in the~case of  
$(\mathcal M_*,\boxtimes)$.

As in the~classical theory of convolutions, a~measure $\mu$ is called 
$\circ$-{\it infinitely divisible} (or i.d. for short) if, 
for every natural number $n$, $\mu$ can be written as 
$\mu=\nu_n\circ\nu_n\circ\dots\circ\nu_n$ ($n$ times) with $\nu_n\in\mathbf M$. 
The~measure $\delta_a$ is necessarily infinitely divisible. 
Note that all i.d. measures are decomposable and Dirac measures
have this property.
As mentioned in the~introduction a~measure $\mu\in\mathbf M$ is to
the~class $I_0$ relative $\circ$ if $\mu$ is $\circ$-infinitely
divisible and has no $\circ$-indecomposable factors. 

Khintchine~\cite{Khi:1937} (1937) was the~first
who studied the~arithmetic of the~semigroup $(\mathcal M,*)$ of 
distribution functions on $\mathbb R$ with respect to the~convolution $*$.
He derived for this semigroup the~three basic results.

1. {\bf Limit of triangular arrays.} Obviously an~i.d. element $\mu$ 
of $(\mathcal M,*)$ can always be 
represented as the~limit of a~convergent infinitesimal triangular array.

2. {\bf Classification.} Any element $\mu$ of $(\mathcal M,*)$ belongs 
to one of the~following classes. Either
\begin{enumerate}
\item $\mu$ is indecomposable, 

\item $\mu$ is decomposable (possibly i.d.) and has
an~indecomposable factor, 

\item $\mu$ is i.d. and has no indecomposable factors.
(This class of p-measures is denoted by $I_0$.)
\end{enumerate}
(Compare the~description of Delphic semigroups in the~introduction.)

3. {\bf Representation.} For each $\mu$ in $\mathcal M$ one may decompose 
$$
\mu=\nu*\mu_1*\mu_2*\dots
$$
in at least one way, where $\nu$ is i.d. and
has no indecomposable sub-factor, and each $\mu_j$ is 
indecomposable. The~convolution product is at most countable and may
be finite or void. 

Note that Gaussian distributions (Cram\'er~(1936)), Poisson 
distributions (Rai\-kov (1937)) and the~convolution of Gaussian and
Poisson distributions (Linnik~(1957)) belong to the~class $I_0$.
Hence in the~semigroup $(\mathcal M,*)$ the~class $I_0$ has  nontrivial
elements (besides the~trivial units $\delta_a$). 

A number of papers have been devoted to the~study of the~arithmetic of
semigroups of p-measures. We refer the~reader
to the~monograph of Linnik and Ostrovskii~\cite{LinOst:1977} (1977), 
the~surveys of Livshic, Ostrovskii and Chistyakov~\cite{LivOst:1975} (1975), 
Ostrovskii~\cite{Ost:1977} (1977), \cite{Ost:1986} (1986).

We shall consider the~semigroups $(\mathbf M,\circ)$ introduced above 
and we study their arithmetic using the~theory of Delphic
semigroups, (see Kendall~\cite{Ken:1968}, Davidson~\cite{Dav:1968(I)}
--\cite{Dav:1969}). 
Kendall and Davidson developed Khintchine's theory for a~wide class
of semigroups (Delphic) where Khintchine's basic theorems
remain valid.

In the~semigroup $(\mathcal M,\boxplus)$ i.d. p-measures were 
first considered in Voiculescu~\cite{Vo:1986}, where compactly supported
$\boxplus$-i.d. measures were characterized.
P-measures with a~finite variance were considered in Maassen~\cite{Maa:1992}
and Bercovichi, Voiculescu~\cite{BercVo:1993} gave a~characterization of
general i.d. p-measures $\mu\in\mathcal M$. There is an~analogue of 
the~L\'evy-Khintchine formula, (see \cite{VoDyNi:1992}, 
\cite{BercVo:1992}, \cite{BercVo:1993})
which states that a~p-measure $\mu$, on $\mathbb R$, is i.d. 
if and only if the~function $\varphi_{\mu}(z)$ has an~analytic
continuation to the~whole of $\mathbb C^+$, with values in $\mathbb C^-
\cup\mathbb R$, and one has
\begin{equation}\label{2.17}
\lim_{y\to +\infty}\frac{\varphi_{\mu}(iy)}y=0.
\end{equation}
By the~Nevanlinna representation for such function, we know that
there exist a~real number $\alpha$, and a~finite nonnegative
measure $\nu$, on $\mathbb R$, such that
\begin{equation}\label{2.18}
\varphi_{\mu}(z)=\alpha+\int\limits_{\mathbb R}\frac{1+uz}{z-u}\,\nu(du).
\end{equation}

\noindent
There is a~one-to-one correspondence between functions
$\varphi_{\mu}(z)$ and couples $(\alpha,\nu)$. For this reason
we shall sometimes write $\varphi_{\mu}=(\alpha,\nu)$.  

Formula~(\ref{2.18}) is an~analogue of the~well-known L\'evy-Khintchine
formula for the~logarithm of characteristic functions 
$\varphi(t;\mu):=\int_{\mathbb R}e^{itu}\,\mu(du),\,t\in\mathbb R$, of 
$*$-i.d. measures $\mu\in\mathcal M$. A~measure $\mu\in\mathcal M$ 
is $*$-i.d. if and only if there exist a~finite 
nonnegative Borel measure $\nu$ on $\mathbb R$, and a~real number
$\alpha$ such that
\begin{equation}\label{2.19}
\log\varphi(t;\mu)=f_{\mu}(t):=\exp\Big\{i\alpha t
+\int\limits_{\mathbb R}\Big(e^{itu}-1
-\frac {itu}{1+u^2}\Big)\frac{1+u^2}{u^2}\,\nu(du)\Big\},
\quad t\in\mathbb R,
\end{equation}
where $(e^{itu}-1-itu/(1+u^2))(1+u^2)/u^2$ will be interpreted as
$-t^2/2$ for $u=0$.
There is a~one-to-one correspondence between functions
$f_{\mu}(t)$ and couples $(\alpha,\nu)$. For this reason
we shall sometimes write $f_{\mu}=(\alpha,\nu)$.  

In the~classical case the~precise formulation of the 
Khintchine limit theorem for $(\mathcal M,*)$ is as follows:

Let $\{\mu_{nk}:\,n\ge 1,\,1\le k\le n\}$ be an~array of 
infinitesimal measures in $\mathcal M$, i.e.,
\begin{equation}\label{2.20}
\lim_{n\to\infty}\max_{1\le k\le n}\mu_{nk}(\{u:|u|>\varepsilon\})=0
\end{equation} 
for every $\varepsilon>0$. 
In order that $\mu\in\mathcal M$ be the~limit in the~weak topology 
of distributions $\mu^{(n)}=\delta_{a_n}*\mu_{n1}
*\mu_{n2}*\dots *\mu_{nn}\to\mu$ for some
suitably chosen constants $a_n$, it is necessary and sufficient that
$\mu$ be i.d.

(Without loss of generality, we shall in the~sequel 
consider arrays of the~length $n$ instead of $k_n$).

The~Khintchine limit theorem in free probability theory 
has the~same form for the~measures $\mu^{(n)}=\delta_{a_n}\boxplus\mu_{n1}
\boxplus\mu_{n2}\boxplus\dots\boxplus\mu_{nn}$.

This theorem was early proved by Bercovici and Pata~\cite{BercPa:2000}
(2000). 
We give another proof of this result, using arguments of 
the~theory of Delphic semigroups.

The~i.d. measures in $(\mathcal M_+,\boxtimes)$
have been characterized by Voicu\-le\-scu~\cite{Vo:1987}, Bercovici 
and Voiculescu~\cite{BercVo:1992}, \cite{BercVo:1993}. There is 
an~analogue of 
the~L\'evy-Khintchine formula which states that a~measure 
$\mu\in\mathcal M_+$ is $\boxtimes$-i.d.
if and only if there exist a~finite nonnegative measure
$\nu$ on $(0,\infty)$ and real numbers $a$ and $b\ge 0$
such that

\begin{equation}\label{2.21}
\Sigma_{\mu}(z)=\exp\Big\{a-bz+\int\limits_{\mathbb R_+}
\frac{1+uz}{z-u}\,\nu(du)\Big\},
\end{equation}  
for $z$, where $\Sigma_{\mu}(z)$ is defined.
For this reason we will write $\Sigma_{\mu}=(a,b,\nu)$.

In other words, a~measure $\mu\in\mathcal M_+$ is $\boxtimes$-i.d. 
if and only if
\begin{equation}\label{2.22}
\Sigma_{\mu}(z)=\exp\{-u(z)\}, 
\end{equation} 
where $u(z)\in\mathcal N$ and $u(z)$ is analytic and real valued 
on the~negative half-line $(-\infty,0)$. 

Khintchine's limit problem for multiplicative free convolution
may be formulated as follows.

Let $\{\mu_{nk}:\,n\ge 1,\,1\le k\le n\}$ be an~array of measures
in $\mathcal M_+$ such that
\begin{equation}\label{2.23}
\lim_{n\to\infty}\max_{1\le k\le n}\mu_{nk}(\{u:\,
|u-1|>\varepsilon\})=0
\end{equation}
for every $\varepsilon >0$. The~measures $\mu_{nk}\in \mathcal M_+$ are
called {\it infinitesimal}. 

We shall characterize in Theorem~2.10 the~class of p-measures
$\mu\in\mathcal M_+$ such that 
$\mu^{(n)}=\delta_{a_n}\boxtimes\mu_{n1}\boxtimes\mu_{n2}
\dots\boxtimes\mu_{nn}\to\mu$ in the~weak topology for some suitably chosen
positive constants $a_n$.

The~i.d. measures of the~semigroup
$(\mathcal M_*,\boxtimes)$ were characterized 
in \cite{Vo:1987}, \cite{BercVo:1992}. 
There is an~analogue of 
the~L\'evy-Khintchine formula which states that a~measure 
$\mu\in\mathcal M_*$ is $\boxtimes$-i.d.
if and only if there exist a~finite nonnegative measure
$\nu$ on $\mathbb T$ and a~real number $a$ 
such that
\begin{equation}\label{2.24}
\Sigma_{\mu}(z)=\exp\Big\{ia+\int\limits_{\mathbb T}
\frac{1+z\xi}{1-z\xi}\,\nu(d\xi)\Big\}, 
\end{equation}
for $z$, where $\Sigma_{\mu}(z)$ is defined.
For this reason we will write $\Sigma=(a,\nu)$.

In other words, a~measure $\mu\in\mathcal M_*$ is $\boxtimes$-i.d. 
if and only if
\begin{equation}\label{2.25}
\Sigma_{\mu}(z)=\exp\{v(z)\}, 
\end{equation} 
where $v(z)\in\mathcal C$.  

Let $\{\mu_{nk}:n\ge 1,\,1\le k\le n\}$ be an~array of measures in 
$\mathcal M_*$. We shall call the~measures $\mu_{nk}$ {\it infinitesimal}
if
\begin{equation}\label{2.26}
\lim_{n\to\infty}\max_{1\le k\le n}\mu_{nk}(\{\xi:\,|\arg\xi|>
\varepsilon\})=0.
\end{equation}
The~Khintchine limit problem for multiplicative free 
convolution for measures $\mu\in\mathcal M_*$ has the~same
form as in the~case of $\mu\in\mathcal M_+$ with constants $a_n\in\mathbb T$.

We give the~solution of the~Khintchine limit problem,
proving the~following result.

\begin{theorem}
Let $\{\mu_{nk}:\,n\ge 1,\,1\le k\le n\}$ be infinitesimal 
probability measures in the~semigroup $(\mathbf M,\circ)$.
The~family of limit measures of sequences 
$\mu^{(n)}=\delta_{a_n}\circ\mu_{n1}\circ\mu_{n2}
\circ\dots\circ\mu_{nn}$ for some suitably chosen
constants $a_n$ coincides with the~family of $\circ$-infinitely 
divisible p-measures.
\end{theorem}
We give the~proof of this theorem, using arguments of the theory
of Delphic semigroups. About the~Khintchine limit problem for multiplicative free 
convolution see Belinschi--Bercovici~\cite{BelBe:2007} as well.

The~arithmetic of the~semigroups $(\mathbf M,\circ)$ is described
in the~following results. 
\begin{theorem}
The~element $\mu$ of $(\mathbf M,\circ)$ can be classified as follows.
Either
\begin{enumerate}
\item $\mu$ is indecomposable, 

\item $\mu$ is decomposable $($possibly infinitely divisible$)$ and has
an~indecomposable factor, 

\item $\mu$ is infinitely divisible and has no indecomposable factors.
$($This class will be denoted by $I_0$.$)$
\end{enumerate}
\end{theorem}

P-measures may be decomposed as follows
\begin{theorem}
Every probability measure $\mu$, which has indecomposable factors,
can be expressed in the~form
\begin{equation}\label{2.27}
\mu=\mu_0\circ\mu_1\circ\mu_2\circ\dots,
\end{equation}
where $\mu_0\in I_0$ and the~probability measures $\mu_1,\mu_2,\dots$ 
are indecomposable $($its number may be finite or enumerable$)$. 
\end{theorem}

We will show that representation (\ref{2.27}) is not unique.

The~class $I_0$ in $(\mathbf M,\circ)$, mentioned above may be described 
as follows.

\begin{theorem}
In Voiculescu's semigroup $(\mathbf M,\circ)$ the~class $I_0$
is trivial, that is $I_0$ is the~class of Dirac measures.
\end{theorem}

In Section~8 we describe wide classes of indecomposable elements
in $(\mathcal M,\circ)$ (see Theorems~\ref{8.1th}--\ref{8.3th}) 
from which it follows the~following result.

\begin{theorem}
The~probability measures with support consisting of a~finite number
of points are indecomposable in $(\mathbf M,\circ)$.

\end{theorem}

\begin{corollary}
The~class of indecomposable elements of $(\bold M,\circ)$ is
dense in $(\mathbf M,\circ)$ in the~weak topology.
\end{corollary}

Theorem~2.14 follows from Belinshi's~\cite{Bel:2006} and 
Bercovici--Wang~\cite{BercWa:2007} results. But, as we note in Section~8,
Theorems~8.1--8.3 do not follow from these results.

Theorem~2.11 describes the~class $I_0$ as the~class
of i.d. elements of $(\mathbf M,\circ)$ which
have i.d. components only. Bercovici and Voiculescu~\cite{BercVo:1995}
proved that a~semicircular measure does not belong to the~class $I_0$
in the~semigroup $(\mathcal M,\boxplus)$. Benaych-Georges~\cite{Ben:2005} proved a~similar
result for the free Poisson measure. These results follow from 
Theorem~2.13.

Speicher and Woroudi~\cite{SpeWo:1997} (1997) introduced a~further 
convolution operation
on $\Cal M$ denoted  $\biguplus$. Let $\mu\in\mathcal M$. Denote as before 
by $G_{\mu}$ the~Cauchy transform of $\mu$ and by 
$F_{\mu}=1/G_{\mu}:\mathbb C^+\to\mathbb C^+$
its reciprocal. We have $\IM z\le\IM F_{\mu}(z)$ so that the~function
$E_{\mu}(z):=z-F_{\mu}(z)$ maps $\mathbb C^+$ to $\mathbb C^-\cup\mathbb R$, and,
in addition, $E_{\mu}(z)/z\to 0$ as $z\to\infty,\,z\in\Gamma_{\alpha}$
for any fixed $\alpha>0$.
Conversely, if $E:\,\mathbb C^+\to\mathbb C^-\cup\mathbb R$ is an~analytic function 
so that $E(z)/z\to 0$ as $z\to\infty,\,z\in\Gamma_{\alpha}$ for any fixed
$\alpha$, then there
exists $\mu\in\mathcal M$ such that $E_{\mu}=E$. This observation leads
to the~formal definition of the~Boolean convolution.
Given $\mu,\nu\in\mathcal M$, there exists
$\rho\in\mathcal M$ such that
$$
E_{\rho}=E_{\mu}+E_{\nu}.
$$
The~measure $\rho$ is called the~{\it Boolean convolution} of $\mu$ and $\nu$,
and it denoted $\mu\biguplus\nu$. Boolean convolution is an~associative
and commutative law, with $\delta_0$ as the~zero element. We have $\delta_s
\biguplus\delta_t=\delta_{s+t}$, but generally $\delta_t\biguplus\mu$ is not 
a~translate of $\mu$. Speicher and Woroudi~\cite{SpeWo:1997} treated 
the~central and Poisson limit theorems, characterized i.d. and
stable distributions and proved analogues of the~classical theorems 
of Cram\'er, Marcinkiewicz, Kac and Lo\`eve. Bercovici and 
Pata~\cite{BercPa:1999} (1999)
established limit laws for Boolean convolutions.

Speicher and Woroudi~\cite{SpeWo:1997} noted that all p-measures $\mu$ are 
i.d. It is easy to see that the~class of indecomposable
elements in the~semigroup $(\mathcal M,\biguplus)$ is empty. Moreover,
the~class of i.d. elements coincides with the~class $I_0$.
Therefore Theorem~2.12 obviously holds for $(\mathcal M,\biguplus)$ .

\section{Auxiliary results}

In the~following we state some results about classes of analytic functions
(see Akhiezer~\cite{Akh:1965} (1965), Section~3, Akhiezer and 
Glazman~\cite{AkhGla:1963} (1963), Section~6, \S 59, and Krein and 
Nudelman~\cite{KrNu:1977} (1977), Appendix).

By $\mathcal C$ we denote C. Carath\'eodory's class of analytic
functions $F(z):\Bbb D\to \{z:\,\RE z\ge 0\}$. A~function $F$ is in
$\mathcal C$ if and only if it admits the~following representation 
(Herglotz,~G., Riesz,~F.)
\begin{equation}\label{3.1}
F(z)=ia+\int\limits_{\mathbb T}\frac{\xi+z}{\xi-z}
\,\sigma(d\xi),
\end{equation}
where $a=\IM F(0)$, $\mathbb T$ is the~unit circle,
and $\sigma$ is finite nonnegative measure. The~number $a$ and 
the~measure $\sigma$ are uniquely determined by $F$.

By $\Cal S$ we denote J. Schur's class of analytic functions
$\varphi(z):\mathbb D\to \overline{\mathbb D}$. 
The~classes $\mathcal C$ and $\mathcal S$ are connected via
\begin{equation}\label{3.2}
\varphi(z)=\frac 1z\frac{F(z)-F(0)}{F(z)+\overline{F(0)}},
\end{equation}
which induces a~one-to-one correspondence between $\mathcal C$ and $\mathcal S$.

Finally we denote by $\mathcal N$ R. Nevanlinna's class of analytic 
functions $f(z):\mathbb C^+\to \mathbb C^+\cup\mathbb R$.
A~function $f$ is in $\Cal N$ if and only if it admits an~integral 
representation
\begin{equation}\label{3.3}
f(z)=a+bz+\int\limits_{\mathbb R}\frac{1+uz}{u-z}\,\tau(du)=
a+bz+\int\limits_{\mathbb R}\Big(\frac 1{u-z}-\frac u{1+u^2}\Big)
(1+u^2)\,\tau(du),
\end{equation}
where $b\ge 0$, $a\in\mathbb R$, and $\tau$ is finite nonnegative 
measure. Here $a,b$ and $\tau$ are uniquely determined by $f$.
More precisely we have $a=\RE f(i)$ and $\tau(\mathbb R)=\IM f(i)-b$.   
From this formula it follows that
\begin{equation}\label{3.4}
f(z)=(b+o(1))z
\end{equation}
for $z\in\mathbb C^+$ such that $|\RE z|/\IM z$ stays bounded as $|z|$
tends to infinity. Hence if $b\ne 0$, then $f$ has a~right inverse 
$f^{(-1)}$ defined on the~region $\Gamma_{\alpha,\beta}$ 
defined (in Section~2) for any $\alpha>0$ and some positive
$\beta=\beta(f,\alpha)$.

A function $f\in\Cal N$ admits the~representation
\begin{equation}\label{3.4a}
f(z)=\int\limits_{\mathbb R}\frac{\sigma(du)}{u-z},\quad z\in\mathbb C^+,
\end{equation}
where $\sigma$ is a~finite nonnegative measure, if and only if
$\sup_{y\ge 1}y|f(iy)|<\infty$.

Note that the~class $\mathcal F$ coincides with the~subclass of Nevanlinna
functions for which $f(z)/z\to 1$ as $z\to\infty$ nontangentially.
Indeed, reciprocal Cauchy transforms of p-measures have obviously
such property. Let $f\in\mathcal N$ and $f(z)/z\to 1$ as $z\to\infty$ 
nontangentially. Then, by (\ref{3.4}), $f$ admits the~representation 
(\ref{3.3}), where $b=1$. By (\ref{3.4}) and (\ref{3.4a}), we see that 
$-1/f(z)$ admits the~representation (\ref{3.4a}), where $\sigma\in\Cal M$.
It follows from (\ref{3.3}) with $b=1$ that a~function $f\in\mathcal F$
admits the~inequality
\begin{equation}\label{3.4b}
\Im f(z)\ge\Im z,\qquad z\in\mathbb C^+.
\end{equation}

The~Stieltjes-Perron inversion formula for a~function $f\in\Cal N$ 
has the~following form.

Let $\psi(u):=\int_0^u(1+t^2)\,\tau(dt)$. Then
\begin{equation}\label{3.5}
\psi(u_2)-\psi(u_1)=\lim_{\eta\to 0}\frac 1{\pi}\int\limits_{u_1}
^{u_2}\IM f(\xi+i\eta)\,d\xi,
\end{equation}
where $u_1<u_2$ are continuity points of the~function $\psi(u)$.

The following two results are due to Krein, M.
\begin{theorem}\label{3.1th}
The~function $f(z)$ admits the~representation
\begin{equation}\label{3.6}
f(z)=a+\int\limits_{\mathbb R_+}\frac{\tau(du)}{u-z},\quad
0<\arg z<2\pi,
\end{equation}                                
where $a\ge 0$ and $\tau$ is a~nonnegative measure such that
\begin{equation}\label{3.7}
\int\limits_{\mathbb R_+}\frac{\tau(du)}{1+u}<\infty,
\end{equation}
if and only if $f(z)\in\mathcal N$ and $f(z)$ is analytic and nonnegative
on $(-\infty,0)$.
\end{theorem}

\begin{theorem}\label{3.2th}
A~function $f(z)\in \mathcal N$ be analytic and nonnegative
on $(-\infty,0)$ if and only if $zf(z)$ is in $\mathcal N$. 
\end{theorem}

\begin{corollary}\label{3.3co}
\begin{enumerate}
\item A~function $f(z)\in \mathcal N$ is analytic and nonpositive
on $(-\infty,0)$ if and only if $f(z)/z$ is in $\mathcal N$. 
\item A~function $f(z)/z$ admits a~representation $(\ref{3.6})$ with 
a~nonnegative measure $\tau$ satisfying assumptions $(\ref{3.7})$ and 
$\tau(\{0\})=0$ if and only if $f(z)\in\mathcal K$.
\end{enumerate}
\end{corollary}
\begin{Proof}
At first we shall prove the~assertion (1).
Let the~function $f(z)\in \mathcal N$ be analytic and nonpositive 
on $(-\infty,0)$. Then
the~function $-1/f(z)\in\Cal N$ and satisfies the~assumptions
of Theorem~\ref{3.2th}. By this theorem, $-z/f(z)\in\mathcal N$ and therefore
$f(z)/z\in \mathcal N$. The~converse assertion follows by repeating
the~previous arguments. 

Let us now prove the~assertion (2). It easy to see that
$f(z)\in\mathcal K$ if the~function $f(z)/z$ admits a~representation 
$(\ref{3.6})$ with a~nonnegative measure $\tau$ which satisfies assumption 
$(\ref{3.7})$ and $\tau(\{0\})=0$. The~converse assertion follows from
the~assertion (1) of Corollary~\ref{3.3co} and Theorem~\ref{3.1th}.
\end{Proof}

\begin{theorem}\label{3.3th}
Any function $f(z)$ of class $\mathcal N$ non identically zero admits a~unique
multiplicative representation
$$
f(z)=C\exp\Big\{\int\limits_{-\infty}^{\infty}\Big(\frac 1{t-z}-\frac t{1+t^2}
\Big)p(t)\,dt,
$$
where $C>0,\,p(t)$ is a~summable function such that $0\le p(t)\le 1$
almost every where and
$$
\int\limits_{-\infty}^{\infty}\frac{p(t)}{1+t^2}\,dt<\infty.
$$
If $f(z)$ is analytic and nonnegative on the~negative real axis $(-\infty,0)$,
then $p(t)=0$ for $t<0$.
\end{theorem}

\begin{proposition}\label{3.4pro}
Let $f(z)\in\mathcal N$ be analytic and nonpositive on $(-\infty,0)$. 
Then $f(z)$ is univalent in the~left half-plane $i\mathbb C^+$.
\end{proposition}
\begin{Proof}
Since the~function $f(z)/z$ admits representation $(\ref{3.6})$, 
we see that
$$
f(z_1)-f(z_2)=(z_1-z_2)\Big(a+\int\limits_{(0,\infty)}\frac u{(u-z_1)
(u-z_2)}\,\tau(du)\Big)
$$
for all $z_1,z_2\in\mathbb C^+$. Since for $z_1\ne z_2$ and $z_1,z_2\in
\mathbb C^+\cap i\mathbb C^+$  
\begin{align}
\IM \int\limits_{(0,\infty)}&\frac u{(u-z_1)(u-z_2)}\,\tau(du)\notag\\&=
\int\limits_{(0,\infty)}\frac{-(\RE z_1\IM z_2+\RE z_2\IM z_1)
+u(\IM z_1+\IM z_2)}{|u-z_1|^2|u-z_2|^2}\,u\tau(du)\ne 0,\notag
\end{align}
we conclude from the preceding formula that $f(z_1)\ne f(z_2)$ for
the considered $z_1$ and $z_2$. Hence $f(z)$ is
univalent in $\mathbb C^+\cap i\mathbb C^+$. Since $f(\mathbb C^+)
\subset\mathbb C^+$ and $f$ is strictly increasing on $(-\infty,0)$,
the~univalence of $f$ on $i\mathbb C^+$ follows from the~identity
$f(\bar z)=\overline{f(z)}$. 
\end{Proof}
We need as well the~following well-known result about Schur functions 
$Q_{\mu}\in\Cal S_*$ (the~definition of  $Q_{\mu}$ see in (\ref{2.13})).
\begin{proposition}\label{3.5pro}
Let $\mu\in\mathcal S_*$, then
$$
Q_{\mu}(z)=z\ffrac{Q_{\mu}'(0)+z\varphi_1(z)}{1
+\overline{Q_{\mu}'(0)} z\varphi_1(z)},
$$
where $\varphi_1(z)\in\mathcal S$. 
\end{proposition}

We need the~following result about the~behavior of nondecreasing
functions (see \cite{LinOst:1977}, Ch.~4, \S 16).

\begin{proposition}\label{3.6pro}
Let $\nu$ be a~nonnegative finite measure. If $\nu([a,b])>0$ 
and $\nu([a,x))$, $a<x<b$, is a~continuous function, 
then there exists $[\alpha,\beta]\subset [a,b]$
such that $\nu ([\alpha,\alpha+h])\,\,>ch$ and $\nu ([\beta-h,\beta])>ch$
for all $0<h\le h_0$, where $c>0$ and $h_0>0$ depend on the~measure $\nu$.
\end{proposition}

We also need some results of the~theory of real and complex variable
(see \cite{Go:1969} (1969), \cite{Mar:1965} (1965)).

\begin{theorem}\label{3.7th} $($Weierstrass' preparation theorem$)$
Let $F(z,w)$ be a~function of two complex variables which is analytic
in neighborhood $|z-z_0|<r,\,|w-w_0|<\rho$ of the~point $(z_0,w_0)$,
and suppose that
\begin{equation}\label{3.8}
F(z_0,w_0)=0,\qquad F(z_0,w)\not\equiv 0.
\end{equation}
Then there is a~neighborhood $|z-z_0|<r'<r,\,|w-w_0|<\rho'<\rho$
where $F(z,w)$ may be written as
\begin{equation}\label{3.9}
F(z,w)=\big(A_0(z)+A_1(z)w+\dots+A_{k-1}(z)w^{k-1}+w^k\big)G(z,w),
\end{equation}
$k$ is determined by
$$
\frac{\partial F(z_0,w_0)}{\partial w}=\dots
\frac{\partial^{k-1} F(z_0,w_0)}{\partial w^{k-1}}=0,\qquad
\frac{\partial^k F(z_0,w_0)}{\partial w^k}\ne 0,
$$
the~functions $A_0(z),A_1(z),\dots,A_{k-1}(z)$ are analytic 
for $|z-z_0|<r'$, and the~function $G(z,w)$ is analytic and nonzero
on the~set $|z-z_0|<r',\,|w-w_0|<\rho'$.
\end{theorem}

\begin{theorem}\label{3.8th}
Let $\{u_n(z)\}_{n=1}^{\infty}$ denote a~sequence of functions
which are harmonic on a~domain $B$. If $\{u_n(z)\}_{n=1}^{\infty}$
is uniformly bounded in the~interior of $B$, it contains
a~subsequence that converges uniformly in the~interior of $B$
to a~harmonic function on $B$.
\end{theorem}

\begin{theorem}\label{3.8ath}
Let $\{u_n(z)\}_{n=1}^{\infty}$ denote a~sequence of functions
which are harmonic on a~domain $B$ such that $u_n(z)\le u_{n+1}(z)$
for all $z\in B$. Suppose that the~sequence converges at some
point $z_0\in B$. Then the~sequence converges uniformly in the~interior 
of $B$ to a~harmonic function on $B$.
\end{theorem}

\begin{theorem}\label{3.9th}
Given a~domain $B$ and a~sequence $\{f_n(z)\}_{n=1}^{\infty}$
of regular functions on $B$, suppose the~sequence 
$\{u_n(z)\}_{n=1}^{\infty}=\{\RE f_n(z)\}_{n=1}^{\infty}$
converges uniformly on every compact subset of $B$, and suppose
$\{f_n(z)\}_{n=1}^{\infty}$ converges at some point $z_0\in B$.
Then $\{f_n(z)\}_{n=1}^{\infty}$ converges uniformly on every compact 
subset of $B$ to a~regular function.
\end{theorem}

\begin{theorem}\label{3.10th} $($Vitali$)$ 
Let $\{f_n(z)\}_{n=1}^{\infty}$ denote a~sequence
function that are re\-gu\-lar on $B$. Suppose that $\{f_n(z)\}_{n=1}^{\infty}$
is uniformly bounded in the~interior of $B$ and converges
on a~set of points $z_k\in B,\,k=1,2,\dots$, that has a~cluster
point in the~interior of $B$. Then, the~sequence $\{f_n(z)\}_{n=1}^{\infty}$
converges uniformly in the~interior of $B$. 

\end{theorem}

In the~sequel we shall need the~following result.

\begin{proposition}\label{3.11pro}
Let $\big\{\mu_n\big\}_{n=1}^{\infty}$ be a~sequence of 
p-measures on $\mathbb R$. The~fo\-llowing assertions are equivalent.
\begin{enumerate}
\item  The~sequence $\big\{\mu_n\big\}_{n=1}^{\infty}$ converges  
weakly to a~p-measure $\mu$.
\item  There exist $\alpha,\beta>0$ such that the~sequence
$\big\{\varphi_{\mu_n}\big\}_{n=1}^{\infty}$ converges uniformly
on the~compact subsets of $\Gamma_{\alpha,\beta}$ to a~function 
$\varphi$, and $\varphi_{\mu_n}(z)=o(|z|)$ uniformly in $n$ as
$z\to \infty,\,z\in \Gamma_{\alpha,\beta}$. 
\item
There exist $\alpha',\beta'>0$ such that the~sequence
$\big\{\varphi_{\mu_n}\big\}_{n=1}^{\infty}$ converges uniformly
on the~compact subsets of $\Gamma_{\alpha',\beta'}$ to a~function 
$\varphi$, and $\Im\varphi_{\mu_n}(iy)=o(y)$ uniformly in $n$ as
$y\to+\infty$. 
\end{enumerate}
Moreover, if $(1)$ and
$(2)$ are satisfied, we have $\varphi=\varphi_{\mu}$ in 
$\Gamma_{\alpha,\beta}$. 
\end{proposition}

This result without the~assertion (3) was proved by
Bercovici and Voiculescu~\cite{BercVo:1993}. 
We should prove that the~assertion (3) implies (1) only.
\begin{proof}
Let $(3)$ hold. We should prove that there exists $\mu\in\mathcal M$ such that
$\varphi(z)=\varphi_{\mu}(z)$ and $\{\mu_n\}_{n=1}^{\infty}$ converges weakly to
$\mu$. Since $\Im \varphi_{\mu_n}(iy)=o(y)$ iniformly in $n$ as $y\to+\infty$,
we obtain the~relation
\begin{equation}\label{3.11,1}
F_{\mu_n}(iy+\varphi_{\mu_n}(iy))=iy
\end{equation}
for sufficiently large $y\ge y_0>0$ and $n\ge 1$. The~functions
$F_{\mu_n}(z)\in\mathcal F$ therefore
\begin{equation}\label{3.11,2}
F_{\mu_n}(z)=z+a_n+\int\limits_{\Bbb R}\frac{1+uz}{u-z}\,\sigma_n(du),
\quad z\in\mathbb C^+,
\end{equation}
where $a_n\in\mathbb R$ and $\sigma_n$ are finite nonnegative measures.
Rewrite (\ref{3.11,1}) in the~form
\begin{equation}\label{3.11,3}
\varphi_{\mu_n}(iy)+a_n+\int\limits_{\Bbb R}\frac{1+u(iy+\varphi_{\mu_n}(iy))}
{u-(iy+\varphi_{\mu_n}(iy))}\,\sigma_n(du)=0, \quad y\ge y_0.
\end{equation}
From this relation it follows, for $y\ge y_0$,
\begin{equation}\label{3.11,4}
-\frac{\Im \varphi_{\mu_n}(iy)}{y+\Im\varphi_{\mu_n}(iy)}=\int\limits_{\Bbb R}
\frac{1+u^2}{(u-\Re \varphi_{\mu_n}(iy))^2+
(y+\Im\varphi_{\mu_n}(iy))^2}\,\sigma_n(du).
\end{equation}
Since $y+\Im\varphi_{\mu_n}(iy)=y(1+o(1))$ uniformly in $n$
for sufficiently large $y$
and $\varphi_{\mu_n}(iy)\to\varphi(iy)$ as $n\to\infty$, we obtain from 
(\ref{3.11,4}) for sufficiently large $y$: $\sigma_n(\mathbb R)
\le c_1(y,\varphi)$, $n\ge 1$. Then we deduce from (\ref{3.11,3}) that 
$|a_n|\le c_2(y,\varphi),\,n\ge 1$. Here $c_j(y,\varphi),\,j=1,2$, 
are positive constants depended on $y,\varphi$ only. 
By the~vague compactness theorem (see \cite{Lo:1963}, p.179),
there exists a~subsequence $\{n'\}$ such that
$$
a_{n'}\to a,\quad \sigma_{n'}(\mathbb R)\to b,
$$
where $a\in\mathbb R, b\in\mathbb R_+$ and $\{\sigma_{n'}\}$ converges 
in the~vague topologue to some nonnegative mesure $\sigma$ such that
$\sigma(\mathbb R)\le b$. Using the~Helly-Brey lemma (\cite {Lo:1963}, p.181)
we obtain from (\ref{3.11,2}) 
\begin{align}\label{3.11,5}
F_{\mu_{n'}}(z)\to F(z)&:=z+a+\int\limits_{\mathbb R}\Big(\frac{1+uz}{u-z}-z\Big)
\,\sigma(du)+bz\notag\\
&=z+a+\int\limits_{\mathbb R}\frac{1+uz}{u-z}\,\sigma(du)
+(b-\sigma(\mathbb R))z,\quad n'\to\infty,
\end{align}
uniformly on every compact set in $\mathbb C^+$. It is easy to see from
(\ref{3.11,3}) that
\begin{equation}\label{3.11,6}
\varphi(iy)+a+\int\limits_{\mathbb R}\frac{1+u(iy+\varphi(iy))}
{u-(iy+\varphi(iy))}\,\sigma(du)+(b-\sigma(\Bbb R))(iy+\varphi(iy))=0, 
\quad y\ge y_0.
\end{equation}
We deduce from this equality that $(b-\sigma(\mathbb R))y+(1+b-\sigma(\mathbb R))
\Im\varphi(iy))\le 0$ for $y\ge y_0$. This implies $\sigma(\Bbb R)=b
=\lim_{n'\to\infty}\sigma_{n'}(\mathbb R)$. Hence $\{\sigma_{n'}\}$ converges 
in the~weak topologue to the~mesure $\sigma$. In addition
there exists $\mu\in\mathcal M$
such that $F(z)=F_{\mu}(z),\,z\in\mathbb C^+$, and we have
$$
F_{\mu_{n'}}(z)=z(1+o(1)),\quad F_{\mu}(z)=z(1+o(1))\quad \text{as}\quad
z\to\infty, \,\,z\in\Gamma_{\alpha,\beta}, 
$$
uniformly in $n'$ for some $\alpha>0,\beta>0$. Therefore $\varphi_{\mu_{n'}}
(z)=o(z)$ uniformly in $n'$ as $z\to\infty,\,z\in\Gamma_{\alpha,\beta}$.
Hence assumptions of the~assertion (2) hold and, by the~equivalence of (1)
and (2), $\{\mu_{n'}\}$ converges 
weakly to $\mu$. In other words we proved that under assumptions of 
the~assertion (3) we can choose a~subsequence $\{n'\}$ such that
$\{\mu_{n'}\}$ converges weakly to some p-measure $\mu$.
It remains to show that $\{\mu_{n}\}$ converges 
weakly to $\mu$. Let to the~contrary $\{\mu_{n}\}$ does not converge
weakly to $\mu$. This means that there exists $\{n''\}$ such that
$\{\mu_{n''}\}$ converges weakly to $\nu\in\mathcal M$ such that $\nu\not\equiv
\mu$. On the~other hand, as it is follows from (\ref{3.11,1}), 
$F_{\mu}(z)$ and $F_{\nu}(z)$ satisfy the~equation
$$
F_{\mu}(iy+\varphi(iy))=iy,\quad F_{\nu}(iy+\varphi(iy))=iy,\quad
y\ge y_0,
$$
and we have $F_{\mu}(iy+\varphi(iy))=F_{\nu}(iy+\varphi(iy))$ for
$y\ge y_0$. This implies $F_{\mu}(z)=F_{\nu}(z),\,z\in\Bbb C^+$, and
$\mu\equiv \nu$, a~contradiction. The~proposition is proved.
\end{proof}

We shall need the~following result
of Bercovici and Voiculescu~\cite{BercVo:1993} as well.

\begin{proposition}\label{3.12pro}
Let $\{\mu_n\}_{n=1}^{\infty}$ be a~tight sequence of 
p-measures on $\Bbb R_+$ such that $\delta_0$ is not in the~weak closure
of $\{\mu_n\}_{n=1}^{\infty}$. The~following assertions are equivalent.
\begin{enumerate}
\item The~sequence $\{\mu_n\}_{n=1}^{\infty}$ converges in 
the~weak topology to a~measure $\mu$.

\item There exist numbers $\alpha\in (0,\pi)$ and
$0<\beta<\Delta$ such that 
the~sequence $\{\Sigma_{\mu_n}\}_{n=1}^{\infty}$ converges uniformly
on $\Gamma_{\alpha,\beta,\Delta}^+$ to a~function $\Sigma$.
\end{enumerate}
Moreover, if $(1)$ and $(2)$ are satisfied, we have 
$\Sigma=\Sigma_{\mu}$ in $\Gamma_{\alpha,\beta,\Delta}^+$.
\end{proposition}

We need the~following result of Bercovici and Voiculescu~\cite{BercVo:1992}
as well.

For positive number $\alpha$ denote $\mathbb D_{\alpha}
:=\{z\in\mathbb C:|z|<\alpha\}$.
\begin{proposition}\label{3.13pro}
Consider a~measure $\mu\in\Cal M_*$ and a~sequence $\mu_j\in\mathcal M_*,\,
j=1,\dots$. The~following assertions are equivalent.
\begin{enumerate}
\item
The~sequence $\{\mu_n\}_{n=1}^{\infty}$ converges in 
the~weak topology to a~measure $\mu\in\mathcal M_*$.
\item
There exists a~positive number $\alpha$ such that the~sequence
$\{\Sigma_{\mu_j}\}_{j=1}^{\infty}$ converges uniformly 
on $\Bbb D_{\alpha}$ to a~function $\Sigma_{\mu}$.
\end{enumerate}
Moreover, if $(1)$ and $(2)$
are satisfied, we have $\Sigma=\Sigma_{\mu}$ in $\mathbb D_{\alpha}$.
\end{proposition}

\section{ Additive free convolution}

In this section we prove Theorem~2.1 and its consequences.
We need the~following auxiliary results.

\begin{lemma}\label{4.1l}
Let $g:\mathbb C^+\to \mathbb C^-$ be analytic with
\begin{equation}\label{4.1}
\liminf_{y\to+\infty}\frac {|g(iy)|}y=0.
\end{equation}
Then the~function $f:\mathbb C^+\to\mathbb C$ defined via $z\mapsto z+g(z)$
takes every value in $\mathbb C^+$ precisely once. The~inverse
$f^{(-1)}:\mathbb C^+\to \mathbb C^+$ thus defined is in the~class $\mathcal F$.
\end{lemma}

This lemma generalizes a~result of Maassen~\cite{Maa:1992} (see Lemma~2.3).
Maassen proved Lemma~\ref{4.1l} under the~additional restriction $|g(z)|
\le c/\IM z$ for $z\in\mathbb C^+$, where $c$ is a~constant.
\begin{Proof}
Since $g:\mathbb C^+\to\mathbb C^-$ is analytic, 
it can be written in Nevanlinna's integral form 
(see (\ref{3.3}) in Section~3) 
\begin{equation}\label{4.2}
g(z)=a-bz-\int\limits_{\mathbb R}\frac{1+uz}{u-z}
\,\sigma(du),
\end{equation}
where $a,b\in\mathbb R$, $b\ge 0$ and $\sigma$ is a~finite 
nonnegative measure. By (\ref{4.1}), $b=0$. 
Denote $\alpha_c:=\sigma(\{|u|>c\})$, where $c>0$ is chosen
such that $\alpha_c<1$. We may decompose $1+uz$ in the~integral~(4.2)
as $1+u^2+u(z-u)$ for $|u|\le c$ and as $1+z^2+z(u-z)$ for $|u|>c$.
Hence we get
$
f(z)=a_c+f_1(z)+f_2(z)
$
for $z\in\mathbb C^+$, where $a_c:=a+\int_{[-c,c]}u\,\sigma(du)$ and
$$
f_1(z):=(1-\alpha_c)z-\int\limits_{[-c,c]}\frac {1+u^2}{u-z}
\,\sigma(du),\quad 
f_2(z):=-(1+z^2)\int\limits_{|u|>c}\frac {\sigma(du)}{u-z}.
$$
Let $w\in\mathbb C^+$, denote $w_1:=w-a_c$.
For every fixed $w\in\mathbb C^+$
we consider a~closed rectifiable curve $\gamma_1=\gamma_1(w)$ 
(see Figure~1)
consisting of some smooth curve $\gamma_{1,1}$ connecting $w_1-R$
to $w_1+R$ inside the~strip $0<\IM z<\IM w$,
the~arc $\gamma_{1,2}:0<\arg(z-w_1)\le \eta$ on the~circle $|z-w_1|=R$
connecting $w_1+R$ to $w_1+Re^{i\eta}$, the~arc 
$\gamma_{1,3}:\eta<\arg(z-w_1)\le \pi-\eta$ on the~circle $|z-w_1|=R$
connecting $w_1+Re^{i\eta}$ to $w_1-Re^{-i\eta}$, and
the~arc $\gamma_{1,4}:\pi-\eta<\arg(z-w_1)
\le \pi$ on the~circle $|z-w_1|=R$ connecting $w_1-Re^{-i\eta}$ 
to $w_1-R$. Here $\eta$ is given by $\eta:=10^{-2}\min\{\arg w_1,-\arg 
\overline w_1\}$.
We also assume that $R>0$ is sufficiently large.


\begin{figure}[htbp]
 
\input{bild1.pstex_t}
 
\caption{}
\end{figure}


Note by (\ref{4.2}) with $b=0$ that $\max_{z\in\gamma_{1,3}}
|g(z)|/|z|\to 0$ as $R\to\infty$. We also see that on $\gamma_{1,3}$ 
$\IM z\ge \IM w+R\sin\eta$. Since $-\Im g(z)=o(R),\,z\in\gamma_{1,3}$, 
we have $\IM f(z)\ge \IM w+R\sin\eta-o(R)>\IM w,\,z\in\gamma_{1,3}$. 
Therefore, if $z$ runs  through $\gamma_{1,3}$ the~image $f(z)$ lies in 
the~half-plane $\IM z>\Im w$.

For $z\in\gamma_{1,1}$, $\IM(z+g(z))<\IM z\le\IM w$. Therefore,
if $z$ runs through $\gamma_{1,1}$ the~image $f(z)$
lies in the~half-plane $\IM z<\IM w$. 

Let $z\in\gamma_{1,2}$. It is easy to see that, for $|z|>c$, 
$$
f_1(z)=(1-\alpha_c)z\Big(1+\frac 1{1-\alpha_c}
\sum_{k=1}^{\infty}\frac 1{z^{k+1}}
\int\limits_{[-c,c]}u^{k-1}(1+u^2)\,\sigma(du)\Big)
$$
Therefore we obtain the~following formula, for $|z|>2(c^2+1)
(\sigma(\mathbb R)+1)/(1-\alpha_c)$,
$$
\log f_1(z)=\log (1-\alpha_c)+\log z+\sum_{m=1}^{\infty}
\frac{a_m}{z^{m+1}},
$$
where $a_m$ are real coefficients such that $|a_m|\le K^m,\,m=1,\dots$,
with some positive constant $K$. Here and in the~sequel we choose
the~principle branch of the~logarithm. Hence, for $|z|>2K$, 
\begin{align}\label{4.3}
\arg f_1(z)&=\arg z-\sum_{m=1}^{\infty}
\frac{a_m\sin \big((m+1)\arg z\big)}{|z|^{m+1}}\\
&=\Big(1+\Theta\sum_{m=1}^{\infty}\frac{K^m(m+1)}{|z|^{m+1}}\Big)\arg z
=\Big(1+6\Theta\frac K{|z|^2}\Big)\arg z,\notag
\end{align}
where $\Theta$ denotes a~real-valued quantity such that $|\Theta|\le 1$.
On the~other hand we easily obtain, for $|z|>2$,
$$
\arg (1+z^2)=2\Big(1+\frac{2\Theta}{|z|^2}\Big)\arg z.
$$
Therefore we conclude from the~definition of $f_2(z)$, taking into account
that $-\int_{|u|>c}\frac{\sigma(du)}{u-z}\in\mathbb C^-$,
\begin{equation}\label{4.4}
-\pi+2\Big(1+\frac{2\Theta}{|z|^2}\Big)\arg z\le\arg f_2(z)
\le 2\Big(1+\frac{2\Theta}{|z|^2}\Big)\arg z,\quad z\in\gamma_{1,2}.
\end{equation}
From (\ref{4.3}) and (\ref{4.4}) it follows that, for $z\in\gamma_{1,2}$, 
$|\arg f_1(z)-\arg f_2(z)|<\pi$ and therefore
\begin{equation}\label{4.5}
-\pi+2\Big(1+\frac{2\Theta}{|z|^2}\Big)\arg z\le\arg (f_1(z)+f_2(z))
\le 2\Big(1+\frac{2\Theta}{|z|^2}\Big)\arg z,\quad z\in\gamma_{1,2}.
\end{equation}

We conclude from (\ref{4.5}) that the~image $\zeta=f(z)$ lies
in the~domain $D_1:=\{\zeta\in\mathbb C:-\pi<\arg \zeta<3\eta\}$ when 
$z$ runs through $\gamma_{1,2}$.
In addition the~point $w$ does not lie in $D_1$ by the~choice of
the~parameter $\eta$. 
In the~same way we deduce that the~image $\zeta=f(z)$ lies in 
$D_2:=\{\zeta\in\mathbb C:\pi-3\eta<\arg \zeta<2\pi\}$ and 
$w\not\in D_2$ when $z$ runs through $\gamma_{1,4}$. 

Hence $f(z)$ winds around $w$ once, 
and it follows from the~argument
principle that 
inside the~curve $\gamma_1$ there is  a~unique 
point $z_0$ such that $f(z_0)=w$. Since this relation holds for all
sufficiently large $R>0$ and all curves $\gamma_{1,1}$, we deduce 
that the~point $z_0$ is unique in $\mathbb C^+$.

Hence the~inverse function $f^{(-1)}:\mathbb C^+\to\Bbb C^+$ exists 
and is analytic in $\mathbb C^+$. By condition (\ref{4.1}), $\lim_{y\to
+\infty}(iy/f^{(-1)}(iy))=1$ and therefore $f^{(-1)}\in\Cal F$.
This proves the~lemma. 
\end{Proof}

Let $z_1\in\mathbb C^+$ and $z_2\in\mathbb C^+$, and introduce
the functions
$$
w_1(z_1,z_2):=z_1+z_2-F_{\mu_2}(z_2),\qquad 
w_2(z_1,z_2):=z_1+z_2-F_{\mu_1}(z_1).
$$ 
\begin{lemma}\label{4.2l}
For every $z\in\Bbb C^+$ there exist unique points $z_1\in\mathbb C^+$ 
and $z_2\in\mathbb C^+$ such that
\begin{equation}\label{4.6}
z=w_1(z_1,z_2)\quad\text{and}\quad z=w_2(z_1,z_2).
\end{equation}
\end{lemma}
\begin{Proof}
Let us fix $z\in\mathbb C^+$. For every $z_2\in\mathbb C^+$ we define 
$z_1:=z-(z_2-F_{\mu_2}(z_2))$. Recall that $F_{\mu_j}\in\mathcal F,\,j=1,2$. 
Since, by (\ref{3.4b}), $(z_2-F_{\mu_2}(z_2))\in\mathbb C^-\cup\mathbb R$, 
it follows that $z_1\in\mathbb C^+$. Hence it suffices to prove that 
the~equation
\begin{equation}\label{4.7}
F_{\mu_2}(z_2)=F_{\mu_1}\big(z-(z_2-F_{\mu_2}(z_2))\big)
\end{equation}
has a~unique solution $z_2\in\mathbb C^+$. Rewrite (\ref{4.7}) in the~form
$$
z=z_2-g_{\mu_1}\big(z-(z_2-F_{\mu_2}(z_2))\big),
$$
where $g_{\mu_1}(w):=F_{\mu_1}(w)-w$ for $w\in\mathbb C^+$. 

By (\ref{3.3}), the~functions $g_{\mu_j}(w),\,j=1,2$, admit the representation
$$
g_{\mu_j}(w)=a_j+\int\limits_{\Bbb R}\frac{1+uw}{u-w}\,\tau_j(du),
\quad w\in\mathbb C^+,
$$
where $a_j\in\mathbb R$ and $\tau_j$ are finite nonnegative measures.
Therefore
$$
\Im g_{\mu_j}(iy)=y \int\limits_{\mathbb R}\frac{1+u^2}{u^2+y^2}
\,\tau_j(du) 
$$
and
$$
\Re g_{\mu_j}(iy)=a_j+\int\limits_{\mathbb R}\frac{u(1-y^2)}
{u^2+y^2}\,\tau_j(du). 
$$
For $y\ge 1$, we obtain the following estimates
$$
\Im g_{\mu_j}(iy)\ge \frac 1{2y}\int\limits_{[-y,y]}
(1+u^2)\,\tau_j(du)+\frac y2\int\limits_{|u|>y}\tau_j(du)
$$
and
$$
|\Re g_{\mu_j}(iy)|\le |a_j|+\int\limits_{[-y,y]}
|u|\,\tau_j(du)+y^2\int\limits_{|u|>y}\frac{\tau_j(du)}{u}.
$$
We conclude from the~last two inequalities and Lyapunov's
inequality that
\begin{equation}\label{4.7aa}
|\Re g_{\mu_j}(iy)|\le c\big(1+(y \Im g_{\mu_j}(iy))^{1/2}
+\Im g_{\mu_j}(iy)\big),\quad y\ge 1,
\end{equation}
where $c$ is a~positive constant.

We shall prove that
\begin{equation}\label{4.7a}
|g_{\mu_1}\big(z-(i\IM z_2-F_{\mu_2}(i\IM z_2))\big)|/\IM z_2\to 0\quad
\text{as}\quad\Im z_2\to +\infty.
\end{equation}
Let to the~contrary 
\begin{equation}\label{4.7b}
|g_{\mu_1}\big(z-(iy_k-F_{\mu_2}(iy_k)\big)|/y_k\ge c>0
\end{equation}
for some sequence $\{y_k\}_{k=1}^{\infty}$ such that 
$y_k\to+\infty$ and for a~positive constant $c$.

If $\liminf_{y_k\to+\infty}|iy_k-F_{\mu_2}(iy_k)|<\infty$, then
there exists a~subsequence $\{y_k'\}_{k=1}^{\infty}
\subset\{y_k\}_{k=1}^{\infty}$ such that $\lim_{y_k'\to+\infty}
|iy_k'-F_{\mu_2}(iy_k')|<\infty$. 
It is easy to see, $\IM (z-(iy_k'-F_{\mu_2}(iy_k'))\ge\IM z$ and
$|z-(iy_k'-F_{\mu_2}(iy_k')|\le |z|+c_1$ for all $y_k$ and for some
constant $c_1>0$. Hence in this case we have 
$|g_{\mu_1}\big(z-(iy_k'-F_{\mu_2}(iy_k')\big)|\le c_2$ for all $y_k'$ 
and for some constant $c_2>0$. This estimate contradicts to (\ref{4.7b}).

Let $\liminf_{y_k\to+\infty}|iy_k-F_{\mu_2}(iy_k)|=\infty$ and 
$\liminf_{y_k\to+\infty}|\Re(iy_k-F_{\mu_2}(iy_k))|=\infty$.
We see that
\begin{align}
&|g_{\mu_1}(z+g_{\mu_2}(iy_k))|\le |a_1|+\frac{\tau_1(\mathbb R)}{\Im z}+
|z+g_{\mu_2}(iy_k)|\int\limits_{\mathbb R}\frac{|u|\,\tau_1(du)}
{|u-z-g_{\mu_2}(iy_k)|}\notag\\
&\le |a_1|+\frac{\tau_1(\mathbb R)}{\Im z}+|z+g_{\mu_2}(iy_k)|\big(I_1(iy_k)
+I_2(iy_k)+I_3(iy_k)\big),\notag
\end{align}
where
$$
I_1(iy_k):=\int\limits_{|u|\le x_k/2}\frac{|u|\,\tau_1(du)}
{|u-z-g_{\mu_2}(iy_k)|}\le c\frac{|\Re g_{\mu_2}(iy_k)|}{|g_{\mu_2}(iy_k)|}
\le c,
$$
$$
I_2(iy_k):=\int\limits_{x_k/2<|u|\le 2x_k}
\frac{|u|\,\tau_1(du)}{|u-z-g_{\mu_2}(iy_k)|}\le c\frac{|\Re g_{\mu_2}(iy_k)|}
{\Im z+\Im g_{\mu_2}(iy_k)}\tau_1(\{|u|>x_k/2\}),
$$
and
$$
I_3(iy_k):=\int\limits_{|u|>2x_k}
\frac{|u|\,\tau_1(du)}{|u-z-g_{\mu_2}(iy_k)|}\le c
$$
with some positive constant and $x_k:=|\Re g_{\mu_2}(iy_k)|$.
Using (\ref{4.7aa}), we finally obtain
$$
|g_{\mu_1}(z+g_{\mu_2}(iy_k))|\le  |a_1|+\frac{\tau_1(\Bbb R)}{\Im z}+
cy_k\tau_1(\{|u|>x_k/2\})+ c|z+g_{\mu_2}(iy_k)|.
$$
From this estimate it follows immediately that 
$$
\frac{|g_{\mu_1}(z-(iy_k-F_{\mu_2}(iy_k)))|}{y_k}\to 0\quad
\text{as}\quad y_k\to +\infty,
$$
a~contradiction with (\ref{4.7b}).

Let $\liminf_{y_k\to+\infty}|iy_k-F_{\mu_2}(iy_k)|=\infty$ and 
$\liminf_{y_k\to+\infty}|\Re(iy_k-F_{\mu_2}(iy_k))|<\infty$.
Without loss of generality we can assume that 
$\limsup_{y_k\to+\infty}|\Re(iy_k-F_{\mu_2}(iy_k))|<\infty$.

Since
$|g_{\mu_1}(z+z_2)|/\IM z_2\to 0$ as $z_2 \to +\infty$ nontangentially to $\mathbb R$
and
$|i\IM z_2-F_{\mu_2}(i\IM z_2)|/\IM z_2\to 0$ as $\Im z_2\to +\infty$, 
we see that  
$$
\frac{|g_{\mu_1}(z-(iy_k-F_{\mu_2}(iy_k)))|}{|z-(iy_k-F_{\mu_2}(iy_k))|}
\frac{|z-(iy_k-F_{\mu_2}(iy_k))|}{y_k}\to 0\quad
\text{as}\quad y_k\to +\infty,
$$
a~contradiction with (\ref{4.7b}). Hence (\ref{4.7a}) is proved.

Consider the~function $\tilde{f}(z_2):=z_2-\tilde{g}(z_2)$, where
$\tilde{g}(z_2):=g_{\mu_1}\big(z-(z_2-F_{\mu_2}(z_2))\big)$, 
$z_2\in\Bbb C^+$. By the~definition of the~function $g_{\mu_1}$, we see 
that $-\tilde{g}:\mathbb C^+\to\mathbb C^-$. By (\ref{4.7a}), the~$-\tilde{g}$
satisfies the~condition (\ref{4.1}).
Applying Lemma~\ref{4.1l} to the~function $\tilde{f}(z_2)$, we obtain 
that (\ref{4.7}) has a~unique solution 
$z_2\in\mathbb C^+$ for every fixed $z\in\mathbb C^+$, thus proving the~lemma.
\end{Proof}

{\it Proof of Theorem $2.1$.}
Consider the~function $F(z,z_2):=F_{\mu_2}(z_2)-F_{\mu_1}
\big(z-(z_2-F_{\mu_2}(z_2))\big)$ as a~function of the~two complex
variables $z$ and $z_2$. It is analytic on 
$\Bbb C^+\times\mathbb C^+$. By Lemma~\ref{4.2l}, for every fixed 
$z=z^0\in\mathbb C^+$
the~equation (\ref{4.7}) has an~unique solution $z_2^0\in\mathbb C^+$.
Hence $F(z^0,z_2^0)=0$. We shall verify that $F(z^0,z_2)\not\equiv 0$
for $z_2\in\Bbb C^+$.
Assume that $F(z^0,z_2)\equiv 0$ holds in $z_2\in\mathbb C^+$. 
Since $\Im F_{\mu_2}(iy)/y\to 1$ and $\Im F_{\mu_1}(z-(iy-
F_{\mu_2}(iy)))/y\to 0$ as $y\to+\infty$, by arguments as in Lemma~\ref{4.2l},
we arrive at a~contradiction.
Therefore the~function $F(z,z_2)$ satisfies the~assumption (\ref{3.8}) of
Theorem~\ref{3.7th} (Weierstrass' preparation theorem) at the~point 
$(z^0,z_2^0)$. 
Moreover, by this theorem and Lemma~\ref{4.2l}, 
the~function $F(z,z_2)$ admits
the~representation (\ref{3.9}) in a~neighborhood $|z-z^0|<r',\,
|z_2-z_2^0|<\rho'$ with the~positive integer $k=1$. 
Let us show that there exists $0<r''\le r'$ such that the~equation
$F(z,z_2)=0$ has a~unique root $z_2$ in $|z_2-z_2^0|<\rho'$ for any
given $z$ with $|z-z^0|<r''$. Since $F(z,z_2)=0$, (\ref{3.9}) implies 
that
$$
P(z,z_2):=A_0(z)+z_2=0
$$
for $z,z_2$ from the~above neighborhood. Here the~functions $A_0(z)$
is analytic in the~domain $|z-z^0|<r''$ and $A_0(z^0)=-z_2^0$.
For every $z$ with $|z-z^0|<r''$ and a~sufficiently small 
$r''=r''(z^0)\in(0,r')$,
this equation has a~root $z_2(z;z_2^0):=-A_{0}(z)$ 
in $|z_2-z_2^0|<\rho'$.

Thus, we have proved that for every given point $z^0\in\mathbb C^+$
there exists a~neighborhood $|z-z^0|<r''(z^0)$ such that (\ref{4.7})
has an~unique regular solution $z_2=z_2(z;z^0)$ with values in
$\Bbb C^+$.
Note that for points $z'\in\mathbb C^+$ and $z''\in\mathbb C^+,\,z'\ne z''$,
we have $z_2(z;z')=z_2(z;z'')$ for all $z$ in $\{|z-z'|<r''(z')\}\cap
\{|z-z''|<r''(z'')\}$.
By the~mo\-no\-dromy theorem (see \cite{Mar:1965}, v. 3, p. 269, 
\cite {NePa:1969}, p. 217), there exists a~regular function $Z_2(z),\,
z\in\Bbb C^+$, such that, for every point $z^0\in\mathbb C^+$,
$Z_2(z)=z_2(z;z^0)$ for $|z-z^0|<r''(z^0)$.
Therefore $Z_2(z)\in\Cal N$ and it is an~unique solution of (\ref{4.7}).

It is easy to see from (\ref{4.7}) that $Z_2(z)$ is in the~class $\mathcal F$. 
Indeed, we note that the~function $F_{\mu_2}(Z_2(z))-Z_2(z)\in\mathcal N$
and, by (\ref{3.4}), $F_{\mu_2}(Z_2(z))-Z_2(z)=(b+o(1))z$, where $b\ge 0$
is some constant, for $z\in\Bbb C^+$ such that $|\Re z|/\Im z$ stays 
bounded as $|z|$ tends to infinity. Hence the~function $z
+F_{\mu_2}(Z_2(z))-Z_2(z)\in\mathcal N$ and $|z+F_{\mu_2}(Z_2(z))-Z_2(z)|\to
\infty$ for the~same $z$. In addition $|\Re(z+F_{\mu_2}(Z_2(z))-Z_2(z))|/
\Im (z+F_{\mu_2}(Z_2(z))-Z_2(z))$ remains bounded as $|z|$ tends to infinity.
Therefore, by (\ref{3.4}),
$$
F_{\mu_1}(z+F_{\mu_2}(Z_2(z))-Z_2(z))=(z+F_{\mu_2}(Z_2(z))-Z_2(z)(1+o(1))
$$ 
for the~considered $z$. Using this relation, we conclude from (\ref{4.7})
and (\ref{3.4}) that
$$
Z_2(z)=(1+o(1))z+o(1)(F_{\mu_2}(Z_2(z))-Z_2(z))=(1+o(1))z
$$
for the~same $z$. Thus $Z_2\in\Cal F$. The~desired result is proved.

Choosing $Z_1(z):=z-Z_2(z)+F_{\mu_2}(Z_2(z)), \,z\in\mathbb C^+$, we see
that $Z_1$ and $Z_2$ are unique solutions of (\ref{4.6}) in the~class 
$\Cal F$. Hence Theorem~2.1 is proved.
$\square$

{\it Proof of Corollary $2.2$.}
For simplicity we shall prove this corollary in the~case $n=3$.
The~general case follows by induction.

Denote $\mu_{2,3}:=\mu_2\boxplus\mu_3$. We have, by associativity, 
$\mu_1\boxplus\mu_2\boxplus\mu_3=\mu_1\boxplus\mu_{2,3}$.
By Theorem~2.1, there exist unique functions $W_1$ and $W_{2,3}$
in the~class $\mathcal F$ such that, for $z\in\mathbb C^+$,
\begin{equation}\label{4.8}
z=W_1(z)+W_{2,3}(z)-F_{\mu_1}(W_1(z))\quad\text{and}\quad
F_{\mu_1}(W_1(z))=F_{\mu_{2,3}}(W_{2,3}(z)).
\end{equation}
On the~other hand, again by Theorem~2.1, there exist unique functions
$W_2\in\mathcal F$ and $W_3\in\mathcal F$ such that, for $z\in\mathbb C^+$,
$$
z=W_2(z)+W_3(z)-F_{\mu_2}(W_2(z))\quad\text{and}\quad
F_{\mu_2}(W_2(z))=F_{\mu_3}(W_3(z)).
$$
Hence, replacing $z$ by $W_{2,3}(z)$ in the~last equation we get, for all
$z\in\mathbb C^+$,
\begin{equation}\label{4.10}
W_{2,3}(z)=W_2(W_{2,3}(z))+W_3(W_{2,3}(z))-F_{\mu_2}(W_2(W_{2,3}(z)))
\end{equation}
and
\begin{equation}\label{4.11}
F_{\mu_2}(W_2(W_{2,3}(z)))=F_{\mu_3}(W_3(W_{2,3}(z))).
\end{equation}
Comparing (\ref{4.8}) and (\ref{4.10}), (\ref{4.11}), we obtain 
the~assertion of Corollary~2.2 with 
$Z_1(z)=W_1(z)$, $Z_2(z)=W_2(W_{2,3}(z))$ and
$Z_3(z)=W_3(W_{2,3}(z))$.
$\square$

Corollary~2.3 is an~obvious consequence of Corollary~2.2.
Note that the~continuous semigroup version of Corollary~2.3
with $t\ge 1$ replacing $n$ is proved using Lemma~4.1 and Lemma~4.2, 
and repeating the~arguments of Theorem~2.1. 

\section{Multiplicative free convolution}

1. Consider the~case of multiplicative convolution for p-measures 
of class $\mathcal M_+$. In order to prove Theorem~2.4
we need the following two auxiliary results.

Let $\mu_1,\mu_2\in\mathcal M_+$ and $w\in\mathbb C^+$. Introduce 
the~functions
$$
f_1(z)=w\frac{K_{\mu_1}(z)}{z},\quad
f_2(z)=\frac{K_{\mu_2}(f_1(z))}{f_1(z)},
\quad\text{and}\quad f_3(z)=\frac z{f_2(z)},
\quad z\in\mathbb C^+.
$$
\begin{lemma}
The~function $f_3:\mathbb C^+\to\mathbb C$  
takes the~value $w\in\mathbb C^+$ precisely once. Moreover, $f_3$ takes 
this value in $D_w:=\{z\in\mathbb C:\arg w\le\arg z<\pi\}$.
\end{lemma}
\begin{Proof}
We shall fix $w\in\mathbb C^+$. Let $\alpha\in(0,\arg w)$. 
Denote by $\gamma_2=\gamma_2(\alpha)$ 
(see Figure~2) the~closed rectifiable curve
consisting of the~line segment $\gamma_{2,1}: te^{i\alpha},\,1/R\le t\le R$, 
connecting $e^{i\alpha}/R$ to $Re^{i\alpha}$,  
the~arc $\gamma_{2,2}:\alpha<\arg z<\pi$ on the~circle $|z|=R$
connecting $Re^{i\alpha}$ to $-R$, 
the~line segment $\gamma_{2,3}:-R\le t\le -1/R$ connecting $-R$ to $-1/R$,
and the~arc $\gamma_{2,4}:\alpha<\arg z<\pi$ on the~circle $|z|=1/R$
connecting $-1/R$ to $e^{i\alpha}/R$. 
Here the~parameter $R>0$ and will be chosen later sufficiently large. 
Let $z$ run through $\gamma_2$ in the~counter clockwise direction.

\begin{figure}[htbp]
 
\input{bild2.pstex_t}
 
\caption{}
\end{figure}


Since the~function $K_{\mu_1}$ is in the~Krein class $\mathcal K$, 
we note that $\arg w\le\arg f_1(z)<\pi$ for all $z\in\mathbb C$ with 
$\arg w\le\arg z\le\pi$. Let $z\in\mathbb C$ 
such that $\arg w<\arg z<\pi+\arg w$. Since $K_{\mu_1}(\bar z)=\overline
{K_{\mu_1}(z)}$ and $0\le\arg(K_{\mu_1}(z)/z)<\pi-\arg z$ for $z\in\mathbb C^+$,
we conclude that $\IM f_1(z)>0$ on the~angular domain
$\{z\in\mathbb C:\arg w<\arg z<\pi+\arg w\}$. Therefore 
the~function $f_2(ze^{i\arg w})$ is in the~class $\mathcal N$. 

On the~other hand we see that $\arg w\le\arg f_1(z)
<\pi+\arg w$ for all $z\in\mathbb C^+$. 
Since $K_{\mu_2}(\bar z)=\overline{K_{\mu_2}(z)}$ and 
$0\le\arg(K_{\mu_2}(z)/z)<\pi-\arg z$ for $z\in\mathbb C^+$, we deduce that
$e^{i\arg w}f_2(z)$ is in the~class $\mathcal N$. 

Using the~representation (\ref{3.3}) for functions $f\in\mathcal N$, we note
that if $|f(re^{i\alpha})|/r\to 0$ as $r\to\infty$ for any fixed
$\alpha\in(0,\pi)$, then $|f(z)|/r\to 0$ as $r\to\infty$ uniformly
in the~angle $\delta\le\arg z\le\pi-\delta$ with any 
fixed $\delta\in(0,\pi)$. From this and from the~relations
$f_2(ze^{i\arg w})\in\mathcal N$ and $e^{i\arg w}f_2(z)\in\mathcal N$ 
we easily obtain that the~estimate
\begin{equation}\label {5.1}
\max_{z\in\gamma_{2,2}}\frac{|f_2(z)|}R\to 0,\quad R\to\infty,
\end{equation}
follows from the~estimate
\begin{equation}\label {5.2}
\frac{|f_2(-R)|}R\to 0,\quad R\to\infty.
\end{equation}

Let us prove (\ref{5.2}). At the~first step we shall consider the~behavior
of the~function $f_1(-R)$ for $R\ge 1$. The~functions
$K_{\mu_j}\in\mathcal K,\,j=1,2$, admit (see (\ref{2.8})) the~representation
\begin{equation}\label{5.3}
K_{\mu_j}(z)=a_jz+zg_j(-z)
:=a_jz+z\int\limits_{(0,\infty)}\frac{\tau_j(du)}{u-z},
\quad z\in\mathbb C^+,
\end{equation}
where $a_j\ge 0$ and the~nonnegative measures $\tau_j$ satisfy 
condition (\ref{2.9}). Hence we have
$$
f_1(-R)=w\Big(a_1+\int\limits_{(0,\infty)}\frac{\tau_1(du)}{u+R}\Big)=
w(a_1+g_1(R)).
$$
Note that $g_1(R)\to 0$ as $R\to\infty$ and, if $\tau_1\not\equiv 0$, then
\begin{equation}\label{5.4}
g_1(R)\ge\frac{\tau_1((0,R))}R\ge \frac cR 
\end{equation}
for sufficiently large $R\ge 1$. Here and in the~sequel 
we shall denote by $c$ positive constants which do not depend on $R$. 
In view of (\ref{5.3}), we obtain the~formula
\begin{equation}\label{5.5}
f_2(-R)=a_2+\int\limits_{(0,\infty)}\frac{\tau_2(du)}{u-w(a_1+g_1(R))}
=a_2+g_2(-w(a_1+g_1(R))).
\end{equation} 
If in (\ref{5.5}) $a_1>0$, then, it is easy
to see, that $|f_2(-R)|\le c$ and we arrive at (\ref{5.2}).
Let $a_1=0$. In this case we need an~upper bound for the~function 
$|g_2(-wg_1(R))|$ for large $R\ge 1$.

Since $a_1=0$, we conclude $\tau_1\not\equiv 0$ and write, using (\ref{5.4}),
\begin{align}\label{5.6}
&\IM g_2(-wg_1(R))=\int\limits_{(0,\infty)}\frac{g_1(R)\IM w}
{(u-g_1(R)\RE w)^2+(g_1(R)\IM w)^2}\,\tau_2(du)\notag\\
&\le\frac{\tau_2((0,-g_1(R)\log g_1(R)))}{g_1(R)\IM w}+
cg_1(R)\int\limits_{[-g_1(R)\log g_1(R),\infty)}
\frac{\tau_2(du)}{u^2}\notag\\
&\le c\Big(\frac{\tau_2((0,-g_1(R)\log g_1(R)))}{g_1(R)}
+\frac 1{g_1(R)(\log g_1(R))^2}\Big)
=o\Big(\frac 1{g_1(R)}\Big)=o(R)
\end{align}
as $ R\to\infty$. In addition we have
\begin{align}\label{5.7}
&|\RE g_2(-wg_1(R))|=\Big|\int\limits_{(0,\infty)}\frac{(u-g_1(R)\RE w)}
{(u-g_1(R)\RE w)^2+(g_1(R)\IM w)^2}\,\tau_2(du)\Big|\\
&\le c\Big(\frac{\tau_2((0,2|w|g_1(R)))}{g_1(R)}+
\int\limits_{[2|w|g_1(R),-g_1(R)\log g_1(R))}
\frac{\tau_2(du)}u\notag\\
&+\int\limits_{[-g_1(R)\log g_1(R),\infty)}\frac{\tau_2(du)}u\Big)
=o\Big(\frac 1{g_1(R)}\Big)=o(R),\quad R\to\infty.\notag
\end{align}
The~estimate~(\ref{5.2}) and, hence (\ref{5.1}), follows immediately 
from (\ref{5.6}) and (\ref{5.7}).

Now we shall prove that
\begin{equation}\label{5.8}
R\max_{z\in\gamma_{2,4}}|f_2(z)|\to \infty,\quad R\to\infty.
\end{equation}

We shall write $f_1(z)$ in the~form
\begin{align}\label{5.9}
f_1(z)&=w(\Re (K_{\mu_1}(z)/z)+i\IM (K_{\mu_1}(z)/z))\notag\\
&=w\Big(a_1+\int\limits_{(0,\infty)}\frac{(u-\RE z)}{(u-\RE z)^2+
(\IM z)^2}\,\tau_1(du)\notag\\
&+i\int\limits_{(0,\infty)}\frac{\IM z}
{(u-\RE z)^2+(\IM z)^2}\,\tau_1(du)\Big).
\end{align}

Let $|z|=1/R$ and $\eta|\RE z|\le \IM z$, where $\eta:=\min
\{\tan\alpha, 1/10\}$. 
For these $z$, using (\ref{5.9}), we obtain for sufficiently large $R\ge 1$
the~inequality
\begin{equation}\label{5.10}
\Big|\Re \frac{K_{\mu_1}(z)}z\Big|\le a_1
+\frac{\tau_1((0,2\IM z/\eta))}
{\IM z}+2\int\limits_{[2\IM z/\eta,\infty)}\frac{\tau_1(du)}u 
=o(R).
\end{equation}
On the~other hand  the~following lower bound holds
\begin{equation}\label{5.11}
\Big|\Re \frac{K_{\mu_1}(z)}z\Big|\ge a_1+\frac 12\int\limits_
{[2\IM z/\eta,\infty)}\frac{\tau_1(du)}u-
\frac{\tau_1((0,2\IM z/\eta))}{\IM z}. 
\end{equation}

In addition, using (\ref{5.9}), we deduce, for the~same $z$,
\begin{equation}\label{5.12}
\frac {\eta^2}{10} q_1(z)\le\IM \frac{K_{\mu_1}(z)}z\le 4q_1(z),
\end{equation}
where
\begin{equation}\label{5.13}
q_1(z):=\frac{\tau_1((0,2\IM z/\eta))}{\IM z}+ \IM z\int\limits_
{[2\IM z/\eta,\infty)}\frac{\tau_1(du)}{u^2}=o(R). 
\end{equation}
Comparing (\ref{5.11}) and the left-hand side of the~inequality (\ref{5.12}),  
we conclude that $|K_{\mu_1}(z)/z|\ge c$ for some positive constant $c$
and the~$z$ considered above. From (\ref{5.10}) 
and the right-hand side of inequality (\ref{5.12}) we see that $|f_1(z)|=o(R)$ 
for those $z$. Moreover by (\ref{5.3}) and by the~definition of $f_1(z)$,
it follows that $\arg w\le \arg f_1(z)\le \pi+\arg w$ for $z\in\mathbb C^+$. 
Hence, we get, for $|z|=1/R, \eta|\Re z|\le \Im z$, 
\begin{equation}\label{5.14}
c\le |f_1(z)|=o(R)\quad\text{and}\quad \arg w\le \arg f_1(z)\le \pi+\arg w,
\end{equation}
with some positive constant $c$.

Let $|z|=1/R$ and $\RE z\le 0,-\eta\RE z>\IM z$. Repeating
the~previous arguments we obtain the~estimate (\ref{5.14}) for such $z$.
Thus, (\ref{5.14}) holds for all $z\in\gamma_{2,4}$.

Using (\ref{5.3}), it is not difficult to conclude that $|K_{\mu_2}(z)|
\ge c(\delta_1,\delta_2)$ for $|z|\ge \delta_1$ and $\delta_2\le
\arg z\le 2\pi-\delta_2$, where $\delta_1>0$, $\delta_2\in(0,\pi)$ are
constants and $c(\delta_1,\delta_2)$ is a~positive constant depending
on $\delta_1$ and $\delta_2$. Using this estimate, we arrive at
the~lower bound  
$$
\Big|\frac{K_{\mu_2}(f_1(z))}{f_1(z)}\Big|\ge \frac c{|f_1(z)|}
\ge \frac{N(R)}R, \quad z\in\gamma_{2,4}, 
$$
where $N(R)\to\infty$ as $R\to\infty$. The~relation~(\ref{5.8}) follows 
from this bound.

Now we let $z$ run through $\gamma_2$ in the~counter clockwise direction.

We see from the~arguments at the~beginning of the~proof of the~lemma
that $\alpha-\arg w<\arg f_2(z)<\pi-\arg w$ for $z\in\Bbb C$
such that $\alpha \le \arg z<\pi$. If $z$ traverses 
$\gamma_{2,1}$, the~image $\zeta=f_3(z)$ lies in the~angular region 
$\arg w+\alpha-\pi <\arg \zeta<\arg w$ in the~$\zeta$-plane.

If $z$ traverses $\gamma_{2,2}$, by (\ref{5.1}), 
the~image $\zeta=f_3(z)$ 
lies in the~domain $|\zeta|>N_1(R),\,\arg w+\alpha-\pi
\le\arg \zeta\le \arg w+\pi-\alpha$, where $N_1(R)\to\infty$ as $R\to\infty$.

If $z$ traverses $\gamma_{2,3}$, the~image $\zeta=f_3(z)$ 
lies in the~angular domain $\arg w<\arg \zeta\le \arg w+\pi-\alpha$.

Finally, by (\ref{5.8}), when $z$ moves in $\gamma_{2,4}$,   
the~image $\zeta=f_3(z)$ lies in the~disk 
$|\zeta|<1/N_2(R)$, where $N_2(R)\to \infty$ as $R\to\infty$.

In view of these results about the image $\zeta=f_3(z)$, 
we see that the~winding number of $f_3(z)$
around $w$, when $z$ traverses $\gamma_2$  in 
the~counter clockwise direction, is equal to~one. Hence, 
by the~argument principle, 
the~function $f_3(z)$ takes the~value $w$ precisely once inside 
$\gamma_2$. Since this assertion holds for all sufficiently large $R>1$
and all $0<\alpha<\arg w$, $f_3(z)$ takes the~value $w$ precisely 
once in $\Bbb C^+$. In addition, it takes this value in 
the~domain $\arg w\le\arg z<\pi$. This proves the~lemma.
\end{Proof}

\begin{lemma}
For every $z\in\mathbb C^+$ there exist unique points $z_1\in\Bbb C^+$ and 
$z_2\in\mathbb C^+$ such that $z_1z^{-1}\in\mathbb C^+\cup(0,\infty)$  
and $z_2z^{-1}\in\mathbb C^+\cup(0,\infty)$, and 
\begin{equation}\label{5.15}
z_1z_2=zK_{\mu_1}(z_1)\quad\text{and}\quad z_1z_2=zK_{\mu_2}(z_2).
\end{equation}
\end{lemma}
\begin{Proof}
Fix $z\in\mathbb C^+$. In view of the~second relation of (\ref{5.15}), 
we have $z_1=zK_{\mu_2}(z_2)/z_2$.
Since $K_{\mu_2}(z_2)$ belongs to the~class $\mathcal K$, we obtain  
$z_1\in\mathbb C^+$ and $z_1z^{-1}\in\mathbb C^+\cup(0,\infty)$ provided that 
$z_2\in \mathbb C^+$ and $z_2z^{-1}\in\mathbb C^+\cup(0,\infty)$. 
It remains to solve the~functional equation
\begin{equation}\label{5.16}
K_{\mu_2}(z_2)=K_{\mu_1}\big(zK_{\mu_2}(z_2)/z_2\big).
\end{equation} 
Rewrite it in the~form $z=z_2/K(z_2;z)$, where 
\begin{equation}\label{5.17}
K(z_2;z)
:=K_{\mu_1}\big(zK_{\mu_2}(z_2)/z_2\big)/(zK_{\mu_2}(z_2)/z_2)
=K_{\mu_1}(z_1)/z_1.
\end{equation} 
We see that, for fixed $z\in\mathbb C^+$, the~function $\tilde f_3(z_2)
:=z_2/K(z_2;z)$ has the~same type as the~function $f_3$ in Lemma~5.1. 
Applying Lemma~5.1 to the~function $\tilde f_3$, we obtain 
the~assertion of the~lemma.
\end{Proof}

We need as well the~following auxiliary lemma which is an~analogue 
of Lemma~4.1. Denote $\mathbb S_{\pi}:=\{z\in\mathbb C:\,0<\IM z<\pi\}$.
\begin{lemma}
Assume that $f_4$ has the~representation
\begin{equation}\label{5.18}
f_4(z):=-a_1z+\frac {a_2}z+\int\limits_{(0,\infty)}\frac{1+uz}{z-u}
\,\sigma(du),\quad 0<\arg z<2\pi,
\end{equation}
where $a_j\ge 0,\,j=1,2$, and $\sigma\,(\sigma\not\equiv 0)$ 
is a~finite nonnegative measure.
Then the~function $f_5(z):=\log z+f_4(z)$, $f_5:\mathbb C^+\to\mathbb C$ 
takes every value in $\mathbb S_{\pi}$ precisely once. 
The~inverse function $f_5^{(-1)}:\mathbb S_{\pi}\to\mathbb C^+$ thus defined  
has the~property that $f_5^{(-1)}(\log z)$ belongs to 
the class $\mathcal K$.
\end{lemma}
\begin{Proof}
We shall prove that for every $w\in\mathbb S_{\pi}$ there exists 
an~unique $z\in\mathbb C^+$ such that
\begin{equation}\label{5.19}
w=f_5(z)=\log z+f_4(z).
\end{equation}
Recall that we take the~principle branch of the~logarithm only.
Let $\alpha\in(0,\IM w)$ and let $\gamma_2$ be the~curve defined 
in the~proof of Lemma~5.1.

Let $z$ traverse $\gamma_{2,3}$. The~image $\zeta=f_5(z)$ lies
on the~line $\IM \zeta=\pi$. Furthermore, we have 
$\log (-R)+f_4(-R)=\log R+f_4(-R)+i\pi$ with 
$\log R+f_4(-R)\ge (\log R)/2$, and $\log (-1/R)+f_4(-1/R)=-\log R
+f_4(-1/R)+i\pi$ with $-\log R+f_4(-1/R)\le -(\log R)/2$.

Let $z$ move in $\gamma_{2,4}$. Since
$$
\int\limits_{(0,\infty)}\frac{u|z|}{|z-u|}\,\sigma(du)\le
\sigma((0,\infty))|z|d(z),\quad z\in\mathbb C^+,
$$ 
where $d(z):=1$ for $\RE z\le 0$ and $d(z):=|z|/\IM z$, 
for $\RE z>0$, we easily conclude, for $z\in\gamma_{2,4}$, 
$\RE z\le 0$, 
$$
\RE\Big(f_4(z)+a_1z-\frac {a_2}z\Big)\le \frac cR ,\qquad R\to\infty,
$$
with some positive constant $c$, and, for $z\in\gamma_{2,4},\,\RE z>0$, 
$$
\arg\Big(f_4(z)+a_1z-\frac {a_2}z\Big)=\arg
\Big(-\int\limits_{(0,\infty)}\frac{u\sigma(du)}{|z-u|^2}+
\bar z\int\limits_{(0,\infty)}\frac{\sigma(du)}{|z-u|^2}+O(1/R)\Big)
\le 2\pi-\frac{\alpha}2.
$$
Hence the~image $\zeta=f_5(z)$ either lies in the~domain 
$\{\zeta\in\mathbb C:\alpha<\Im \zeta<\pi,\,|\zeta|\ge c(\alpha)\log R\}$, 
where $c(\alpha)$ is a~positive constant, 
or in the~half-plane $\{\zeta\in\mathbb C:\IM \zeta\le\alpha\}$.

Let $z$ traverse $\gamma_{2,1}$. Note that the~image $\zeta=f_5(z)$ 
lies in the~half-plane $\IM \zeta\le \alpha$.  

Let finally $z$ traverse $\gamma_{2,2}$. Since $\IM f_4(z)\le 0$, 
we see that the~image $\zeta=f_5(z)$ lies in the~half-plane 
$\{\zeta\in\mathbb C:\IM \zeta\le\alpha\}$ or in the~domain
$\{\zeta\in\mathbb C:\alpha<\Im \zeta<\pi\}$. 

Assume that $\zeta=f_5(z)$ lies in $\{\zeta\in\mathbb C:\alpha<\Im \zeta<\pi\}$.
In this case we easily see that $-\IM f_4(z)\le \pi-\alpha$ and
we obtain the inequality
\begin{equation}\label{5.20}
-\IM f_4(z)=a_1\IM z+a_2\frac {\IM z}{|z|^2}+\IM z\int\limits_{(0,\infty)}
\frac{1+u^2}{|z-u|^2}\,\sigma(du)\le \pi-\alpha,
\end{equation}  

On the~other hand note that, for $z\in\gamma_{2,2}$, 
\begin{equation}\label{5.21}
\RE f_4(z)=-a_1\RE z+\int\limits_{(0,\infty)}
\frac{u(|z|^2-u\RE z)}{|z-u|^2}\,\sigma(du)+O(1),\quad R\to\infty.
\end{equation}  

Since, for $z\in\gamma_{2,2}$, 
$
\big|f_4(z)+a_1z-a_2/z\big|=o(R),
$
we deduce from (\ref{5.21}) in the~case $a_1\ne 0$ the bound
$|\RE f_4(z)|\ge cR$. This means that $\zeta=f_5(z)$ lies in  
the~domain
$\{\zeta\in\mathbb C:\alpha<\Im \zeta<\pi,\,|\zeta|>\frac 12\log R\}$.
If $a_1=0$, then we have from (\ref{5.20}) and (\ref{5.21})
\begin{align}
\RE f_4(z)&\ge \int\limits_{(0,\infty)}
\frac{uR^2}{|z-u|^2}\,\sigma(du)-(1+sign (\RE z))
\frac{\RE z}{\IM z}(\pi-\alpha)+O(1)\notag\\
&\ge R^2\int\limits_{(0,\infty)}\frac{u\sigma(du)}
{|z-u|^2}+O(1),\quad R\to\infty.\notag
\end{align}
We conclude again that $\zeta=f_5(z)$ lies in the~domain
$\{\zeta\in\mathbb C:\alpha<\Im \zeta<\pi,\,|\zeta|>\frac 12\log R\}$.

Hence the~image $\zeta=f_5(z)$ for $z\in\gamma_{2,2}$
lies in the~half-plane  
$\{\zeta\in\mathbb C:\IM \zeta\le\alpha\}$ or in the~domain
$\{\zeta\in\mathbb C:\alpha<\Im \zeta<\pi,\,|\zeta|>\frac 12\log R\}$.

Therefore we conclude that the~image $\zeta=f_5(z)$ winds around
$w$ once when $z$ runs through $\gamma_2$. By the~argument principle,
there is an~unique point $z$ inside 
$\gamma_2$ such that (\ref{5.19}) holds.
This relation is valid for all sufficiently large $R>1$ and sufficiently
small $\alpha>0$. Thus for every fixed $w\in\mathbb S_{\pi}$ there is 
an~unique point $z\in\mathbb C^+$ such that (\ref{5.19}) holds. This implies
that the~inverse function $q=f_5^{(-1)}:\mathbb S_{\pi}\to\mathbb C^+$
exists and is analytic on $\mathbb S_{\pi}$. 

Let us show that $q(z)$ admits an~analytic continuation 
on the~half-line $\gamma_-:\IM z=\pi,\Re z<0$, and that its value on 
this half-line is negative.  
It is easy to~see that 
$$
f_5'(x)=\ffrac 1x -a_1-\frac {a_2}{x^2}-\int\limits_{(0,\infty)}
\frac{(1+u^2)\sigma(du)}{(x-u)^2}<0,\quad x<0.
$$ 
Since $f_5(z)$ is analytic 
on $(-\infty,0)$, we conclude that $f_5^{(-1)}$ exists and 
is analytic on $\gamma_-$ as well. Since, as shown above,
for every fixed $w\in\mathbb S_{\pi}$ there is  
an~unique point $z\in\mathbb C^+$ such that (\ref{5.19}) holds,
this function coincides  for $z\in\mathbb C^+$ with the~function
$q(z)$ obtained early. Introduce the~function 
$f_6(z):=q(\log z),\,z\in\mathbb C^+$. Note that $f_6\in\mathcal N$ and
$f_6^{(-1)}(z)=z\exp\{f_4(z)\}$, on the~domain $\mathbb C^+$, 
where $f_6^{(-1)}$ exists. Moreover, the~function $f_6(z)$ admits
an~analytic continuation on $(-\infty,0)$.

From the~definition of 
$f_6(z)$ it follows that $f_6(x)<0$ for $x<0$ and 
$f_6(x)\to 0$ as $x\uparrow 0$. By Corollary~\ref{3.3co}, 
$f_6$ belongs to the~class $\mathcal K$ and  
we obtain the~assertion of the~lemma.
\end{Proof}

{\it Proof of Theorem $2.4$.}
Consider the~function 
$$
F(z,z_2):=K_{\mu_2}(z_2)-K_{\mu_1}\big(zK_{\mu_2}(z_2)/z_2\big)
$$ 
which, by (2.8), is analytic on $\Bbb C^+\times\mathbb C^+$. 
By Lemma~5.2, for every
fixed $z=z^0\in\Bbb C^+$ the~equation~(5.16)
has an~unique solution $z_2=z_2^0\in\{z_2\in\mathbb C:\arg z^0
\le\arg z_2<\pi\}$. Let us show that the~function $F(z,z_2)$
satifies (\ref{3.8}) at the~point $(z^0,z_2^0)$.
Note that if $K_{\mu_2}(z_2)=K_{\mu_1}(z^0
K_{\mu_2}(z_2)/z_2)$ holds for $z_2\in\mathbb C^+$ such that
$|z_2-z_2^0|<r'$ with $r'>0$, then this equality holds  for all
$z_2\in\Bbb C^+$. By (5.2), the~function $K(-r;z^0)$ introduced
in (\ref{5.17}) has the~property $K(-r;z^0)/r\to 0$ as
$r\to \infty$. Since the~relation $z^0=-r/R(-r;z^0)$ holds for $r>0$, 
we arrive at a~contradiction. Hence the~point $(z^0,z_2^0)$
satisfies the~assumptions (\ref{3.8}) of Theorem~3.7. Now repeating almost
word for word the~arguments of the~proof of Theorem~2.1, using Lemma~5.2
instead of Lemma~4.2, and Theorem~\ref{3.1th}, Corollary~\ref{3.3co},
we obtain the~assertion of the~theorem. 
$\square$

We prove Corollaries~2.5 and 2.6 in the~same way as Corollaries~2.2 
and 2.3. 

Rewrite the~relation (\ref{2.12}) with $n-1=t\in[0,\infty)$
in the~form
\begin{equation}\label{5.21a}
Z(z)=z\Big(\frac{K_{\mu}(Z(z))}{Z(z)}\Big)^t,\quad z\in\mathbb C^+,\quad t\ge 0.
\end{equation}
Consider the~function $F(z):=K_{\mu}(z)/z,\,z\in\mathbb C^+$. This function
belongs to the~class $\mathcal N$, is analytic and nonnegative on 
the~negative real axis $(-\infty,0)$. By Theorem~\ref{3.3th}, $F(z)$
admits the~representation (\ref{5.18}) with $a_1=a_2=0$. By Lemma~5.3, 
there exists a~function $Z(z)\in\mathcal K$ such that (\ref{5.21a}) holds.

The~relation (\ref{5.21a}) implies that there exists a~semigroup
$\mu_t\in\mathcal M_+,\,t\ge 1$, such that $(\Sigma_{\mu}(z))^t=\Sigma_{\mu_t}(z),
\,t\ge 1$.

\vskip 0,5cm

2. Consider the~case of multiplicative convolution for p-measures 
of the~class $\mathcal M_*$. In order to prove Theorem~1.7
we need some auxiliary results.

\begin{lemma}
Let $Q\in \mathcal S$ and $Q(0)\ne 0$. Then the~function $g:\mathbb D\setminus\{0\}
\to\mathbb C:z\mapsto Q(z)/z$ 
takes every value in $\Bbb C\setminus\overline{\Bbb D}$ precisely once.
\end{lemma}
\begin{Proof} 
Without loss of generality we assume that $Q(0)>0$. 
Let $w\in\Bbb C\setminus\overline{\Bbb D}$ and let $r\in(0,1)$. We assume
that $1-r$ is sufficiently small. It is easy to see that
there exist $\theta=\theta(w,r)\in[0,2\pi)$ and $\theta\ne \arg w$ 
such that $Q(z)/z\ne w$ for $z=\rho e^{i\theta},\,1-r\le\rho\le r$. 
Let $\gamma_3=\gamma_3(w,\theta)$ 
(see Figure~3) be the~closed rectifiable curve consisting of 
the~two circles $\gamma_{3,2}:|z|=1-r$ and $\gamma_{3,4}:|z|=r$, and
the~segments $\gamma_{3,1},\gamma_{3,3}$ which coincide with
the~segment $\gamma:te^{i\theta},1-r\le t\le r$, traversed twice 
in opposite directions.

\begin{figure}[htbp]
 
\input{bild3.pstex_t}
 
\caption{}
\end{figure}

We note that 
the~number of loops around $w$ by the~map $z\mapsto Q(z)/z$ as 
$z$ traverses the~curve $\gamma_3$ once
in the~counter clockwise direction is equal to 1. Indeed, if $z$
traverses $\gamma_{3,1}$  the~number of windings around $w$
is equal to some number $\alpha$. If $z$ traverses $\gamma_{3,4}$ the~number of
windings around $w$ is equal to $0$. Finally, moving along 
$\gamma_{3,3}$, the~number of windings around $w$ is equal to $-\alpha$.
Hence traversing $\gamma_{3,1},\,\gamma_{3,4}$ and $\gamma_{3,3}$
in the~counter clockwise direction the~number of windings around $w$ is equal 
to $0$. Traversing $\gamma_{3,2}$, we wind around $w$ once.   
Hence, by the~argument principle,
the~function $Q(z)/z$ takes the~value $w$ in the~domain $\{z:1-r<|z|
<r\}$ precisely once. Since this assertion holds for all
$0<r<1$, $Q(z)/z$ takes every value $w\in \mathbb C\setminus\overline{\mathbb D}$ 
in $\mathbb D$ precisely once.
\end{Proof}

Let $\mu_1$ and $\mu_2$ belong to the~class $\mathcal M_*$.
Denote $a_j:=\int_{\mathbb T}\xi\,\mu_j(d\xi),\,j=1,2$.
Recall that $Q_{\mu_j}(0)=0,\,a_j\ne 0,\,j=1,2$.
\begin{lemma}
For every $z\in\mathbb D\setminus\{0\}$, there exist unique points 
$z_1\in\mathbb D$ and $z_2\in\Bbb D$ such that $(z_1,z_2)\ne(0,0)$
and
\begin{equation}\label{5.22}
z_2=z\frac{Q_{\mu_1}(z_1)}{z_1}\quad\text{and}\quad 
z_1=z\frac{Q_{\mu_2}(z_2)}{z_2}.
\end{equation}
\end{lemma}
\begin{Proof} Fix $z\in\mathbb D\setminus\{0\}$. 
Assume for definiteness $|a_2|\le |a_1|$.
It follows from Proposition~\ref{3.5pro} that
$Q_{\mu_1}(a_2z)\ne 0$ for every $z\in\mathbb D$. 

By Schwarz's lemma,
$\big| Q_{\mu_2}(z_2)/z_2\big|\le 1$, $z_2\in\mathbb D$.
The~second relation of (\ref{5.22}) implies that if $z_2\in\Bbb D$,
then $z_1\in\Bbb D$ as well. Hence, by the~first relation of 
(\ref{5.22}), we need to solve the~functional equation 
\begin{equation}\label{5.23}
Q_{\mu_2}(z_2)=Q_{\mu_1}\big(zQ_{\mu_2}(z_2)/z_2\big),\quad z_2\in\mathbb D.
\end{equation}
Rewrite this relation in the~form $z_2/z=Q(z_2)$, where 
$$
Q(z_2)=Q_{\mu_1}\big(zQ_{\mu_2}(z_2)/z_2\big)/(zQ_{\mu_2}(z_2)/z_2).
$$
The~function $Q$ satisfies $Q(0)=Q_{\mu_1}(a_2z)/(a_2z)\ne 0$ 
and $Q\in \Cal S$.
By Lemma~5.4, we obtain the~assertion of the~lemma for all 
$z\in\mathbb D\setminus\{0\}$.
\end{Proof}

We need as well the~following auxiliary lemma which is an~analogue 
of Lemma~4.1 for the~semigroup $(\mathcal M_*,\boxtimes)$.
\begin{lemma}
Define
\begin{equation}
Q_1(z):=\int\limits_{\mathbb T}\frac{\xi+z}{\xi-z}
\,\sigma(du),\quad z\in\mathbb D,
\end{equation}
where $\sigma\,(\sigma\not\equiv 0)$ is finite nonnegative measure.
Then the~function $Q_2:\mathbb D\to\mathbb C$, $Q_2(z):=z\exp\{Q_1(z)\}$, 
takes every value in $\mathbb D$ precisely once. 
The~inverse $Q_2^{(-1)}:\mathbb D\to \mathbb D$
thus defined is in the class $\mathcal S_*$.
\end{lemma}
\begin{Proof}
Let $w\in\mathbb D$. Introduce the~curve $\gamma_4:\,|z|=r$, where $r\in(0,1)$
and $1-r$ is sufficiently small. Note that $\Im Q_1(0)=0$. Then, by 
the~minimum and maximum principles for harmonic functions, 
for every $r\in(0,1)$ there exist points $\theta'$ and 
$\theta''\in[0,2\pi)$ such that 
$\Im Q_1(re^{i\theta'})<0$ and $\Im Q_1(re^{i\theta''})>0$. 
Since for every fixed $r\in(0,1)$ 
$\theta\mapsto \IM Q_1(re^{i\theta})$ is analytic on the~segment $[0,2\pi]$, 
there are only a~finite number of points $\{re^{i\theta_l}\}_{l=1}^k$
with $0\le\theta_1<\theta_2<\dots<\theta_k<2\pi$,
such that $\IM Q_1(re^{i\theta_l})=0,\,l=1,2,\dots,k$. 
Consider the~arcs $\gamma_{4,l}:\theta_l<\arg z\le \theta_{l+1},\,|z|=r$,  
$l=1,\dots,k$, with $\theta_{k+1}=\theta_1+2\pi$.

\enlargethispage{1cm}

\begin{figure}[htbp]
 
\input{bild4.pstex_t}
 
\caption{}
\end{figure}

Let $z$ run transverse $\gamma_{4,l}$ in the~counter clockwise direction. 
Since $|\exp\{Q_1(z)\}|\ge 1$, we note that the~change in 
$Arg (Q_2(z)-w)$ is equals to 
$$
(\arg(re^{i\theta_{l+1}}-w)-\arg(re^{i\theta_l}-w))/(2\pi).
$$ 
Hence the~image $\zeta=Q_2(z)$ winds around $w$ once when $z$ transverses
$\gamma_4$ once in the~counter clockwise direction. By the~argument principle,
the~function $Q_2:\Bbb D\to\Bbb C$ takes the~value $w\in\mathbb D$ in 
$\mathbb D$ precisely once. The~inverse function $Q_2^{(-1)}:\mathbb D\to\mathbb D$ 
thus defined is analytic on $\mathbb D$ and obviously belongs to the~class 
$\mathcal S_*$. The~lemma is proved. 
\end{Proof}

\medskip
{\it Proof of Theorem $2.7$.}
We assume for definiteness that $|a_2|\le |a_1|$.
Consider the~function $F(z,z_2):=Q_{\mu_2}(z_2)-
Q_{\mu_1}\big(zQ_{\mu_2}(z_2)/z_2\big)$ which is analytic on
$\mathbb D\times\mathbb D$. For a~fixed $z=z^0\in\mathbb D$ 
the~equation~(\ref{5.23}) has, by Lemma~5.5, an~unique solution 
$z_2=z_2^0\in\mathbb D$. 
Moreover, if $F(z^0,z_2)=0$ for $z_2$ in some neighborhood of 
the~point $z_2^0$,
then $F(z^0,z_2)=0$ for all $z_2\in\mathbb D$. This relation is equivalent
to $z_2/z^0=Q(z_2),\,z_2\in\Bbb D$, where $Q\in\mathcal S$. 
Choosing here $|z_2|>|z^0|$ we arrive at a~contradiction. 
Hence the~function $F(z,z_2)$ satisfies the~assumptions (\ref{3.8}) 
of Theorem~3.7 at the~point $(z^0,z_2^0)$. Repeating almost
word for word the~arguments of the~proof of Theorem~2.1, using Lemma~5.5
instead of Lemma~4.2, we obtain that there exists a~neighborhood
$|z-z^0|<r(z^0)$ such that the~equation $F(z,z_2)=0$ has an~unique
regular solution $z_2=z_2(z;z^0)$ with values in $\mathbb D$. 
Choosing $z_1(z;z^0):=zQ_{\mu_2}(z_2(z;z^0))/z_2(z;z^0)$, there are 
two regular functions $z_1(z;z^0)$ and $z_2(z;z^0)$ with values 
in $\mathbb D$ which satisfy the~equations
$$
z_1(z;z^0)z_2(z;z^0)=zQ_{\mu_1}(z_1(z;z^0))
\quad\text{and}\quad z_1(z;z^0)z_2(z;z^0)=zQ_{\mu_2}(z_2(z;z^0)).
$$  
By Lemma~5.5, for the~points $z'\in\mathbb D$ and $z''\in\mathbb D,\,z'\ne z''$,
we obtain the~identities $z_1(z;z')=z_1(z;z''),\,z_2(z;z')=z_2(z;z'')$ for
$z$ of the~domain $\{|z-z'|<r(z')\}\cap\{|z-z''|<r(z'')\}$. 
By the~monodromy theorem, there exist regular functions $Z_1(z)$ 
and $Z_2(z),\,z\in\Bbb D$, such that, for every point $z^0\in\mathbb D$,
$Z_j(z)=z_j(z;z^0),\,j=1,2$, for $|z-z^0|<r(z^0)$. These functions
belong to the~class $\mathcal S$ and are unique solutions of 
equations~(\ref{2.14}) for $\mu_j,\,j=1,2$.
Moreover, it is easy to see that $Z_1,Z_2\in\mathcal S_*$. This proves
the~theorem.
$\square$

We prove Corollaries~2.8 and 2.9 in the~same way as Corollaries~2.2 
and 2.3. 

Now we prove the~relation (\ref{2.16a}) with $n=t\in[1,\infty)$.
In the ~first step we note from (\ref{2.16a}) that $Q_{\mu\boxtimes\mu}(z)/z
\ne 0,\,z\in\mathbb D$. Since this function belongs to the~class $\mathcal S$,
it has the~form
$$
\frac{Q_{\mu\boxtimes\mu}(z)}z=\exp\Big\{ai-\int\limits_{\mathbb T}\frac{\zeta+z}
{\zeta-z}\,\sigma(d\zeta)\Big\},\quad z\in\mathbb D,
$$
where $a\in\mathbb R$ and $\sigma$ is a~finite nonnegative measure.
By Lemma~5.6, we conclude that there exists $Z(z)\in\mathcal S_*$ such that
\begin{equation}\label{5.24}
Z(z)=z \Big(\frac{Q_{\mu\boxtimes\mu}(Z(z))}{Z(z)}\Big)^{t},\quad z\in\mathbb D.
\quad t\ge 0,
\end{equation}
We obtain from (\ref{5.24}) the~existence a~semigroup $\mu_u\in\mathcal M_*,\,
u\ge 2$, such that 
\begin{equation}\label{5.25}
\big(\Sigma_{\mu}(z)\big)^u=\Sigma_{\mu_u}(z),\quad u\ge 2.
\end{equation}
If $Q_{\mu}(z)/z\ne 0,\,z\in\mathbb D$, then, as before, the~following relation
holds
\begin{equation}\label{5.26}
Z(z)=z \Big(\frac{Q_{\mu}(Z(z))}{Z(z)}\Big)^{t},\quad z\in\mathbb D,
\quad t\ge 0,
\end{equation}
and (\ref{5.25}) is true for $u\ge 1$.

\section{
Basic arithmetic 
and Khintchine's limit theorems 
for the~semigroups $(\bold M,\circ)$}

In this section we shall prove Therems~2.10, 2.11 and 2.12.

Introduce sets $\mathbf M^{\gamma}$ of measures of $\mathbf M$
in the~following way.

Consider the~semigroup $(\mathcal M,\boxplus)$.
Let $\overline{\Gamma}_{\alpha,\beta}$ denote the~closure of 
$\Gamma_{\alpha,\beta}$. Denote by $\mathcal M^{(\alpha,\beta)}$
the~set of those $\mu\in\mathcal M$ such that
$F_{\mu}(z)$ is univalent on $\Omega_{\mu}$, where
$\Omega_{\mu}\subseteq\mathbb C^+$ is a domain, for which $F_{\mu}(\Omega_{\mu})
\supset \overline{\Gamma}_{\alpha,\beta}$. The~function $F_{\mu}(z)$
has an~inverse $F_{\mu}^{(-1)}(z)$ defined on $\overline{\Gamma}
_{\alpha,\beta}$ such that $\IM \varphi_{\mu}(z)\le 0$.   
We denote $\mathbf M^{\gamma}:=\mathcal M^{(\alpha,\beta)}$, where $\gamma
:=(\alpha,\beta)\in\mathcal G:=\{(\alpha,\beta):\alpha>0,\beta>0\}$, 
a~subset of the~semigroup $(\mathcal M,\boxplus)$.

Consider $(\mathcal M_+,\boxtimes)$ and let $\overline{\Gamma}
_{\alpha,\beta,\Delta}^{\,+}$ denote the~closure of 
$\Gamma_{\alpha,\beta,\Delta}^+$. 
Denote by $\mathcal M^{(\alpha,\beta,\Delta)}_+$ a~subset of  
the~set of those $\mu\in\mathcal M_+$ such that $K_{\mu}(z)$ is univalent on 
$\Omega_{\mu}$, where $\Omega_{\mu}\subseteq\mathbb C\setminus [0,+\infty)$
is a~domain, with $K_{\mu}(\Omega_{\mu})\supset 
\overline{\Gamma}_{\alpha,\beta,\Delta}^{\,+}$. The~function $K_{\mu}(z)$
has an~inverse $K_{\mu}^{(-1)}(z)$ defined on $\overline{\Gamma}
_{\alpha,\beta,\Delta}^{\,+}$ such that $\arg \Sigma_{\mu}(z)\le 0$ for $z\in
\overline{\Gamma}_{\alpha,\beta,\Delta}^{\,+}\cap\mathbb C^+$ 
and $\arg \Sigma_{\mu}(z)\ge 0$ for $z\in\overline{\Gamma}
_{\alpha,\beta}^{\,+}\cap\mathbb C^-$. 
As above we write $\mathbf M^{\gamma}
:=\mathcal M^{(\alpha,\beta,\Delta)}_+$, where 
$\gamma:=(\alpha,\beta,\Delta)\in\mathcal G:=\{(\alpha,\beta,\Delta):
\alpha\in(0,\pi),0<\beta<\Delta<\infty\}$,
for a~subset of the~semigroup $(\mathcal M_+,\boxtimes)$.

Consider $(\mathcal M_*,\boxtimes)$ and let $\overline{\mathbb D}_{\alpha}$ denote
the~closure of $\mathbb D_{\alpha}:=\{z\in\mathbb C:|z|<\alpha\}$. 
Denote by $\mathcal M^{\alpha}_*$
the~set of those $\mu\in\mathcal M_*$ such that $Q_{\mu}(z)$ is univalent on 
$\Omega_{\mu}$, where $\Omega_{\mu}\subseteq\mathbb D$ is a~domain,
with $Q_{\mu}(\Omega_{\mu})\supset\overline{\mathbb D}_{\alpha}$. 
The~function $Q_{\mu}(z)$ has an~inverse $Q_{\mu}^{(-1)}(z)$ 
defined on $\overline{\mathbb D}_{\alpha}$ such that 
$|\Sigma_{\mu}(z)|\ge 1$. Introduce $\mathbf M^{\gamma}:=\mathcal M^{\alpha}_*$, 
where $\gamma:=\alpha\in\mathcal G:=\{\alpha:0<\alpha<1\}$, for a~subset of 
the~semigroup $(\mathcal M_*,\boxtimes)$. Note that $\mathbf M=
\cup_{\gamma\in\mathcal G }\mathbf M^{\gamma}$. 

A~next proposition show that the~set $\mathbf M^{\gamma}$ is hereditary 
in $(\mathbf M,\circ)$.
\begin{proposition}\label{6.1pr}
Given $\gamma\in\mathcal G$, let $\mu\in\mathbf M^{\gamma}$ and
$\mu=\mu_1\circ\mu_2$, where $\mu_j\in\mathbf M,\,j=1,2$. Then
$\mu_j\in\mathbf M^{\gamma},\,j=1,2$.  
\end{proposition}
\begin{Proof} 
At first we shall prove the~assertion of the~proposition for 
$(\mathcal M,\boxplus)$.
Let $\mu\in\mathcal M^{(\alpha,\beta)}$ and $\mu=\mu_1\boxplus\mu_2$, where
$\mu_1,\mu_2\in\mathcal M$. We should show that
$\mu_1,\mu_2\in\mathcal M^{(\alpha,\beta)}$. By Theorem~2.1,
there exist functions $Z_1,Z_2\in\mathcal F$ such that
\begin{equation}\label{6.1}
F_{\mu}(z)=F_{\mu_1}(Z_1(z))=F_{\mu_2}(Z_2(z)),
\quad z\in\Bbb C^+.
\end{equation}
By the~definition of $\mathcal M^{(\alpha,\beta)}$, there exists
a~domain $\Omega_{\mu}$ such that $F_{\mu_j}(Z_j(z))$, $j=1,2$, 
are univalent on $\Omega_{\mu}$ and 
$F_{\mu_j}(Z_j(\Omega_{\mu}))\supset 
\overline{\Gamma}_{\alpha,\beta}$, respectively.  
Consider the~domains $\Omega_{\mu_j}:=Z_j(\Omega_{\mu}),\,j=1,2$.
It is clear that $Z_j(z),\,j=1,2$, are univalent on $\Omega_{\mu}$
and $F_{\mu_j}(z),\,j=1,2$, are univalent on
$\Omega_{\mu_j},\,j=1,2$, respectively,
thus proving the~assertion.

Consider $(\mathcal M_+,\boxtimes)$. Let $\mu\in\mathcal M_+^{(\alpha,\beta,\Delta)}$
and $\mu=\mu_1\boxtimes\mu_2$, where $\mu_j\in\mathcal M_+,\,j=1,2$. 
Let us show that $\mu_1,\mu_2\in\mathcal M_+^{(\alpha,\beta,\Delta)}$. 
By Theorem~1.4, there exist functions $Z_1,Z_2\in\Cal K$ such that
\begin{equation}\label{6.2}
K_{\mu}(z)=K_{\mu_1}(Z_1(z))=K_{\mu_2}(Z_2(z)),
\quad z\in\Bbb C^+.
\end{equation}
By the~definition of $\mathcal M^{(\alpha,\beta,\Delta)}$, there exists
a~domain $\Omega_{\mu}$ such that $K_{\mu_j}(Z_j(z))$
are univalent on $\Omega_{\mu}$ and  
$K_{\mu_j}(Z_j(\Omega_{\mu}))\supset 
\overline{\Gamma}_{\alpha,\beta,\Delta},\,j=1,2$, respectively. 
Denote $\Omega_{\mu_j}:=Z_j(\Omega_{\mu}),\,j=1,2$. 
We see that $Z_j(z),\,j=1,2$, are univalent on $\Omega_{\mu}$
and $K_{\mu_j}(z),\,j=1,2$, are univalent on
$\Omega_{\mu_j}$, respectively,
thus proving the~assertion.

Repeating the~previous arguments, using Theorem~2.7 instead 
of Theorem~2.4, we obtain the~assertion
of the~proposition for $\mathcal M^{\alpha}_*$ as well
which completes the~proof of Proposition~\ref{6.1pr}.
\end{Proof}

Define an~equivalence relation $E$ on $\mathbf M$ via $\mu E\nu$
if $\mu=\nu\circ \delta_a$, where $\delta_a$ is Dirac measure with 
the~support at the~point $a$. Here $a$ denotes a~real number for
$(\mathcal M,\boxplus)$, is a~positive number for
$(\mathcal M_+,\boxtimes)$ and $a\in\mathbb T$ for 
$(\mathcal M_*,\boxtimes)$.
It is clear that $E$ is an~equivalence relation on the~set 
$\mathbf M^{\gamma}$ as well. We shall denote the~equivalence class
of $\mu\in\mathbf M$ as $E\mu$, and the~element of $\mathbf M/E$ by
$\mu^*,\nu^*$ etc. Note that $\mathbf M/E$ is a~semigroup under 
$\circ$ multiplication of representatives.
This semigroup has an~identity element $e:=\delta_a^*$,
where $a=0,1,1$ for the~semigroups $(\mathcal M,\boxplus)$,
$\mathcal M_+,\boxtimes)$, $(\mathcal M_*,\boxtimes)$, respectively.
In the~semigroup $(\mathbf M/E,\circ)$ an~element $\mu^*$ is indecoposable
if $\mu^*=E\mu$, where $\mu$ is indecomposable in $(\mathbf M,\circ)$.


For $\mu^*,\nu^*\in \mathbf M/E$ define a~distance $L^*(\mu^*,\nu^*)$
via
$$
L^*(\mu^*,\nu^*)=inf\{L(\mu,\nu):E\mu=\mu^*,\,E\nu=\nu^*\},
$$
where $L(\mu,\nu)$ denotes the~L\'evy metric. It is easy to check that
$L^*(\mu^*,\nu^*)$ is a~metric, and that $\{\mu_n^*\}$ converges to
$\mu^*\in\mathbf M/E$ with respect to $L^*$ if and only if there is 
$\{\nu_n\}$ with
$\nu_n\in\mathbf M$ and $\mu_n^*=E\nu_n$ for each $n$, and $\nu\in
\mathbf M$ with $\mu^*=E\nu$, such that $\{\nu_n\}$ converges to $\nu$ 
with respect to $L$.

By $M$ we denote the~set of sequences of elements of $\mathbf M$ 
with weak limits in $\mathbf M$, and $l$ denotes the~operation of taking 
the~weak limit. We denote by $M/E$ the~set of sequences of elements of
$\mathbf M/E$ that admit limits with respect to $L^*$ and by $l$ again 
the~operation of taking the~limit.

Define a~homomorphism $D^{\gamma}$ of $\mathbf M^{\gamma}/E$
to non-negative real numbers with addition in the~following way.

For the~semigroup $(\mathcal M,\boxplus)$:

$D^{\gamma}(\mu^*):=-\IM\varphi_{\nu}(i(\beta+1))$,
where $\nu=E\mu$ and $\nu\in\mathcal M^{(\alpha,\beta)}$. 

For the~semigroup $(\mathcal M_+,\boxtimes)$:

$D^{\gamma}(\mu^*):=-\IM \log\Sigma_{\nu}((\beta+\Delta)
e^{i(\pi+\alpha)/2}/2)$, 
where $\nu=E\mu$ and $\nu\in\mathcal M^{(\alpha,\beta,\Delta)}_+$. 

For the~semigroup $(\mathcal M_*,\boxtimes)$:

$D^{\gamma}(\mu^*):=\RE \log\Sigma_{\nu}(\alpha/2)$, where
$\nu=E\mu$ and $\nu\in\mathcal M^{\alpha}_*$. 

Note that $D^{\gamma}$ is finite-valued and
does not depend on the~particular $\nu$.

\begin{proposition}\label{6.2apr}
Let $\mu\in\mathbf M^{\gamma}$.
If the~relation $D^{\gamma}(\mu^*)=0$ holds, then $\mu^*=e$.
\end{proposition}
\begin{proof}
We shall prove the~proposition for the~semigroup $(\mathcal M,\boxplus)$.
Let $D^{(\alpha,\beta)}(\mu^*)=0$. Then $=-\IM\varphi_{\mu}(i(\beta+1))=0$, 
where $\mu^*=E\mu$. Since the~function
$-\IM\varphi_{\mu}(z)$ is harmonic and nonnegative in 
$\Gamma_{\alpha,\beta}$ and the~point $i(\beta+1)$ lies in
$\Gamma_{\alpha,\beta}$, we conclude, by minimum principle for
harmonic functions, that $\IM\varphi_{\mu}(z)=0$ as 
$z\in \Gamma_{\alpha,\beta}$. From this it follows that 
there exists a~real number $a$ such that $\varphi_{\mu}(z)=a$ 
as $z\in \Gamma_{\alpha,\beta}$ and we have $F_{\mu}(z)=z-a,\,
z\in\Bbb C^+$. Hence $\mu=\delta_a$ and the~assertion is proved.

For the other semigroups $(\mathcal M_+,\boxtimes)$
and $(\mathcal M_*,\boxtimes)$ the~proof of the~proposition is similar.
\end{proof}

\begin{proposition}\label{6.2bpr}
Let $\mu$ and $\mu_n,\,n=1,\dots$, belong to the~set $\mathbf M^{\gamma}$
for a~fixed $\gamma\in \mathcal G$. Let $\{\mu_n^*\}\in M/E$ and  
$l\{\mu_n^*\}=\mu^*$. Then $D^{\gamma}(\mu^*)=\lim_{n\to\infty}
D^{\gamma}(\mu^*_n)$. 
\end{proposition}
We omit the~proof of this simple proposition.

\begin{proposition}\label{6.2pr}
The~operation $\circ$ is continuous with respect to L\'evy's metric.
\end{proposition}

This result is due to Bercovici and Voiculescu~\cite{BercVo:1992}, 
\cite{BercVo:1993}.
For convenience of the~reader we include a~proof. 

{\it Proof of Proposition~$\ref{6.2pr}$.} 
Let $\{\mu_n\}$ and $\{\nu_n\}$ be in $M$. Let us prove that
$\{\mu_n\boxplus\nu_n\}\in M$ and $l\{\mu_n\boxplus\nu_n\}
=l\{\mu_n\}\boxplus l\{\nu_n\}$.

At first consider the~semigroup $(\Cal M,\boxplus)$.
By Proposition~\ref{3.11pro}, the~sequences $\{\varphi_{\mu_n}\}$ and
$\{\varphi_{\nu_n}\}$ converge
uniformly on compact subsets of some $\Gamma_{\alpha,\beta}$ to 
the~functions $\varphi_{\mu}$ and  $\varphi_{\nu}$, 
and $\varphi_{\mu_n}(iy)=o(y)$, $\varphi_{\nu_n}(iy)=o(y)$
uniformly in $n$ as $y\to +\infty$. Since $\varphi_{\mu_n\boxplus
\nu_n}(z)=\varphi_{\mu_n}(z)+\varphi_{\nu_n}(z)$ and  
$\varphi_{\mu\boxplus\nu}(z)=\varphi_{\mu}(z)+\varphi_{\nu}(z)$
for $z\in\Gamma_{\alpha,\beta}$, note that $\{\varphi_{\mu_n\boxplus
\nu_n}\}$ converges uniformly on compact subsets of some 
$\Gamma_{\alpha,\beta}$ to the~function $\varphi_{\mu\boxplus\nu}$.
Furthemore, $\varphi_{\mu_n\boxplus\nu_n}(iy)=o(y)$ uniformly in $n$ as 
$y\to +\infty$. The~assertion for $(\Cal M,\boxplus)$ now
follows immediately from Proposition~\ref{3.11pro}.

Consider the~semigroup $(\mathcal M_+,\boxtimes)$.
Assume that the~sequences $\{\mu_n\}$ and $\{\nu_n\}$ converge 
in the~weak topology
to measures $\mu\ne\delta_0$ and $\nu\ne\delta_0$, respectively. 
Since $\Sigma_{\mu_n\boxtimes\nu_n}=\Sigma_{\mu_n}\Sigma_{\nu_n}$ 
on the~set, where these functions 
are defined, we easily conclude from Proposition~\ref{3.12pro} that
the~sequence $\{\mu_n\boxtimes\nu_n\}$ converges as well to a~p-measure
which is not the~Dirac measure $\delta_0$. Indeed, by 
Proposition~\ref{3.12pro}, there exist numbers $\alpha\in(0,\pi)$ and
$0<\beta<\Delta$ such that
the~sequences $\{\Sigma_{\mu_n}\}$ and $\{\Sigma_{\nu_n}\}$ converge 
uniformly on $\Gamma_{\alpha,\beta,\Delta}^+$ to the~functions $\Sigma_{\mu}$
and $\Sigma_{\nu}$, respectively. Hence, $\{\Sigma_{\mu_n}\Sigma_{\nu_n}\}$
converges uniformly on $\Gamma_{\alpha,\beta,\Delta}^+$ 
to the~function 
$\Sigma_{\mu}\Sigma_{\nu}$ and, by Proposition~\ref{3.12pro}, we conclude that
$\{\mu_n\boxtimes\nu_n\}$ converges weakly to $\mu\boxtimes\nu\ne\delta_0$,
which is the~desired result.

It remains to consider the~semigroup $(\mathcal M_*,\boxtimes)$.
By Proposition~\ref{3.13pro}, the~sequences $\{\Sigma_{\mu_n}\}$ and 
$\{\Sigma_{\nu_n}\}$ converge uniformly in some neighborhood
of zero, say $\mathbb D_{\alpha},\,\alpha\in(0,1]$, to the~functions 
$\Sigma_{\mu}$ and $\Sigma_{\nu}$, 
respectively. Since $\Sigma_{\mu_n\boxtimes\nu_n}=\Sigma_{\mu_n}
\Sigma_{\nu_n}$ in $\mathbb D_{\alpha}$, we conclude that 
$\{\Sigma_{\mu_n\boxtimes\nu_n}\}$ converge uniformly to 
$\Sigma_{\mu}\Sigma_{\nu}$. Applying Proposition~\ref{3.13pro}, we arrive at
the~assertion of the~proposition for $(\mathcal M_*,\boxtimes)$.
$\square$

In the~sequel we denote 
$$
a(\mu):=-\RE\varphi_{\mu}(i(\beta+1)) \quad\text{for} 
\quad \mu\in\mathcal M^{(\alpha,\beta)},
$$ 
$$
a(\mu):=1/|\Sigma_{\mu}((\beta+\Delta)e^{i(\pi+\alpha)/2}/2)|
\quad\text{ for}\quad\mu\in\mathcal M_+^{(\alpha,\beta,\Delta)},
$$ 
and 
$$
a(\mu):=\exp\{-i\arg\Sigma_{\mu}(\alpha/2)\}\quad\text{ for}
\quad \mu\in\mathcal M_*^{\alpha}.
$$ 

\begin{proposition}\label{6.3pr}
If $\{\mu_n^*\}$ is any element of $M/E$, and if, for each
$n\in\mathbb N$, $\nu_n^*$ is a~factor of $\mu_n^*$, then 
there is a~subsequence $\{n'\}$ of $\{n\}$ such that 
$\{\nu_{n'}^*\}\in M/E$ and $l\{\nu_{n'}^*\}$ is a~factor of
$l\{\mu_n^*\}$.
\end{proposition}
\begin{Proof}
Let $\{\mu_n^*\}$ converges to $\mu^*\in \mathbf M/E$ in the~metric $L^*$, 
and let $\nu_n^*$ be a~factor of $\mu_n^*$ for each $n$. Then for each $n$
there exist probability measures $\mu_n$ and $\nu_n$, such that
$\mu_n^*=E\mu_n$ and $\nu_n^*=E\nu_n$, $\nu_n$ is a $\mathbf M$-factor
of $\mu_n$, and $\{\mu_n\}$ converges to $\mu$ in the~metric
$L$, where $\mu\in\mathbf M$ and $\mu^*=E\mu$.

At first let us prove the~assertion of the~proposition
for $(\Cal M,\boxplus)$.
By Proposition~\ref{3.11pro}, the~sequence $\{\varphi_{\mu_n}\}$ converges
uniformly on compact subsets of some set $\Gamma_{\alpha,\beta}$ to 
the~function $\varphi_{\mu}$, and $\varphi_{\mu_n}(iy)=o(y)$
uniformly in $n$ as $y\to +\infty$.
Since, for $\mu\in \mathcal M^{(\alpha,\beta)}$, $-\Im\varphi_{\mu}(z)\ge 0$, 
$z\in\Gamma_{\alpha,\beta}$, we have $-\Im\varphi_{\nu_n}(z)\le
-\Im\varphi_{\mu_n}(z),\,z\in\Gamma_{\alpha,\beta}$. By Theorem~\ref{3.8th},
there exists a~subsequence $\{n'\}$ of $\{n\}$ such that
$\{-\Im\varphi_{\nu_{n'}}(z)\}$ converges uniformly on compact 
subsets of $\Gamma_{\alpha,\beta}$ to a~harmonic function $U(z)\ge 0$.
Moreover $-\Im\varphi_{\nu_{n'}}(iy)=o(y)$ uniformly in $n$ 
as $y\to+\infty$.
Consider the~sequence of analytic functions  
$\{\varphi_{\nu_{n'}}(z)+a(\nu_{n'})\}
=\{\varphi_{\nu_{n'}\boxplus\delta_{a(\nu_{n'})}}(z)\}$ 
in $\Gamma_{\alpha,\beta}$. 
This sequence converges at the~point $i(\beta+1)$ and
hence, by Theorem~\ref{3.9th}, $\{\varphi_{\nu_{n'}}(z)+a(\nu_{n'})\}$
converges uniformly on every compact subset of $\Gamma_{\alpha,\beta}$.  
Hence the~assumptions of Proposition~\ref{3.11pro} hold for the~sequence 
$\{\nu_{n'}\boxplus\delta_{a(\nu_{n'})}\}$. By this proposition, 
this sequence of p-measures converges in the~weak topology 
to a~p-measure $\nu$. 
Ite\-ra\-ting this argument, we see that 
$\{\rho_{n'}\}$ converges weakly to $\rho$, say, 
where $\rho_{n'}$ is some
cofactor of $\nu_{n'}\boxplus\delta_{a(\nu_{n'})}$ in $\mu_{n'}$. By 
continuity of the~additive free convolution (Proposition~\ref{6.2pr}), 
we have $\nu\boxplus\rho=\mu$, so $\nu$ is
a~factor of $\mu$. Thus Proposition~\ref{6.3pr} 
is proved for $(\mathcal M,\boxplus)$.

Now we shall prove the~assertion for 
$(\mathcal M_+,\boxtimes)$. By Proposition~\ref{3.12pro}, the~sequence 
$\{\Sigma_{\mu_n}\}$ converges uniformly on compact subsets 
of some $\Gamma_{\alpha,\beta,\Delta}^+$ to the~function $\Sigma_{\mu}$.
Recall that for every $\mu\in\mathcal M^{(\alpha,\beta,\Delta)}_+$
the~function $\log\Sigma_{\mu}(z)$ is analytic in 
$\Gamma_{\alpha,\beta,\Delta}^+$ and, in addition, 
$-\Im\log\Sigma_{\mu}(z)\ge 0$ for
$z\in\Gamma_{\alpha,\beta,\Delta}^+\setminus\Bbb C^-$ and 
$\Im\log\Sigma_{\mu}(\bar z)=-\Im\log\Sigma_{\mu}(z)$ for
$z\in\Gamma_{\alpha,\beta,\Delta}^+$. Since $-\IM\log\Sigma_{\nu_n}(z)
\le-\IM\log\Sigma_{\mu_n}(z)$ for $z\in\Gamma_{\alpha,\beta,\Delta}^+
\setminus\mathbb C^-$, we may apply Theorem~\ref{3.8th} to 
$\{-\IM\log\Sigma_{\nu_n}(z)\}$
and obtain that there exists a~subsequence $\{n'\}$ of $\{n\}$ such that
$\{\IM\log\Sigma_{\nu_{n'}}(z)\}$ converges uniformly on the~compact 
subsets of $\Gamma_{\alpha,\beta,\Delta}^+$ to a~harmonic 
function $U(z)$.
Consider the~sequence of analytic functions  
$\{\log(\Sigma_{\nu_{n'}}(z)a(\nu_{n'}))\}
=\{\log\Sigma_{\nu_{n'}\boxtimes\delta_{a(\nu_{n'})}}(z)\}$ in 
$\Gamma_{\alpha,\beta,\Delta}^+$. 
This sequence converges at the~point $(\beta+\Delta)e^{i(\pi+\alpha)/2}/2$ 
and then, by Theorem~\ref{3.9th}, 
$\{\Sigma_{\nu_{n'}\boxtimes\delta_{a(\nu_{n'})}}(z)\}$
converges uniformly on every compact subset of $\Gamma_{\alpha,\beta,
\Delta}^+$. 
Hence the~assumptions of Proposition~\ref{3.12pro} hold for the~sequence 
$\{\nu_{n'}\boxtimes\delta_{a(\nu_{n'})}\}$ and  
the~sequence of p-measures converges in the~weak topology 
to a~p-measure $\nu$. 
Iterating this argument, we see that 
$\{\rho_{n'}\}$ converges weakly to $\rho$, say, 
where $\rho_{n'}$ is some
cofactor of $\nu_{n'}\boxtimes\delta_{a(\nu_{n'})}$ in $\mu_{n'}$. By 
continuity of the~multiplicative free convolution (Proposition~\ref{6.2pr}), 
we have $\nu\boxplus\rho=\mu$, so $\nu$ is
a~factor of $\mu$. Thus Proposition~\ref{6.3pr} 
is proved for $(\mathcal M_+,\boxtimes)$.

It remains to consider the~case of the~semigroup $(\mathcal M_*, \boxtimes)$.
By Proposition~\ref{3.13pro}, the~sequence 
$\{\Sigma_{\mu_n}\}$ converges uniformly on compact subsets 
of some $\mathbb D_{\alpha}$ to the~function $\Sigma_{\mu}$.
Recall that for every $\mu\in\mathcal M^{\alpha}_*$
the~function $\log\Sigma_{\mu}(z)$ is analytic in 
$\mathbb D_{\alpha}$ and $\RE\log\Sigma_{\mu}(z)=\log |\Sigma_{\mu}(z)|\ge 0$ 
for $z\in\mathbb D_{\alpha}$. Since 
$\log|\Sigma_{\nu_n}(z)|\le\log|\Sigma_{\mu_n}(z)|$ for 
$z\in\mathbb D_{\alpha}$, we may apply 
Theorem~\ref{3.8th} to $\{\log|\Sigma_{\nu_n}(z)|\}$
and obtain that there exists a~subsequence $\{n'\}$ of $\{n\}$ such that
$\{\log|\Sigma_{\nu_{n'}}(z)|\}$ converges uniformly on compact 
subsets of $\mathbb D_{\alpha}$ to a~harmonic function $U(z)\ge 0$.
Consider the~sequence of analytic functions in $\mathbb D_{\alpha}$ 
$\{\log(\Sigma_{\nu_{n'}}(z)a_{n'})\}
=\{\log\Sigma_{\nu_{n'}\boxtimes\delta_{a(\nu_{n'})}}(z)\}$. 
This sequence converges at the~point $\alpha/2$ and
thus, by Theorem~\ref{3.9th}, the~sequence $\{\Sigma_{\nu_{n'}\boxtimes
\delta_{a(\nu_{n'})}}(z)\}$
converges uniformly on every compact subset of $\mathbb D_{\alpha}$. 
Hence the~assumptions of Proposition~\ref{3.13pro} hold for the~sequence 
$\{\nu_{n'}\boxtimes\delta_{a(\nu_{n'})}\}$. By this proposition, 
this sequence of p-measures converges in the~weak topology 
to a~p-measure $\nu$. 
Iterating the~argument, we see that 
$\{\rho_{n'}\}$ converges weakly to $\rho$, say, 
where $\rho_{n'}$ is some
cofactor of $\nu_{n'}\boxtimes\delta_{a(\nu_{n'})}$ in $\mu_{n'}$. By 
continuity of the~multiplicative free convolution (Proposition~\ref{6.2pr}), 
we have $\nu\boxplus\rho=\mu$, so $\nu$ is
a~factor of $\mu$, thus proving Proposition~\ref{6.3pr} 
for $(\Cal M_*,\boxtimes)$.
 
Proposition~\ref{6.3pr} is completely proved.
\end{Proof}

\begin{proposition}\label{6.3apr}
If $\{\mu_n^*\}$ is any element of $M/E$, and if, for each
$n\in\Bbb N$, $\mu_n^*=\nu_n^*\circ\rho_n^*$, where $\nu_n^*,\rho_n^*\in\mathbf M/E$ and
$\nu_n^*$ is a~factor of $\nu_{n+1}^*$, then 
$\{\nu_{n}^*\}$ and $\{\rho_{n}^*\}$ in $M/E$ and $l\{\nu_{n}^*\}$ and
$l\{\rho_{n}^*\}$ are factors of $l\{\mu_n^*\}$.
\end{proposition}

\begin{proof}
We prove this proposition in the~same way as Proposition~\ref{6.3pr}.
We keep all notations  and demonstrate the~difference between proofs of
Proposition~\ref{6.3pr} and Proposition~\ref{6.3apr} in connection with
the~semigroup $(\mathcal M,\boxplus)$. We shall repeat all arguments  
accepting the following. By the~assumptions of Proposition~\ref{6.3apr}
we have
$$
-\Im \varphi_{\nu_n}(z)\le-\Im \varphi_{\nu_{n+1}}(z)\le 
-\Im \varphi_{\mu_{n+1}}(z),\quad z\in\Gamma_{\alpha,\beta},
$$
for all $n\ge 1$. Hence, by Theorem~\ref{3.8ath}, $\{-\Im \varphi_{\nu_n}(z)\}$
converges uniformly on compact subsets of $\Gamma_{\alpha,\beta}$ to 
a~harmonic function $U(z)\ge 0$. Then we repeat the~arguments of
Proposition~\ref{6.3pr} for $\{\varphi_{\nu_n}(z)+a(\nu_n)\}$ and obtain that
$\{\nu_n\boxplus\delta_{a(\nu_n)}\}$ converges weakly to some p-measure $\nu$.
Furthermore, as in the~proof of Proposition~\ref{6.3pr}
we see that 
$\{\rho_{n}\}$ converges weakly to $\rho$, say, 
where $\rho_{n}$ is 
the~cofactor of $\nu_{n}\boxplus\delta_{a(\nu_{n})}$ in $\mu_{n}$
such that $\rho_n^*=E\rho_n$,
and $\nu\boxplus\rho=\mu$, so $\nu$ and $\rho$ are
factors of $\mu$. 
Thus, Proposition~\ref{6.3apr} is proved for
the~semigroup $(\mathcal M,\boxplus)$.

The proof for the~semigroups $(\mathcal M_+,\boxtimes)$ and  
$(\mathcal M_*,\boxtimes)$ is similar therefore we omit it.
\end{proof}  
\begin{proposition}\label{6.4pr}
Let $\{\nu_n\}$ be a~sequence of p-measures 
$\nu_n\in\mathbf M^{\gamma}$ and let $D^{\gamma}(\nu_n^*)\to 0$ as 
$n\to\infty$. For every fixed $\gamma'\in\mathcal G$ 
there exists $n(\gamma')$ such that 
$\nu_n\circ \delta_{a(\nu_n)}\in\mathbf M^{\gamma'}$ 
for $n\ge n(\gamma')$.
\end{proposition}
\begin{Proof}
Let us prove the~assertion for the~semigroup
$(\Cal M,\boxplus)$.
Let $\{\nu_n\}$ denote a~sequence of p-measures 
$\nu_n\in\mathcal M^{(\alpha,\beta)}$ and 
let $\IM\varphi_{\nu_n}(z_0)\to 0$ as $n\to\infty$,
where $z_0:=i(\beta+1)$. 
We shall prove that for every fixed $\alpha'>\alpha$ and
$0<\beta'<\beta$, there exists $n(\alpha',\beta')$ such that 
$\nu_n\boxplus\delta_{a(\nu_n)}\in\mathcal M^{(\alpha',\beta')}$ for 
$n\ge n(\alpha',\beta')$.

The~functions $\IM\varphi_{\nu_n}(z)$ are harmonic non positive-valued
functions in $\Gamma_{\alpha,\beta}$. The\-re\-fore, by the~theorem
of the~mean for harmonic functions,
\begin{equation}\label{6.6}
|\IM\varphi_{\nu_n}(z_0)|=\frac 1{2\pi}\int\limits_0^{2\pi}
|\IM\varphi_{\nu_n}(z_0+\varepsilon e^{i\theta})|\,d\theta
\end{equation}
for sufficiently small $\varepsilon>0$.
Since $\RE\varphi_{\nu_n\boxplus\delta_{a(\nu_n)}}(z_0)=0$,
we have, by Schwarz's formula (see Markushevich (1965), v. 2, p. 151), 
$$
\varphi_{\nu_n\boxplus\delta_{a(\nu_n)}}(z)=
\frac i{2\pi}\int\limits_0^{2\pi}\IM\varphi_{\nu_n}(z_0
+\varepsilon e^{i\theta}/2)\,\frac{\varepsilon e^{i\theta}/2+
(z-z_0)}{\varepsilon e^{i\theta}/2-(z-z_0)}\,d\theta
$$
for $z\in\mathbb C$ such that $|z-z_0|<\varepsilon/2$. Using (\ref{6.6}),
we obtain from this formula, for $|z-z_0|\le \varepsilon/4$,
\begin{align}\label{6.7}
|\varphi_{\nu_n\boxplus\delta_{a(\nu_n)}}(z)|&\le
\frac 1{2\pi}\int\limits_0^{2\pi}|\IM\varphi_{\nu_n}(z_0 
+\varepsilon e^{i\theta}/2)|\frac{\varepsilon/2+|z-z_0|}
{\varepsilon/2-|z-z_0|}\,d\theta\notag\\
&\le\frac 3{2\pi}\int\limits_0^{2\pi}|\IM\varphi_{\nu_n}(z_0 
+\varepsilon e^{i\theta}/2)|\,d\theta=
3 |\IM\varphi_{\nu_n}(z_0)|. 
\end{align}
Since $|\IM\varphi_{\nu_n}(z_0)|\to 0$ as $n\to\infty$, we conclude
from (\ref{6.7}) that $\varphi_{\nu_n\boxplus\delta_{a(\nu_n)}}(z)\to 0$ 
uniformly for $|z-z_0|\le\varepsilon/4$ as $n\to\infty$. 
Recall that $F_{\nu_n\boxplus\delta_{a(\nu_n)}}^{(-1)}(z)=
z+\varphi_{\nu_n\boxplus\delta_{a(\nu_n)}}(z)$. Therefore we easily 
conclude that $F_{\nu_n\boxplus\delta_{a(\nu_n)}}=z+o(1)$ and hence
$G_{\nu_n\boxplus\delta_{a(\nu_n)}}(z)=1/z+o(1)$ uniformly for 
$|z-z_0|\le\varepsilon/8$ as $n\to\infty$. From the~last relation we see that
\begin{equation}\label{6.8}
\IM \Big(G_{\nu_n\boxplus\delta_{a(\nu_n)}}(z_0)-\frac 1{z_0}\Big)
=\frac 1{\IM z_0}\int\limits_{\Bbb R}
\frac{u^2}{u^2+(\IM z_0)^2}\,
\nu_n\boxplus\delta_{a(\nu_n)}(du)\to 0,
\,\, \text{as\,\,}n\to \infty.
\end{equation}
From (\ref{6.8}) it follows that $\nu_n\boxplus\delta_{a(\nu_n)}$
converges weakly to $\delta_0$ as $n\to\infty$. The~last relation 
easily implies the~assertion of Proposition~\ref{6.4pr} 
for $(\mathcal M,\boxplus)$.

We shall consider the~assertion for the~semigroup $(\mathcal M_+,\boxtimes)$.
Let $\{\nu_n\}$ denote a~sequence of p-measures 
$\nu_n\in\mathcal M^{(\alpha,\beta,\Delta)}_+$ and assume that
$\IM\log\Sigma_{\nu_n}(z_0)$$\to 0$ as $n\to\infty$, where  
$z_0:=(\beta+\Delta)e^{i(\pi+\alpha)/2}/2$. We shall prove that 
for eve\-ry fixed $0<\alpha'<\alpha$,
$0<\beta'<\beta$, and $\Delta'>\Delta$, there exists 
$n(\alpha',\beta',\Delta')$ such that 
$\nu_n\boxplus\delta_{a(\nu_n)}\in\mathcal M^{(\alpha',\beta',
\Delta')}$ for $n\ge n(\alpha',\beta',\Delta')$.

The~functions $\IM\log\Sigma_{\nu_n}(z)$ are harmonic non positive-valued
functions in the domain $\Gamma_{\alpha,\beta,\Delta}^+\cap\mathbb C^+$ such
that $\log|\Sigma_{\nu_n\boxtimes\delta_{a(\nu_n)}}(z_0)|=0$. Repeating
the~previous arguments for the~analytic functions  
$\log\Sigma_{\nu_n\boxtimes\delta_{a(\nu_n)}}(z)$ on 
$\Gamma_{\alpha,\beta,\Delta}^+$,
we see that these functions 
$\log\Sigma_{\nu_n\boxtimes\delta_{a(\nu_n)}}(z)\to 0$ 
uniformly for $|z-z_0|\le\varepsilon/4$ as $n\to\infty$, where
as above $\varepsilon>0$ is chosen sufficiently small.
Hence $\Sigma_{\nu_n\boxtimes\delta_{a(\nu_n)}}(z)$ converges to $1$ 
uniformly for $|z-z_0|\le\varepsilon/4$ as $n\to\infty$.
Recalling the~definition of $\Sigma_{\nu_n\boxtimes\delta_{a(\nu_n)}}(z)$,
we conclude from this relation that 
\begin{equation}\label{6.9}
K_{\nu_n\boxtimes\delta_{a(\nu_n)}}(z)\to z
\end{equation} 
uniformly for $|z-z_0|\le\varepsilon/8$ as $n\to\infty$.
Taking into account the~representation (\ref{2.8}) for the~functions
$K_{\nu_n\boxtimes\delta_{a(\nu_n)}}$, we have 
\begin{equation}\label{6.10}
\frac 1z K_{\nu_n\boxtimes\delta_{a(\nu_n)}}(z)=a_n+\int\limits
_{(0,\infty)}\frac{\tau_n(dt)}{t-z},\qquad z\in\mathbb C^+,
\end{equation}
where $a_n\ge0$ and the~nonnegative measures $\tau_n$ satisfy (\ref{2.9}).
We easily see from (\ref{6.9}) and (\ref{6.10}) that
$$
a_n+\int\limits_{(0,\infty)}\frac{\tau_n(dt)}{t+1}\le c(z_0)<\infty,
\qquad n=1,2,\dots,
$$
where $c(z_0)$ is a~positive constant depending on $z_0$.
This estimate implies the~inequality $|K_{\nu_n\boxtimes
\delta_{a(\nu_n)}}(z)/z|\le c(\alpha',\beta',\Delta'),\,n=1,2,\dots$, 
in every domain $\Gamma^+_{\alpha',\beta',\Delta'}$ with $0<\alpha'<\alpha,
0<\beta'<\beta$, and $\Delta'>\Delta$. Here $ c(\alpha',\beta',\Delta')$
denotes a~positive constant depending on $\alpha',\beta'$, and $\Delta'$. 
Hence, by (\ref{6.10}), we conclude with 
the~help of Vitali's theorem (see Theorem~\ref{3.10th} of Section~3 that
relation (\ref{6.9}) holds for $z\in\Gamma^+_{\alpha',\beta',\Delta'}$
with any fixed $\alpha',\beta',\Delta'$. 
From this relation we immediately obtain the~assertion 
for $(\mathcal M_+,\boxtimes)$.

Let us prove the~assertion for $(\mathcal M_*,\boxtimes)$.
Let $\{\nu_n\}$ be a~sequence of p-measures 
$\nu_n\in\mathcal M^{\alpha}_*$ such that 
let $\log|\Sigma_{\nu_n}(z_0)|\to 0$ as $n\to\infty$. Here 
$z_0:=\alpha/2$. We shall prove that for every fixed $\alpha'>\alpha$ 
there exists $n(\alpha')$ 
such that $\nu_n\boxplus\delta_{a(\nu_n)}\in\mathcal M^{\alpha'}$ 
for $n\ge n(\alpha')$.

The~functions $\log|\Sigma_{\nu_n}(z)|$ are harmonic non negative-valued
functions on $\mathbb D_{\alpha}$ such that
$\IM\log\Sigma_{\nu_n\boxtimes\delta_{a(\nu_n)}}(z_0)=0$. Repeating
the~previous arguments for the~ana\-lytic functions 
$\log\Sigma_{\nu_n\boxtimes\delta_{a(\nu_n)}}(z)$ on $\mathbb D_{\alpha}$,
we see that $\log\Sigma_{\nu_n\boxplus\delta_{a(\nu_n)}}(z)\to 0$ 
uniformly for $|z-z_0|\le\varepsilon/4$ as $n\to\infty$. Hence 
$\Sigma_{\nu_n\boxtimes\delta_{a(\nu_n)}}(z)\to 1$ 
uniformly for $|z-z_0|\le\varepsilon/4$ as $n\to\infty$.
Recalling the~definition of $\Sigma_{\nu_n\boxtimes\delta_{a(\nu_n)}}(z)$,
we conclude from this relation that $Q_{\nu_n\boxtimes\delta_{a(\nu_n)}}(z)
\to z$ uniformly for $|z-z_0|\le\varepsilon/8$ as $n\to\infty$.
Since for $\mu\in\mathcal S_*$ the~function $Q_{\mu}(z)/z$ is analytic and
$|Q_{\mu}(z)/z|\le 1$ for $z\in\mathbb D$, we conclude, by Vitali's theorem,
that $Q_{\nu_n\boxtimes\delta_{a(\nu_n)}}(z)\to z$ uniformly for $z\in
\mathbb D_{\alpha'}$ for any fixed $\alpha'\in(\alpha,1)$. This implies at once
the~assertion for $(\mathcal M_*,\boxtimes)$.

Thus, Proposition~\ref{6.4pr} is proved.
\end{Proof}

\begin{proposition}\label{6.5pr}
Let $\{\mu_{nk}:n\ge 1,1\le k\le n\}$ be an~array of infinitesimal
measures in $\mathbf M^{\gamma}$. Then, for every fixed $\gamma\in\mathcal G$,
$\max_{1\le k\le n}D^{\gamma}(\mu_{nk}^*) \to 0$ as $n\to\infty$.
\end{proposition}

This proposition follows from obvious calculations which  
we omit.

{\bf Proof of Theorem~2.10}.
Let $\{\mu_{j,s}:1\le s\le j<\infty\}$ be a~array of infinitesimal 
measures in $\mathbf M$ and let
$$
l\{\mu^*_j\}=\mu^*\in\mathbf M/E,\quad
\text{ where }\quad\mu^*_j=\mu^*_{j,1}\circ\dots\circ\mu^*_{j,j}. 
$$
It is clear that $\mu$ and $\mu_j,\,j=1,\dots$, such that $\mu^*=E\mu$ and
$\mu_j^*=E\mu_j,\,j=1,\dots$,  
belong to the~set
$\mathbf M^{\gamma}$ with some $\gamma\in\mathcal G$ and, by Proposition~\ref{6.1} and
\ref{6.5pr},
$$
\lim_{j\to\infty}\max_{1\le s\le j} D^{\gamma}(\mu^*_{j,s})=0.
$$

Let $k\in\mathbb N$ be given. Since 
$\max_{1\le s\le j} D^{\gamma}(\mu^*_{j,s})\to 0$ as $j\to\infty$
and, by Propositions~\ref{6.2bpr}, $D^{\gamma}(\mu_j^*)
\to D^{\gamma}(\mu^*)$ as $j\to\infty$,
given $n\in\Bbb N$ there exists $j(n)$ so large that we may
group $\mu^*_{j(n),s},\,1\le s\le j(n)$, via
free convolutions into $\nu^*_{k,l,n}$ say, $1\le l\le k$, such that
$\mu^*_{j(n)}=\nu^*_{k,1,n}\circ\dots\circ\nu^*_{k,k,n}$ and
$| D^{\gamma}(\nu^*_{k,l,n})-D^{\gamma}(\mu^*)/k|\le 1/n$
for $1\le l\le k$. By Proposition~\ref{6.3pr}, there is a~decomposition
$\mu^*=\nu^*_{k,1}\circ\dots\circ\nu^*_{k,k}$ with
$D^{\gamma}(\nu^*_{k,l})=D^{\gamma}(\mu^*)/k$ for
$1\le l\le k$. Hence 
there is a~decomposition $\mu=\nu_{k,1}\circ\dots
\circ\nu_{k,k}$ for every $k\in\mathbb N$ such that 
$\nu_{k,l}\in\mathbf M^{\gamma}$ and $D^{\gamma}(\nu^*_{k,l})=D^{\gamma}(\mu^*)/k$ 
for $1\le l\le k$.
Applying Proposition~\ref{6.4pr} to the~sequence $\{\nu_{k,l}\}$ with 
$k=1,\dots$ and $l=1,\dots,k$, we see that, for any fixed $\gamma\in\mathcal G$,  
there exists $k(\gamma)$ such that $\nu_{k,l}\in\mathbf M^{\gamma}$ for
$k\ge k(\gamma),\,l=1,\dots,k$. 

From this relation we see in the~case
of $(\mathcal M,\boxplus)$ that $\varphi_{\mu}(z)$ admits an~analytic 
continuation in $\mathbb C^+$ with values $\mathbb C^-\cup\mathbb R$ and hence,
by the~Bercovici and Voiculescu result (see Section~2), $\mu$ is infinitely
divisible p-measure. This assertion follows from 
Lemma~\ref{4.1l} as well.

In the~case of $(\mathcal M_+,\boxtimes)$ we note that the~function
$\log \Sigma_{\mu}(z)$ admits an~analytic 
continuation in $\mathbb C^+$ with values $\mathbb C^-\cup\mathbb R$. In addition
this function is analytic and real-valued on the negative half-line.
Again, by the~Bercovici and Voiculescu result (see Section~2), $\mu$ 
is infinitely divisible p-measure. We can obtain this result with the~help
of Lemma~\ref{5.3} as well.

In the~case of $(\mathcal M_*,\boxtimes)$ the~function $\log \Sigma_{\mu}(z)$ 
admits an~analytic continuation in $\mathbb D$ with values in
$i(\mathbb C^+\cup\mathbb R)$. By the~Bercovici and Voiculescu result 
(see Section~2), $\mu$ is infinitely divisible p-measure. This result
follows from Lemma~\ref{5.6} as well.

The~theorem is proved.
$\square$

{\bf Proof of Theorem~2.11}.
Let $\mu$ be not a~Dirac measure and have no indecomposable factors. 
We shall show that $\mu$ is infinitely divisible p-measure.
Without loss of generality we assume that $\mu\in \mathbf M^{\gamma}$
with some $\gamma\in\mathcal G$.
By Proposition~\ref{6.2apr}, $D^{\gamma}(\mu^*)>0$. Let us show
that $\inf_{\nu\in S}D^{\gamma}(\nu^*)=0$, where $S$ is the~set of all
factors $\nu$ of $\mu$ such that $\nu^*\ne e$. For suppose not:
call the~non-zero infimum $b$. By Proposition~\ref{6.1pr}, Proposition~\ref{6.3pr}, and
Proposition~\ref{6.2bpr}, there is $\nu_0^*\in\mathbf M^{\gamma}/E$
such that $\nu_0^*=E\nu_0$ with $\nu_0\in S$ and 
$D^{\gamma}(\nu^*_0)=b$. Since $\nu_0$ 
is a~factor of $\mu$, it has factors from $S$, say $\nu_1$, on which 
$D^{\gamma}(\nu_1^*)$ is less then $b$, a~contradiction.

Consider now all n-fold decompositions $\mathcal D: \,\mu=\nu_1\circ\nu_2\circ
\dots\circ\nu_n$. By Proposition~\ref{6.1pr}, $\nu_j\in \mathbf M^{\gamma}
,\,j=1,\dots,n$. Define $m(\mathcal D):=\min_{j=1,\dots,n} D^{\gamma}(\nu^*_j)$, 
so that $0\le m(\mathcal D)\le (1/n)D^{\gamma}(\mu^*)$. 
Let $h:=\sup_{\mathcal D} m(\mathcal D)$. 
By Proposition~\ref{6.3pr} and Proposition~\ref{6.2bpr}, there is
a~decomposition $\mathcal D'$ such that $m(\mathcal D')=h$. Assume that 
$\nu_1',\nu_2',\dots,\nu_n'$ have been odered  so that 
$D^{\gamma}((\nu'_j)^*)$ increases with $j$; then 
$h=D^{\gamma}((\nu'_1)^*)$. Now either all $D^{\gamma}$'s are
equal or there is a~least $l$ with $1\le l<n$ and
$D^{\gamma}((\nu'_l)^*)<D^{\gamma}((\nu'_{l+1})^*)$. By the~first
paragraph of this proof we can write
$$
\nu_{l+1}'=\rho'_1\circ\rho_2'\circ\dots\circ\rho_l'\circ \nu''_{l+1},
$$
where $D^{\gamma}((\rho'_1)^*),\dots,D^{\gamma}((\rho'_l)^*)$ 
are all arbitrary small but non-zero. If we choose them so that
$D^{\gamma}((\nu''_{l+1})^*)>D^{\gamma}((\nu'_l)^*)$, we shall obtain
a~decomposition
$$
\mathcal D'':\,\mu=(\nu_1'\circ\rho_1')\circ\dots\circ(\nu_l'\circ\rho_l')
\circ \nu_{l+1}''\circ\nu_{l+2}'\dots\circ\nu_n'
$$
which has $m(\mathcal D'')>h$, a~contradiction. So in $\mathcal D'$ all 
the~$D^{\gamma}$'s
must be equal, entailing $D^{\gamma}((\nu_j')^*)=(1/n)D^{\gamma}(\mu^*)$
for each $j,\,1\le j\le n$.

We have  a~$\mathcal D'$ for every $n$: call it $\mathcal D'(n):\,
\mu=\nu'_1(n)\circ\nu'_2(n)
\circ\dots\circ\nu'_n(n)$. The~$\mathcal D'(n)$ forms
a~triangular array whose row $\circ$-products are equal to $\mu$.
Further, $\max_{j=1,\dots,n}D^{\gamma}((\nu_j'(n))^*)$ $=(1/n)D^{\gamma}(\mu^*)$.
It follows from Theorem~2.10 that $\mu$ is infinitely divisible p-measure.
The~theorem is proved completely.
$\square$

{\bf Proof of Theorem~2.12}. It is clear that $\mu\in\mathbf M^{\gamma}$
with some $\gamma\in\mathcal G$. Hence the homomorphism $D^{\gamma}(\mu^*)$
is defined and finite.
By Proposition~\ref{6.1}, all factors of $\mu$ belong to 
$\mathbf M^{\gamma}$.  

Denote by $S_1$ the~set of all indecomposable factors of $\mu$.
By the~assumption of the~theorem, this set is not empty. Denote
$d_1:=\sup_{\nu\in S_1}D^{\gamma}(\nu^*)$. Let an~indecomposable 
factor $\nu_1$ of $\mu$ have the~property 
$D^{\gamma}(\nu^*_1)\ge d_1/2$.
Since $\nu_1$ is a~factor of $\mu$, we have $\mu=\nu_1\circ
\rho_1$, where $\rho_1\in\mathbf M^{\gamma}$.

If $\rho_1$ is an~indecomposable p-measure, the~proof is completed.
In other case we repeat the~previous arguments for $\rho_1$ and
arrive at the~relation $\mu=\nu_1\circ\nu_2\circ\rho_2$,
where $\nu_2$ and $\rho_2$ belong to $\mathbf M^{\gamma}$, and $\nu_2$
is an~indecomposable p-measure such that
$$
D^{\gamma}(\nu^*_2)\ge d_2/2,\qquad d_2:=\sup_{\nu\in S_2}
D^{\gamma}(\nu^*),
$$  
where $S_2\,(\ne\emptyset)$ is the~set of all indecomposable factors 
of $\rho_1$.

Repeating next steps in the~same way we obtain
\begin{equation}\label {2.12,1}
\mu=\nu_1\circ\nu_2\circ\dots\circ\nu_k\circ\rho_k,
\end{equation}
where $\nu_1,\nu_2,\dots,\nu_k$ are indecomposable factors
of $\mu$ and $\rho_k$ is an~indecomposable p-measure, then the~proof 
is completed, or we have an~infinite sequence of the~relations 
(\ref{2.12,1}). Recall that all p-measures $\nu_1,\nu_2,\dots,\nu_k$
and $\rho_k$ belong to $\mathbf M^{\gamma}$. In addition we have
\begin{equation}\label {2.12,2}  
D^{\gamma}(\nu^*_k)\ge d_k/2,\qquad d_k:=\sup_{\nu\in S_k}
D^{\gamma}(\nu^*),
\end{equation}
where $S_k\,(\ne\emptyset)$ is the~set of all indecomposable factors 
of $\rho_{k-1}$ and, for $n>m$,
\begin{equation}\label {2.12,3}
\rho_m=\nu_{m+1}\circ\nu_{m+2}\dots\circ\nu_n\circ\rho_n. 
\end{equation}
By (\ref{2.12,1}) and (\ref{2.12,2}), we see that 
\begin{equation}\label {2.12,4}
\frac 12\sum_{k=1}^{\infty}d_k\le\sum_{k=1}^{\infty}D^{\gamma}(\nu^*_k)
\le D^{\gamma}(\mu^*)<\infty.
\end{equation}

Applying Proposition~\ref{6.3apr} to (\ref{2.12,1}), we see that
$\{\nu^*_1\circ\nu^*_2\circ\dots\circ\nu^*_k\}$ and $\{\rho^*_k\}$ are 
elements of $M/E$ and $l\{\nu^*_1\circ\nu^*_2\circ\dots\circ\nu^*_k\}$
and $l\{\rho^*_k\}$ are factors of $\mu^*$. 

It remains to show that
$l\{\rho^*_k\}$ has no indecomposable factors.
Applying Proposition~\ref{6.3apr} to (\ref{2.12,3}), we obtain
the~relation 
\begin{equation}\label {2.12,5}
\rho^*_m=l\{\nu^*_{m+1}\circ\nu^*_2\circ\dots\circ\nu^*_n\}\circ
l\{\rho^*_n\},\qquad m=1,2,\dots.
\end{equation}
Let $l\{\rho^*_n\}$ have an~indecomposable factor $\kappa^*$ such that
$\kappa^*=E\kappa$ with $\kappa\in\mathbf M^{\gamma}$.
We see from (\ref{2.12,5}) that $\kappa$ is a~factor of $\rho_m$
for every $m=1,2,\dots$ and hence belongs to all sets $S_m,\,m=1,\dots$.
By Proposition~\ref{6.2apr}, $D^{\gamma}(\kappa^*)>0$. On the~other
hand $d_k\ge D^{\gamma}(\kappa^*)>0$ for all $k=1,\dots$ and the~series
$\sum_{k=1}^{\infty}d_k$ should deverge. By (\ref{2.12,4}), we arrive
at a~contradiction. Thus, $l\{\rho^*_n\}$ has no indecomposable factors
and the~theorem is proved completely.
$\square$

\section{Description of the~class $I_0$ in $(\mathbf M,\circ)$ 
and nonuniqueness of the~representation (\ref{2.27})}

Our next step is to prove Theorem~2.13 
for the~semigroups $(\mathbf M,\circ)$.

{\it Proof of Theorem~$2.13$}.
Since $I_0$ is the~set of i.d. elements all 
of whose factors are i.d., we should prove
that any nontrivial i.d. element in $(\mathbf M,\circ)$
has a~non i.d. factor. 

$\mathbf {7.1}$. Consider the~case of $(\mathcal M,\boxplus)$. As we saw 
in Section~2 a~measure $\mu\in\mathcal M$ 
is $\boxplus$-i.d. if and only if $\varphi_{\mu}$ 
admits a~representation~(\ref{2.18}).

Let $a<b$ be real numbers such that $\nu([a,b])>0$. Consider the~function 
\begin{equation}\label{7.1}
\phi(z):=\varepsilon_0\int\limits_{[a,b]}\ffrac{(1+u^2)\nu(du)}{z-u},
\quad z\in\mathbb C^+,
\end{equation}
with sufficiently small $\varepsilon_0>0$, depending on $\nu,a,b$ only,
and introduce the~functions $f_j(z):=2\phi(z)+(-1)^j\varepsilon
\phi^2(z),\,j=1,2$.  
In order to prove the~theorem we need the~following lemma.
\begin{lemma}
For sufficiently small $\varepsilon >0$, there exist $\mu_j\in\mathcal M,
\,j=1,2$, such that $\varphi_{\mu_j}(z)=f_j(z)$ for $z$, 
where $\varphi_{\mu_j}(z)$ is defined.
\end{lemma}

\begin{Proof}
Consider the~function
$$
F(z):=z+\phi(z),\quad z\in\mathbb C^+.
$$
By Lemma~4.1, $F:\mathbb C^+\to\mathbb C$ takes every value in $\mathbb C^+$
precisely once. The~inverse function $F^{(-1)}:\mathbb C^+\to\mathbb C^+$
thus defined is in the~class $\mathcal F$. 

Denote by $\mathbb C^+_{\eta}$
the half-plane $\IM z>\eta$ for any real $\eta$. Consider the~domains 
$\Omega_{\eta}:=F^{(-1)}(\mathbb C^+_{\eta}),\,\eta\in(0,1]$. 
Note that $\Omega_{\eta}$ are domains such that 
$\Omega_{\eta}\subset\mathbb C^+$.

The~boundary of $\Omega_{\eta}$ is a~curve $\gamma_5(\eta)$
characterized by the~functional equation $\IM F(x+iy)=\eta,\,x\in\Bbb R,\,
y>0$. Using (\ref{7.1}) rewrite this equation in the~form
\begin{equation}\label{7.1a}
y\Big(1-\varepsilon_0\int\limits_{[a,b]}\frac{(1+u^2)\,\nu(du)}{(u-x)^2+y^2}
\Big)=\eta.
\end{equation}
It is easy to see that (\ref{7.1a}) has an~unique solution for every 
fixed $x\in\mathbb R$.
Since $F'(z)\ne 0$ for $z\in\Omega_{\eta},\,\eta>0$, we conclude, by 
the~implicit function theorem, that $\gamma_5(\eta),\,\eta>0$, are
smooth Jordan curves.

In view of (\ref{7.1a}), we see that the~curves $\gamma_5(\eta)$ lie 
in the~union of domains $|z-(a+b)/2|< b-a,\,\Im z>0$ and  
$0<\IM z<2\eta$. 
Denote by $\tilde\gamma_5(\eta)$ the~part of the~curve $\gamma_5(\eta)$
contained in the~disk $|z-(a+b)/2|\le b-a$. Let us show that
\begin{equation}\label{7.2}
\sup_{0<\eta\le 1}\sup_{z\in\tilde\gamma_5(\eta)}|\phi(z)|:=A_1<\infty.
\end{equation} 
Assume that $A_1=\infty$. Then for any
$N>1$ there exist $\eta_0\in(0,1]$ and $z_0=z_0(\eta_0)\in
\tilde\gamma_5(\eta_0)$ 
such that $\IM F(z_0)>0$ and $|F(z_0)|>N$. Since $F(x)$ is real-valued 
for $x\in\mathbb R\setminus[a,b]$ and $\IM z\IM F(z)> 0$ for $|z|>|a|+|b|, \,\Im z\ne 0$,  
and $|\phi(z)|\to 0$ for $|z|=R\to\infty$,
we easily see, by Rouch\'e's theorem, that for $w\in\mathbb C^+$ 
with sufficiently large modulus 
the~equation $w=F(z)$ has a~solution $z=z(w)\in\mathbb C^+$ such that 
$z(w)=w+O(1)$. Recalling that $F(z):\mathbb C^+\to\mathbb C$ takes every
value in $\mathbb C^+$ precisely once, we conclude that $z_0=F(z_0)+O(1)$,
a~contradiction for sufficiently large $N$. 

Now let us show, for $j=1,2$, that $F_j(z):=z+f_j(z):\mathbb C^+\to\mathbb C$ 
takes every value in $\mathbb C^+$ precisely once in $\mathbb C^+$.  
By (\ref{7.2}), we see that the~inequality
$1+2(-1)^j\varepsilon\RE \phi(z)>0$ holds for 
$z\in \gamma_5(\eta),\,\eta\in(0,1]$, and for sufficiently small 
$\varepsilon>0$.
Since $\IM \phi(z)\le 0$ and $\IM F(z)=\eta$ for the~same $z$, we
conclude from the~formula
\begin{equation}\label{7.3}
\IM F_j(z)=\IM z+\IM\phi(z)+\IM \phi(z)\Big(
1+2(-1)^j\varepsilon\RE \phi(z)\Big)
\end{equation}
that $\IM F_j(z)\le \eta$ for $z\in \gamma_5(\eta),\,\eta\in(0,1]$.

Choose $R>1$ to be determined later sufficiently large.
Denote by $\tilde a, \,\Re \tilde a<0$, and $\tilde b,\,\RE \tilde b>0$,
points of intersection of the~curve $\gamma_5(\eta)$ with 
the~circle $|z|=R$. 
Consider the~closed rectifiable curve $\gamma_6=\gamma_6(\eta),\,\eta
\in(0,1]$, consisting of $\gamma_{6,1}$:
the~part of $\gamma_5(\eta)$ lying in the~disc $|z|\le R$, connecting
$\tilde a$ to $\tilde b$, and the~arc $\gamma_{6,2}:Re^{i\theta},\,
\arg \tilde b<\theta<\arg\tilde a$, connecting $\tilde b$ to $\tilde a$.

Fix $w\in\mathbb C^+$. Assume that $\eta<\Im w$. If $z$ traverses 
$\gamma_{6,1}$, the~image
$\zeta=F_j(z)$ lies in the~half-plane $\IM \zeta\le \eta$.
If $z$ traverses $\gamma_{6,2}$, the~image $\zeta=F_j(z)$ is equal to
$Re^{i\theta}+o(1),\,\arg \tilde b<\theta<\arg\tilde a$, as 
$R\to\infty$. We conclude that the~image
$\zeta=F_j(z)$ winds around $w$ once when $z$ runs through $\gamma_6$
in the~counter clockwise direction.
By the~argument principle, the~function $F_j(z)$ takes the~value
$w$ precisely once inside $\gamma_6$. 
Since this assertion 
holds for all sufficiently large $R>1$ and for all sufficiently small
$\eta>0$, we obtain the~desired result.

Thus we conclude that $F_j^{(-1)}:\mathbb C^+\to \mathbb C^+,\,j=1,2$, 
exist and are analytic functions. Hence $F_j^{(-1)}\in\mathcal N$. Moreover
it is easy to see that $F_j^{(-1)}\in\mathcal F$, by proving the~lemma.
\end{Proof}
We shall complete the~proof of the~theorem for the~semigroup
$(\mathcal M,\boxplus)$.

Since $\varphi_{\mu}(z)$ admits the~representation (\ref{2.18}),
we may write $\mu$ in the~form $\mu=\mu_1\boxplus\mu_2\boxplus\mu_3$,
where p-measures $\mu_1,\mu_2$, and $\mu_3$ are determined as follows.
The~measures $\mu_1,\mu_2$ are defined in Lemma~7.1. Since, by (\ref{7.1}),
$$
\varphi_{\mu_1}(z)+\varphi_{\mu_2}(z)=4\phi(z)=
4\varepsilon_0\int\limits_{[a,b]}\frac{1+uz}{z-u}\,\nu(du)-
4\varepsilon_0\int\limits_{[a,b]}u\,\nu(du),
$$
and $\varphi_{\mu}(z)=\varphi_{\mu_1}(z)+\varphi_{\mu_2}(z)
+\varphi_{\mu_3}(z)$, we easily conclude that
$$
\varphi_{\mu_3}(z):=\alpha+4\varepsilon_0\int\limits_{[a,b]}u\,\nu(du)
+(1-4\varepsilon_0)\int\limits_{[a,b]}\ffrac {1+uz}{z-u}\,\nu(du)
+\int\limits_{\Bbb R\setminus[a,b]}\ffrac {1+uz}{z-u}\,\nu(du).
$$
We consider the functions $\varphi_{\mu}(z)$ and 
$\varphi_{\mu_j}(z),\,j=1,2,3$, for $z$, where they are defined. 
By the~Bercovici and Voiculescu result~\cite{BercVo:1993} (see Section~2), 
the~measure $\mu_3$ is i.d. for $\varepsilon_0\le 1/4$.

If $\nu(\{a\})>0$,
then we choose $b=a$ in (\ref{7.1}). If $\nu(\{a\})=0$
for all $a\in\Bbb R$, then we choose the~points $a<b$ 
such that $\nu ([a,a+h])\,\,>c_0h$ and $\nu ([b-h,b])>c_0h$
for all $0<h\le h_0$, where $c_0>0$ and $h_0>0$ depend on 
the~measure $\nu$ only. Such points exist by Proposition~\ref{3.6pro}.
Let us show that under these assumptions the~function $\phi(z)$ 
(see (\ref{7.1})), has the~property
\begin{equation}\label{7.5}
-\RE \phi(a-h+ih)\to +\infty,\quad h\downarrow 0,\quad\text{and}\quad
\RE \phi(b+h+ih)\to +\infty,\quad h\downarrow 0.
\end{equation} 
This property is obvious in the~case where $a=b$ and $\nu(\{a\})>0$.
Consider the~case where $\nu(\{a\})=0$ for all $a\in\mathbb R$.
In this case we obtain for small $h>0$ the~estimates
\begin{align}
-\RE \phi(a-h+ih)&=\varepsilon_0
\int\limits_{[a,b]}\frac{(u-a+h)(1+u^2)\nu(du)}
{(u-a+h)^2+h^2}\ge\frac {\varepsilon_0}2\int\limits_{[a,b]}
\frac{\nu(du)}{u-a+h}\notag\\
&=\frac {\varepsilon_0}2\frac{\nu([a,b])}{b-a+h}
+\frac {\varepsilon_0}2\int\limits_a^b\nu([a,u))\frac{du}{(u-a+h)^2}
\notag\\
&\ge \frac{c_0\varepsilon_0}2\int\limits_a^{\min\{a+h_0,b\}}
(u-a)\frac{du}{(u-a+h)^2}\ge-\frac {c_0\varepsilon_0}4\log h.\notag
\end{align}
This lower bound implies the~first assertion of (\ref{7.5}).
In the same way we obtain the~second assertion of (\ref{7.5}).

Note that the~measures $\mu_j,\,j=1,2$, are not i.d. 
Indeed, the~functions $\varphi_j(z),\,j=1,2$, admit an~analytic 
continuation in $\mathbb C^+$. We will denote again by $\varphi_j(z),\,j=1,2$, 
respectively, these analytic continuations. Since 
$$
\Im \varphi_{\mu_j}(z)=2\Im\phi(z)(1+(-1)^j\varepsilon\Re\phi(z))
$$ 
and $\Im \phi(z)<0,\,z\in\mathbb C^+$, we see from (\ref{7.5}) that
there exist points in $\mathbb C^+$ 
where the functions $\Im \varphi_{\mu_1}(z)$ and $\Im \varphi_{\mu_2}(z)$
have positive values. By Bercovici and Voiculescu result~\cite{BercVo:1993}
(see Section~2), the~measures $\mu_j,\,j=1,2$, are not i.d.

Hence the~measure $\mu$ has a~non-i.d. factor and $\mu$ does not 
belong to the~class $I_0$.
$\square$

$\mathbf {7.2}$. Let us consider the~semigroup $(\mathcal M_+,\boxtimes)$. 
In order to prove Theorem~2.13 for this semigroup we need 
the~following auxiliary lemmas.
\begin{lemma}
For every $c>0$ , there exist 
$\mu_j\in\mathcal M_+,\,j=1,2$, such that
$$
\Sigma_{\mu_j}(z)=\exp\Big\{-\frac{cz}2+(-1)^j i\sqrt {cz}\Big\}
$$
for $z$, where $\Sigma_{\mu_j}(z)$ is defined.
\end{lemma}

Here and in the~sequel we shall choose the~principle branch of $\sqrt z$
and $\log z$.

\begin{Proof}
In the~first step we consider the~function
\begin{equation}\label{7.6}
w=g_1(z):=\log z-\frac{cz}2- i\sqrt {cz}.
\qquad z\in\mathbb C^+.
\end{equation} 
Let us show that for every fixed $w\in\mathbb C^+$ such that $0<\IM w<\pi$
there exists an~unique $z\in\Bbb C^+$ such that $w=g_1(z)$. 
To prove this assertion we consider the~closed rectifiable curve 
$\gamma_7$, consisting of the~interval $\gamma_{7,1}:t,-R\le t\le-1/R$, 
connecting $-R$ to $-1/R$, 
the~arc $\gamma_{7,2}:e^{i\varphi}/R,0<\varphi<\pi$, 
connecting $-1/R$ to $1/R$, 
the~interval $\gamma_{7,3}:t,1/R\le t\le R$, connecting $1/R$ to $R$, and 
the~arc $\gamma_{7,4}:Re^{i\varphi},0<\varphi<\pi$, connecting $R$ to $-R$.   
The~parameter $R\ge 1$ will be chosen later on sufficiently large.

\begin{figure}[htbp]
 
\input{bild7.pstex_t}
 
\caption{}
\end{figure}

Let $z$ traverse $\gamma_7$ in the~counter clockwise direction.

If $z$ traverses $\gamma_{7,1}$ the~image $g_1(\gamma_{7,1})$ 
is contained in the~line $\IM \zeta=\pi$, connecting $\log R+cR/2+
\sqrt {cR}+i\pi$ to $-\log R+c/(2R)+\sqrt {c/R}+i\pi$.

If $z$ traverses $\gamma_{7,2}$ the~image $g_1(\gamma_{7,2})$ is contained 
in the~half-plane $\RE \zeta\le -\log R +1$,  
connecting $-\log R+c/(2R)+\sqrt{c/R}+i\pi$ to
$-\log R-c/(2R)-i\sqrt{c/R}$.

If $z$ traverses $\gamma_{7,3}$ the~image $g_1(\gamma_{7,3})$ is 
contained in the~half-plane $\IM \zeta<0$, connecting 
$-\log R-c/(2R)-i\sqrt{c/R}$ to $\log R-cR/2-i\sqrt {cR}$.

Finally, if $z$ traverses $\gamma_{7,4}$ the~image 
$g_1(\gamma_{7,4})$ is contained 
in the~domain $\{|\zeta|\ge cR/4,\, \IM \zeta\le \pi\}$, 
connecting the~points $\log R-cR/2-i\sqrt {cR}$ 
to $\log R+cR/2+\sqrt {cR}+i\pi$.

Therefore we can conclude that the~image $\zeta=g_1(z),z\in\gamma_7$, 
winds around $w$ once when $z$ runs through $\gamma_7$. Hence there 
is a~unique point $z$ inside the~curve $\gamma_7$ 
such that $w=g_1(z)$.  
Thus we have proved  
that for every fixed $w\in \mathbb S_{\pi}=\{w\in\mathbb C:0<\IM w<\pi\}$ 
there exists an~unique $z\in\mathbb C^+$ such that (\ref{7.6}) holds. 
Thus the~inverse
function $g_1^{(-1)}:\mathbb S_{\pi}\to \mathbb C^+$ exists and is analytic 
in $\mathbb S_{\pi}$. Introduce the~function $K_1(z):=g_1^{(-1)}(\log z),\,
z\in\mathbb C^+$. We note that
$K_1\in\mathcal N$ and $K_1^{(-1)}(z)=z\exp\{-cz/2-i\sqrt {cz}\}$
in the~domain of $\mathbb C^+$, where $K_1^{(-1)}$ is uniquely defined. 

Let us show that $K_1$ admits an~analytic continuation 
on $(-\infty,0)$ and its value on $(-\infty,0)$ is negative. 
Moreover $K_1(x)\to 0$ as $x\uparrow 0$. It is easy 
to~see that $(e^{g_1(x)})'>0$ for $x<0$. Since $e^{g_1(z)}$ is analytic 
on $(-\infty,0)$, we conclude that $(e^{g_1})^{(-1)}$ exists and 
is analytic on $(-\infty,0)$ as well. This function coincides 
for $z\in\mathbb C^+$ with the~function
$K_1(z)$ defined earlier. From the~definition of 
$e^{g_1(z)}$ it follows that $(e^{g_1})^{(-1)}(x)<0$ for $x<0$ and 
$(e^{g_1})^{(-1)}(x)\to 0$ as $x\uparrow 0$. By Corollary~\ref{3.3co}, 
$K_1(z)$ belongs to the~class $\mathcal K$. 
Therefore there exists $\mu_1\in\mathcal M_+$ such that
$K_{\mu_1}(z)=K_1(z)$ for $z\in\mathbb C^+$ and we proved 
the~assertion of the~lemma in the~case $j=1$.

Now we shall consider the~function 
\begin{equation}\label{7.7}
w=g_2(z):=\log z-\frac {cz}2+i\sqrt {cz}.
\qquad z\in\mathbb C^+.
\end{equation}
Let us show that for every fixed $w\in\mathbb C^+$ such that
$-2\pi<\IM w<-\pi$ there exists an~unique
$z\in\mathbb C^+$ such that $w=g_2(z)$. 

To prove this assertion we consider as before the~closed rectifiable 
curve $\gamma_7$ and as before $z$ traverses $\gamma_7$ in the~counter 
clockwise direction. 

If $z$ traverses $\gamma_{7,1}$ the~image $g_2(\gamma_{7,1})$ is
contained in the~line $\IM \zeta=\pi$, connecting $\log R+cR/2-
\sqrt {cR}+i\pi$ to $-\log R+c/(2R)-\sqrt {c/R}+i\pi$.

If $z$ traverses $\gamma_{7,2}$ the~image $g_2(\gamma_{7,2})$ is
contained in the~half-plane $\RE \zeta\le -\log R +1$, 
connecting $-\log R+c/(2R)-\sqrt{c/R}+i\pi$ to
$-\log R-c/(2R)+i\sqrt{c/R}$.

If $z$ traverses $\gamma_{7,3}$ the~image $g_2(\gamma_{7,3})$ is 
contained in the~half-plane $\IM \zeta>0$, connecting 
$-\log R-c/(2R)+i\sqrt{c/R}$ to $\log R-cR/2+i\sqrt {cR}$.

If $z$ traverses $\gamma_{7,4}$ the~image $g_2(\gamma_{7,4})$ is 
contained in the~domain $|\zeta|> cR/4$, $\Im\zeta<\pi,\,\Re \zeta\le 0$,
and $\Im\zeta<\sqrt{cR},\,\Re\zeta>0$, 
connecting  
$\log R-cR/2+i\sqrt {cR}$ to 
$\log R+cR/2-\sqrt {cR}+i\pi$.

Therefore we deduce that the~image $\zeta=g_2(z)$ winds around
$w$ once when $z$ runs through $\gamma_7$. Hence there is a~unique
point $z$ inside the~curve $\gamma_7$ such that $w=g_2(z)$. This 
relation holds for all sufficiently large $R>1$. Thus we have proved  
that for every fixed $w\in \mathbb S_{-\pi}:=\{w\in\mathbb C:-2\pi<\IM w<-\pi\}$ 
there exists an~unique $z\in\mathbb C^+$ such that (\ref{7.7}) holds. Then we 
deduce that the~inverse function $g_2^{(-1)}:\mathbb S_{-\pi}\to \mathbb C^+$ 
exists and is analytic in $\mathbb S_{-\pi}$. Introduce the~function
$K_2(z):=g_2^{(-1)}(\log z -2\pi i),\,z\in\mathbb C^+$. Obviously, 
$K_2^{(-1)}(z)=z\exp\{-cz/2+i\sqrt {cz}\}$, where $K_2^{(-1)}$
is well defined.  

As in the~case of $j=1$ we show show that $K_2(z)$ admits an~analytic 
continuation on $(-\infty,0)$ and its value on $(-\infty,0)$ is negative. 
Moreover $K_2(x)\to 0$ as $x\uparrow 0$. Hence, by Corollary~3.3, 
$K_2(z)$ belongs to the~class $\mathcal K$. 
Therefore there exists $\mu_2\in\mathcal M_+$ such that
$K_{\mu_2}(z)=K_2(z)$ for $z\in\mathbb C^+$ and we proved 
the~assertion of the~lemma in the~case $j=2$ as well.
\end{Proof}

\begin{lemma}
For every $c>0$, there exist 
$\mu_j\in\mathcal M_+,\,j=1,2$, such that
$$
\Sigma_{\mu_j}(z)=\exp\Big\{\frac c{2z}+(-1)^ji\sqrt {\frac cz}\Big\}
$$
for $z$, where $\Sigma_{\mu_j}(z)$ is defined.
\end{lemma}

The~proof of this lemma is similar to the~proof of Lemma~7.2
therefore we omit it.

Let the~functions $f_j(z),\,j=1,2$, are defined as in Lemma~7.1.
In addition we assume in (\ref{7.1}) that $0<a\le b<\infty$.
\begin{lemma}
For sufficiently small $\varepsilon>0$, there exist 
$\mu_j\in\mathcal M_+,\,j=1,2$, such that
$$
\Sigma_{\mu_j}(z)=\exp\{f_j(z)\}
$$
for $z$, where $\Sigma_{\mu_j}(z)$ is defined.
\end{lemma}
\begin{Proof}
Consider the~functions
$$
G(z):=\log z+\phi(z),\quad z\in\mathbb C^+.
$$
By Lemma~5.3, $G:\mathbb C^+\to\mathbb C$ takes every value in
$\mathbb S_{\pi}$ precisely once. Recall that $\mathbb S_{\pi}:=\{z\in\mathbb C:
\,0<\IM z<\pi\}$. For any $\eta\in(0,1/10)$ denote by $\Bbb S_{\eta}$
the~strip $\{z\in\mathbb C:\,\eta<\IM z <\pi-\eta\}$. Consider the~domains 
$\Omega_{\eta}:=G^{(-1)}(\mathbb S_{\eta}),\,\eta\in(0,1/10)$. 
Note that $\Omega_{\eta}$ is a~domain such that 
$\Omega_{\eta}\subset\mathbb C^+$.

The~boundary of $\Omega_{\eta}$ is a~curve $\gamma_8(\eta)$ consisting 
of curves $\gamma_{8,1}(\eta)$ and $\gamma_{8,2}(\eta)$
characterized by the~equations $\IM G(x+iy)=\pi-\eta,\,x<0,y>0$,
and $\IM G(x+iy)=\eta,\,x>0,y>0$, respectively. 
Since $G'(z)\ne 0$ for $z\in\Omega_{\eta}$, we conclude, by 
the~implicit function theorem, that $\gamma_{8,1}(\eta)$ and 
$\gamma_{8,2}(\eta)$, $\eta\in(0,1/10)$, are smooth Jordan curves.

We see, by the~definition of $\gamma_{8,1}(\eta)$, that this curve,
connecting $\infty$ to $0$, lies in the~angular region $\pi-2\eta<\arg z<\pi$.
By definition the~curve $\gamma_{8,2}(\eta)$, connecting $0$ to
$\infty$, lies in the~union of domains $|z-(a+b)/2|< b/2,\,\Im z>0$,
and $0<\arg z<2\eta$.
 
Denote by $\tilde\gamma_8(\eta)$ the~part of the~curve $\gamma_8(\eta)$,
lying in the~disk $|z-(a+b)/2|\le b/2$. Let us show that
\begin{equation}\label{7.8}
\sup_{0<\eta\le \frac 1{10}}\sup_{z\in\tilde\gamma_8(\eta)}
|\phi(z)|:=A_2<\infty.
\end{equation} 
Assume that $A_2=\infty$. Then for any
$N>1$ there exist $\eta_0\in(0,1/10]$ and $z_0=z_0(\eta_0)\in
\tilde\gamma_8(\eta_0)$ such that $0<\IM G(z_0)<\pi$ 
and $|G(z_0)|>N$. Note that $G(x)$ is real-valued for $x>0$ and $x\in\mathbb R\setminus[a,b]$.
Let $\IM z\ne 0$, then $\IM z\IM G(z)>0$ for $0<\RE z\le a/2$ and for $\RE z\ge b+a/2$.
In addition $\IM G(x)=\pi$ for $x<0$ and $\IM G(z)>\pi$ for $\RE z<0$ and $\IM z<0$. 
Moreover $|\phi(z)|/\log |z|\to 0$ as $|z|=R, 1/R$ 
and $R\to\infty$. Therefore we easily conclude, 
by Rouch\'e's theorem, that for $w\in\mathbb S_{\pi}$ 
with sufficiently large modulus 
the~equation $w=G(z)$ has a~solution $z=z(w)\in\mathbb S_{\pi}$ such that 
$z(w)=\exp\{w(1+o(1))\}$. Recalling that $G(z):\mathbb C^+\to\Bbb C$
takes every value in $\mathbb S_{\pi}$ precisely once, we conclude 
that $z_0=\exp\{G(z_0)(1+o(1))\}$,
a~contradiction for sufficiently large $N$. 

Now let us show for $j=1,2$ that 
$G_j(z):=\log z+f_j(z):\mathbb C^+\to\mathbb C$ takes every value in  
$\mathbb S_{\pi}$ precisely once in $G^{(-1)}(\mathbb S_{\pi})$.  
By the~formula
\begin{equation}\label{7.8a}
\Im G_j(z)=\arg z+\Im\phi(z)+\Im\phi(z)
(1+2(-1)^j\varepsilon\Re \phi(z))
\end{equation}
and by (\ref{7.8}), we see that $\IM G_j(z)\le \eta$, 
where $z\in \gamma_{8,2}(\eta),\,\eta\in(0,1/10)$, and $\varepsilon>0$ 
is sufficiently small. 

Using (\ref{7.1}), we note that, for $z\in \gamma_{8,1}(\eta),\eta
\in(0,1/10)$,
\begin{align}
-\IM \phi(z)&=\varepsilon_0\IM z\int\limits_{[a,b]}\frac{(1+u^2)\,\nu(du)}
{(u-\RE z)^2+(\IM z)^2}\le \varepsilon_0\frac{1+b^2}a
\nu([a,b])\frac{\IM z}{|\RE z|}
\notag\\
&\le 2\varepsilon_0\frac{1+b^2}a\nu([a,b])\eta.\notag 
\end{align}
Therefore we conclude from (\ref{7.8a}) that $\IM G_j(z)\ge \pi-2\eta$ for
$z\in \gamma_{8,1}(\eta),\eta\in(0,1/10)$, and sufficiently small
$\varepsilon_0,\varepsilon>0$. In addition we see from (\ref{7.1}) that for 
$|z|=1/R$ and $|z|=R$ the~inequality
$|\phi(z)|\le c$ holds, with a~positive constant $c$, 
depending on $\phi$ only, where $R>1$ is sufficiently large.  

Let $R>1$ be sufficiently large. Denote by $c'$ and $c''$ 
the~intersection  points of the~curve 
$\gamma_{8,1}(\eta)$ with the~circles $|z|=R$ and $|z|=1/R$, respectively. 
We denote as well by $d''$ and $d'$ the~intersection points of the~curve 
$\gamma_{8,2}(\eta)$ with the~circles $|z|=1/R$ and $|z|=R$, respectively.

Consider the~closed rectifiable curve $\gamma_9=\gamma_9(\eta),\,\eta
\in(0,1/10)$, consisting of $\gamma_{9,1}$:
the~part of $\gamma_{8,1}(\eta)$ lying in the~disc $|z|\le R$, connecting
$c'$ to $c''$, the~arc $\gamma_{9,2}:\,e^{i\theta}/R,\,
\arg d''<\theta<\arg c''$, connecting $c''$ to $d''$, the~curve 
$\gamma_{9,3}$: the~part of $\gamma_{8,2}(\eta)$ lying in the~disc 
$|z|\le R$, connecting $d''$ to $d'$, and the~arc 
$\gamma_{9,4}:\,Re^{i\theta},\,\arg d'<\theta<\arg c'$, 
connecting $d'$ to $c'$.

Assume that $w\in\mathbb S_{\pi}$ and $\eta<\frac 12\min\{\IM w,\pi-\IM w\}$. 
If $z$ traverses $\gamma_{9,1}$, the~image
$\zeta=G_j(z)$ lies in the~half-plane $\IM \zeta\ge \pi-4\eta$.
If $z$ traverses $\gamma_{9,2}$, it is easy to see that the~image 
$\zeta=G_j(z)$ lies in the~half-plane $\RE \zeta<-\frac 12\log R$.
If $z$ traverses $\gamma_{9,3}$, the~image $\zeta=G_j(z)$ lies in 
the~half-plane $\IM \zeta\le \eta$. Finally, if $z$ traverses $\gamma_{9,4}$, 
the~image $\zeta=G_j(z)$, lies in the~half-plane $\RE \zeta>\frac 12\log R$.

Summarizing, it follows that the~image
$\zeta=G_j(z)$ winds around $w$ once when $z$ traverses $\gamma_9$.
By the~argument principle, the~function $G_j(z)$ takes the~value
$w$ precisely once inside $\gamma_9$. Since this assertion holds 
for all sufficiently large $R>1$ and 
for all sufficiently small $\eta>0$, we obtain the~desired result.   

This result implies that the~inverse functions 
$\rho_j=G_j^{(-1)}:\mathbb S_{\pi}\to\mathbb C^+,\,j=1,2$,
exist and are analytic on $\mathbb S_{\pi}$. 

Let us show that $\rho_j(z)$ admit an~analytic continuation 
on the~half-line $\gamma_-:\IM z=\pi,\Re z<0$, and that their values on 
this half-line are negative.  
It is easy to~see that 
$$
G_j'(x)=\ffrac 1x+2\phi'(x)(1+(-1)^j\varepsilon\phi(x))<0,
\quad x<0,\,\,j=1,2.
$$ 
Since $G_j(z)$ is analytic 
on $(-\infty,0)$, we conclude that $G_j^{(-1)}$ exists and 
is analytic on $\gamma_-$ as well. Since, as shown above,
for every fixed $w\in\mathbb S_{\pi}$ there is  
an~unique point $z\in G^{(-1)}(\mathbb C^+)$ such that $w=\log z+f_j(z)$ holds,
this function coincides  for $z\in\mathbb S_{\pi}$ with the~function
$\rho_j(z)$ obtained early. Introduce the~function 
$\tilde \rho_j(z):=\rho_j(\log z),\,z\in\mathbb C^+$. Note that 
$\tilde \rho_j\in\mathcal N$ and $\tilde \rho_j^{(-1)}(z)=z\exp\{f_j(z)\}$, 
on the~domain $\mathbb C^+$, where $\tilde \rho_j^{(-1)}$ is uniquely defined. 
Moreover, the~function $\tilde \rho_j(z)$ admits
an~analytic continuation on $(-\infty,0)$.

From the~definition of 
$\tilde \rho_j(z)$ it follows that $\tilde \rho_j(x)<0$ for $x<0$ and 
$\tilde \rho_j(x)\to 0$ as $x\uparrow 0$. By Corollary~\ref{3.3co}, 
$\tilde \rho_j,\,j=1,2$, belong to the~class $\Cal K$ and  
we conclude the~assertion of the~lemma.
\end{Proof}

Now we complete the~proof of Theorem~2.13 for the~semigroup 
$(\mathcal M_+,\boxtimes)$.

At first note that the~measure $\mu\in\mathcal M_+$ with 
$\Sigma_{\mu}(z)=\exp\{-cz\}$, where $c>0$, does not belong to 
the~class $I_0$. Indeed, by Lemma~7.2, there exist measures $\mu_1$ 
and $\mu_2$ of $\mathcal M_+$ such that
$$
\Sigma_{\mu_j}(z)=\exp \{\tilde{f}_j(z)\}, \quad\text{where}\quad
\tilde{f}_j(z):=-cz/2+(-1)^ji\sqrt {cz},\quad j=1,2,
$$
for $z$, where $\Sigma_{\mu_j}(z)$ are defined. In addition 
$\tilde{f}_j(z)$ admit an~analytic continuation in $\Bbb C^+$
and $\IM\tilde{f}_j(z)$ have positive values
at some points of the~half-plane $\Bbb C^+$. Hence, by Bercovici and
Voiculescu result~\cite{BercVo:1992}, ~\cite{BercVo:1993}, $\mu_j,\,j=1,2$,
are not i.d. elements in the~semigroup 
$(\mathcal M_+,\boxtimes)$. On the~other hand $\Sigma_{\mu}(z)=
\Sigma_{\mu_1}(z)\Sigma_{\mu_2}(z)$ and we see that $\mu$ has non-i.d.
factors as was to be proved.
In a~similar way, using Lemma~7.3, we prove that the~measure 
$\mu\in\mathcal M_+$ with $\Sigma_{\mu}(z)=\exp\{c/z\}$, where $c>0$, 
does not belong to the~class $I_0$.

Now we consider the~case where $\Sigma_{\mu}(z)$ admits 
a~representation (\ref{2.21}) with $a=b=0$, $\nu(\{0\})=0$ and
$\nu((0,\infty))>0$.

For every fixed $0<a\le b$, we have the~representation
$$
\Sigma_{\mu}(z)=\Sigma_{\mu_1}(z)\Sigma_{\mu_2}(z)\Sigma_{\mu_3}(z),
$$
where
$$
\Sigma_{\mu_3}(z):=\exp\Big\{4\varepsilon_0\int\limits_{[a,b]}u\,\nu(du)
+(1-4\varepsilon_0)\int\limits_{[a,b]}\frac {1+uz}{z-u}\,
\nu(du)\Big\}
+\int\limits_{\mathbb R\setminus[a,b]}\frac {1+uz}{z-u}\,\nu(du)
$$
and $\Sigma_{\mu_j}(z),\,j=1,2$, are defined in Lemma~7.4. By Bercovici
and Voiculescu result (see Section~2), the~p-measure $\mu_3$ is i.d. for
$\varepsilon_0\le 1/4$.

If $\nu(\{a\})>0$, then we assume in (\ref{7.1}) $a=b$. If $\nu(\{a\})=0$
for all $a\in (0,\infty)$, then we choose the~points $0<a<b$ 
such that $\nu ([a,a+h])\,\,>ch$ and $\nu ([b-h,b])>ch$
for all $0<h\le h_0$, where $c>0$ and $h_0>0$ depend on 
the~measure $\nu$ only. Such points exist by Proposition~\ref{3.6pro}.
We showed in Section~7.1 that under these assumptions the~function 
$\phi(z)$, see (\ref{7.1}), has the~limiting behaviour (\ref{7.5}).

Note that the~measures $\mu_j,\,j=1,2$, are not i.d. 
Indeed, by (\ref{7.5}), there exist points in $\Bbb C^+$ 
at which the functions $\Im f_1(z)$ and $\Im f_2(z)$ 
have positive values and the~desired assertion follows from Bercovici
and Voiculescu result (see Section~2). Hence the~measure $\mu$ 
has a~non-i.d. factor and $\mu$ does not belong 
to the~class $I_0$.
$\square$

$\mathbf{7.3}$. It remains to prove Theorem~2.13 for the~semigroup 
$(\Cal M_*,\boxtimes)$. We need the~following auxiliary result.
\begin{lemma}
For sufficiently small $\varepsilon>0$, there exist 
$\mu_j\in\mathcal M_*,\,j=1,2$, such that, for $z$, where $\Sigma_{\mu_j}(z)$
is defined,
$$
\Sigma_{\mu_j}(z)=\exp\{q_j(z)\}.
$$
Here $q_j(z):=2q(z)+(-1)^ji\varepsilon q^2(z),\,j=1,2$, and,
for $\xi_0\in\mathbb T$ and $0<\Delta <1/100$,
\begin{equation}\label{7.9}
q(z):=\varepsilon_0\int\limits_{\{\xi\in\mathbb T:\,|\xi-\xi_0|\le \Delta\}}
\ffrac{\xi+z}{\xi-z}\,\nu(d\xi),\quad z\in\Bbb D.
\end{equation}
In $(\ref{7.9})$ $\varepsilon_0>0$ is a~sufficiently small constant,
depending on $\nu$ and $\Delta$ only, and 
$\nu(\{\xi\in\mathbb T:\,|\xi-\xi_0|\le \Delta\})>0$. 
\end{lemma}
\begin{Proof}
Let us prove the~lemma in the~case $j=1$. One can prove the~lemma
in the~case $j=2$ in the~same way.
Without loss of generality we assume that in the~definition of
the~function $q(z)$ the~parameter $\xi_0$ is equal to 1. 

Consider the~function
$$
Q(z):=\log z+q(z),\quad z\in \tilde{\mathbb D}:=\mathbb D
\setminus[-1,0].
$$
This function is analytic on $\tilde{\mathbb D}$.

Denote by $\gamma_{10,1}$, $\gamma_{10,2}$, and $\gamma_{10,3}$
the~Jordan curves defined by the~parametric equations $\zeta=\log(1-t)+
i\pi+q(t-1),\,0\le t< 1$, $\zeta=it+q(-1),\,
-\pi\le t\le\pi$, and $\zeta=\log t-i\pi+q(-t),\,0<t\le 1$, 
respectively. Note that $\RE q(-1)=0$.
Define a~curvilinear half-strip $S_-$ as a~domain with boundary
$\gamma_{10,1},\gamma_{10,2}$, and $\gamma_{10,3}$. Let us show that
$Q:\tilde{\mathbb D}\to\mathbb C$ takes every value in $S_-$ precisely once.  
For this we need to prove that for any $w\in S_-$ there exists 
an~unique point $z\in\tilde{\mathbb D}$ such that $w=Q(z)$, 
provided that $\varepsilon_0>0$ is sufficiently small.

Consider the~contour $\gamma_3$ with parameter $\theta=\pi$ 
(see Section~5). Choose $r\in(0,1)$ such that $1-r$ is small.
If $z$ traverses
$\gamma_{3,1}$ the~image $\zeta=Q(z)$ lies on the~curve $\gamma_{10,1}$.
If $z$ traverses $\gamma_{3,2}$ the~image $\zeta=Q(z)$ lies in 
the~half-plane $\RE\zeta<\frac 12\log(1-r)$. If $z$ traverses
$\gamma_{3,3}$ the~image $\zeta=Q(z)$ lies on the~curve $\gamma_{10,3}$.
Finally, by the~inequality $\RE q(z)\ge 0,\,z\in\mathbb D$, if $z$ traverses 
$\gamma_{3,4}$ the~image $\zeta=Q(z)$ lies
in the~half-plane $\RE\zeta>\log r$.  

Hence, the~image
$\zeta=Q(z)$ winds around $w$ once when $z$ runs through $\gamma_3$
in the~counter clockwise direction.
By the~argument principle, the~function $Q(z)$ takes the~value
$w\in S_-$ precisely once inside $\gamma_3$ 
for sufficiently small $1-r$ and the~desired result is proved. 
Hence the~inverse function $Q^{(-1)}:S_-\to\tilde {\mathbb D}$ is defined and
is analytic function on $S_-$.

Let $\gamma_{10,2}(\eta)$ is the~half-open interval of 
the~vertical line $\RE z=-\eta,\eta\in(0,1/100)$,
lying between the~curves $\gamma_{10,1}$ and $\gamma_{10,3}$.
Denote by $\gamma_{11}(\eta)$ the~closed rectifiable Jordan curve $Q^{(-1)}
(\gamma_{10,2}(\eta))$.  
By the~definition of $\gamma_{11}(\eta)$, we note that 
$\RE Q(z)=-\eta$ for 
$z\in \gamma_{11}(\eta),\,\eta\in(0,1/100)$. 
Since $Q:\tilde{\mathbb D}\to\mathbb C$ takes eve\-ry value in $S_-$ precisely 
once, we note that if $z$ traverses the~curve $\gamma_{11}(\eta)$
in the~counter clockwise direction, then $\Im Q(z)$ is a~monotone
function such that $-\pi+\alpha(-\eta)\le \Im Q(z)\le\pi+\alpha(-\eta)$,
where, as it is easy to see, $|\alpha(-\eta)-q(-1)|\le 1/10$ 
for $\eta\in(0,1/100)$.   
 
Let us prove
\begin{equation}\label{7.10}
\sup_{0<\eta\le \frac 1{100}}\sup_{z\in\gamma_{11}(\eta)}|q(z)|\le 1/3.
\end{equation}  
We shall assume that $\eta\in(0,1/100)$.
Note that $|q(z)|\le 1/10$ for $z\in\gamma_{11}(\eta)
\cap\{|z-1|\ge 2\Delta\}$.
From the~relation $\RE Q(z)=-\eta,\,z\in \gamma_{11}(\eta)$, and
the~inequality $\log|z|\ge \log (1-2\Delta)$ 
for $z\in\gamma_{11}(\eta)\cap\{|z-1|<2\Delta\}$  
we see that the~bound $|\Re q(z)|\le 1/10$ holds for such $z$.
In addition it is not difficult to see that the~inequality 
$|Q(z')|\le 1/10$ holds at points of intersection $z'$ of
$\gamma_{11}(\eta)$ with the~circle $|z-1|=2\Delta$. Since $\Im Q(z),\,
z\in\gamma_{11}(\eta)$,
is a~monotone function, we have $|\Im Q(z)|\le 1/10$ and hence
$|\Im q(z)|\le -\log (1-2\Delta)+\arctan (2\Delta)+1/10\le 1/5$ for
$z\in\gamma_{11}(\eta)\cap\{|z-1|<2\Delta\}$. Thus (\ref{7.10}) is proved.

Then, by the~relation
$$
\log |z|+\RE q_1(z)=\RE Q(z)+\RE q(z)(1+2\varepsilon\IM q(z))
$$
and (\ref{7.10}), we conclude that 
\begin{equation}\label{7.11}
\log |z|+\RE q_1(z)\ge -\eta,\,z\in\gamma_{11}(\eta),\,
\eta\in(0,1/100).
\end{equation}

Now we shall prove that $\tilde{q}_1(z):=ze^{q_1(z)}:\mathbb D\to\mathbb C$ 
takes every value $w\in\mathbb D$ precisely once in $Q^{(-1)}(S_-)\cup(-1,0]$.

Fix $w\in\mathbb D$. Let $\eta>0$ be sufficiently small.
If $z$ traverses $\gamma_{11}(\eta)$ in the~counter clockwise direction 
we see, by (\ref{7.10}) and (\ref{7.11}), that $|Arg e^{q_1(z)}|\le
1<\pi$ and $|e^{q_1(z)}|\ge e^{-\eta}$, respectively. 
By the~argument principle, the~function $ze^{q_1(z)}$ 
takes the~value $w\in\mathbb D$ precisely once inside 
$\gamma_{11}(\eta)$ for sufficiently small $\eta>0$. 
We obtain the~desired result.

Hence  the~inverse function $\tilde{q}_1^{(-1)}:\mathbb D\to\mathbb D$ thus defined
is analytic on $\mathbb D$ and, as it is easy to see, 
belongs to the~class $\mathcal S_*$. The~lemma is proved.
\end{Proof}

Let us complete the~proof of Theorem~2.13 for the~semigroup
$(\mathcal M_*,\boxtimes)$. Without loss of generality, we assume that
$\nu(A)>0$, where $A:=\{\xi\in\mathbb T,|\xi-1|\le\Delta\}=
\{\xi\in\mathbb T,-\Delta_1\le\arg\xi\le\Delta_1\}$
and write
$$
\Sigma_{\mu}(z)=\Sigma_{\mu_1}(z)\Sigma_{\mu_2}(z)\Sigma_{\mu_3}(z),
$$
where
$$
\Sigma_{\mu_3}(z):=(1-4\varepsilon_0)\int\limits_{A}\frac{\xi+z}
{\xi-z}\,\nu(d\xi)+
\int\limits_{\Bbb T\setminus A}\frac{\xi+z}{\xi-z}\,\nu(d\xi)
$$
and $\Sigma_j(z),\,j=1,2$, are defined in Lemma~7.5. If $\nu(\{1\}>0$,
we assume $A:=\{1\}$. If $\nu(\{\xi\})=0$ for all $\xi\in\mathbb T$, 
we choose the~point $\Delta_1>0$  
so that 
$\nu(\{\xi\in\mathbb T,\Delta_1-h\le\arg\xi\le\Delta_1\}>ch$ for all $0<h\le h_0$, 
where
the~constant $c>0$ and $h_0>0$ depend on the~measure $\nu$ only.
By Proposition~3.6, such points exist. Repeating the~arguments 
used in the~proof of (\ref{7.5}), we obtain
\begin{equation}\label{7.12}
\IM q(e^{i(h+\Delta_1)}(1-h))\to+\infty, \qquad h\downarrow 0.
\end{equation}
Note that the~measures $\mu_j,j=1,2$ are not i.d. Indeed, $q_2(z)$ is
an~analytic function in $\mathbb D$ and,
by (\ref{7.12}), there exist points $z\in\mathbb D$ where
$\RE q_2(z)<0$.
Hence, by Bercovici and Voiculescu result~\cite{BercVo:1992}, 
the~measure $\mu$ has a~non-i.d. factor
and therefore $\mu\notin I_0$.
$\square$

Now we shall prove that there exists a~measure $\mu\in\mathbf M$ for which
the~representation (\ref{2.27}) is not unique. We establish this result 
for the~semigroup $(\mathcal M,\boxplus)$. One can prove this fact 
for the~semigroups $(\mathcal M_+,\boxtimes)$ and $(\mathcal M_*,\boxtimes)$ 
in the~same way.

Assume to the~contrary that every measure $\mu\in\mathcal M$, which has
indecomposable factors, admits an~unique representation
\begin{equation}\label{7.13}
\mu=\mu_1\boxplus\mu_2\dots,
\end{equation}
where $\mu_1,\mu_2,\dots$ are some indecomposable nondegenerate 
elements of the~semigroup $(\mathcal M,\boxplus)$, with respect 
to the~equivalence relation $\mu\sim\nu$ if 
$\mu=\nu\boxplus\delta_a$ for some $a\in\mathbb R$.

Let $\mu$ be an~i.d. p-measure. Hence, for every
$n\in\mathbb N$, $\mu=\nu_n\boxplus\dots\boxplus\nu_n$ ($n$ times),
where $\nu_n\in\Cal M$.
Since the~p-measure $\nu_n$ admits an~unique representation of 
the~form (\ref{7.13}), we see that, for every $n\in\mathbb N$, 
$\mu=\mu_1^{n\boxplus}\boxplus\rho_n$,
where $\mu_1$ denotes the~measure from the~representation (\ref{7.13}) 
and $\rho_n\in\mathcal M$. We
return to the~notation of Section~6. Note that the~measure $\mu$ belongs
to the~set $\mathcal M^{(\alpha,\beta)}$ for some $\alpha>0$ and
$\beta>0$. Hence $0<-\IM\varphi_{\mu}(i(\beta+1))<\infty$. By 
Proposition~\ref{6.1pr}, $\mu_1^{n\boxplus},\,\rho_n\in\mathcal M^{(\alpha,\beta)}$
as well and $-\IM\varphi_{\mu_1^{n\boxplus}}(i(\beta+1))\ge 0$ 
and  $-\IM\varphi_{\rho_n}(i(\beta+1))\ge 0$. 
Since 
\begin{align}
\IM\varphi_{\mu}(i(\beta+1))&=\IM\varphi_{\mu_1^{n\boxplus}}
(i(\beta+1)))+\IM\varphi_{\rho_n}(i(\beta+1))\notag\\
&=n\IM\varphi_{\mu_1}(i(\beta+1)))+\IM\varphi_{\rho_n}(i(\beta+1)),\notag
\end{align}
we have
$\IM\varphi_{\mu_1}(i(\beta+1))\to 0$ as $n\to\infty$ and we get
$\IM\varphi_{\mu_1}(i(\beta+1))=0$. By Proposition~\ref{6.2apr}, 
$\mu_1=\delta_a$ for some $a\in\Bbb R$. We arrive at a~contradiction
which proves the~assertion.
$\square$ 

\section{Dense classes of indecomposable elements in 
$(\mathbf M,\circ)$}

In this section we describe a~wide class of indecomposable elements
in $(\mathbf M,\circ)$. Theorem~2.14 and Corollary~2.15 follow 
immediately from our results.

$\mathbf{8.1}$. At first we shall formulate and prove our result for the~semigroup
$(\mathcal M,\boxplus)$.

\begin{theorem}\label{8.1th}
Let $\mu\in\mathcal M$ such that $\mu(\{a\})>0$, $\mu([b,\infty))>0$ and 
$\mu(\{a\})+\mu([b,\infty))=1$, 
where $a<b$ are real numbers and $\lim_{x\uparrow b}G_{\mu}(x)=0$. 
Then $\mu$ is an~indecomposable p-measure in $(\mathcal M,\boxplus)$.
\end{theorem}{

From Proposition~\ref {3.6pro} it follows that there exist p-measures $\mu$ 
with an~unique atom $a$ satisfying conditions of Theorem~\ref{8.1th}
(see similar arguments in Section 7). Therefore the~statement of this theorem does not
follow from Belinshi~\cite{Bel:2006} and Bercovici and Wang~\cite{BercWa:2007} results.
This conclusion holds for Theorem~\ref{8.2th} and~\ref{8.3th} as well.

\begin{proof}
It follows from the~assumptions of the~theorem that $F_{\mu}(z)$
is an~analytic function on $\mathbb C\setminus [b,\infty)$ and $F_{\mu}(x)$ 
is strictly monotone on $(-\infty,b)$. In addition 
$F_{\mu}(x)\downarrow-\infty$ as $x\downarrow-\infty$ and
$F_{\mu}(x)\uparrow+\infty$ as $x\uparrow b$. 
Using the~Stieltjes-Perron inversion formula (\ref{3.5}), we see that $F_{\mu}(z)$
admits the~representation
\begin{equation}\label{8.0*}
F_{\mu}(z)=z+c+\int\limits_{[b,\infty)}\Big(\frac 1{t-z}-\frac t{1+t^2}
\Big)\,\sigma(dt),\quad z\in\mathbb C^+,
\end{equation}
where $c\in\mathbb R$ and $\sigma$ is a~nonnegative measure
such that $\int_{[b,\infty)}\sigma(dt)/(1+t^2)<\infty$. 

Let $\mu_j\in\mathcal M,\,j=1,2$, and $\mu=\mu_1\boxplus\mu_2$. 
Let us prove that either $\mu_1=\mu\boxplus\delta_{\alpha},\,
\mu_2=\delta_{-\alpha}$ 
or $\mu_1=\delta_{-\alpha},\,\mu_2=\mu\boxplus\delta_{\alpha}$ with 
$\alpha\in\mathbb R$.
By Theorem~2.1, there exist functions $Z_j(z)\in\mathcal F,j=1,2$, such that
(\ref{2.3}) holds and $F_{\mu}(z)=F_{\mu_1}(Z_1(z))=F_{\mu_2}(Z_2(z)),
\,z\in\mathbb C^+$. 
Using the~integral representation (\ref{3.3}) for Nevanlinna functions, 
we rewrite the~relation (\ref{8.0*}) in the~form
\begin{align}\label{8.0}
&z+c+\int\limits_{[b,\infty)}\Big(\frac 1{t-z}-\frac t{1+t^2}\Big)\,\sigma(dt) 
=b_j+Z_j(z)+W_j(z)\notag\\
&:=b_j+Z_j(z)+\int\limits_{\mathbb R}
\Big(\frac 1{t-Z_j(z)}-\frac t{1+t^2}\Big)\sigma_j(dt),
\quad z\in\mathbb C^+,\,\, j=1,2,
\end{align}
where $b_j\in\mathbb R$ and $\sigma_j$ are nonnegative measures such that
$\int_{\mathbb R} \sigma_j(dt)/(1+t^2)<\infty$. Since, again by
the~representation (\ref{3.3}) for functions in $\mathcal F$, we have
$$
Z_j(z)=z+c_j+\int\limits_{\mathbb R}\Big(\frac 1{t-z}
-\frac t{1+t^2}\Big)\nu_j(dt),\quad z\in\mathbb C^+,\,\,j=1,2,
$$
where $c_j\in\mathbb R$ and $\nu_j$ are nonnegative measures such that
$\int_{\mathbb R}\nu_j(dt)/(1+t^2)<\infty$. 
Moreover, by (\ref{3.4}), $|Z_j(iy)-iy|=o(y),\,j=1,2$, as $y\to+\infty$.
We note that the~functions $W_j(z)\in\mathcal N,\,j=1,2$, and, as it is easy to 
see, $|W_j(iy)|/y\to 0,\,j=1,2$, as $y\to+\infty$. Therefore, by (\ref{3.3})
and (\ref{3.4}), $W_j(z),\,j=1,2$, admit the~representation
$$
W_j(z)=\tilde{c}_j+\int\limits_{\mathbb R}\Big(\frac 1{t-z}
-\frac t{1+t^2}\Big)\tilde{\sigma}_j(dt),\quad z\in\mathbb C^+,
$$
where $\tilde{c}_j\in\mathbb R$ and $\tilde{\sigma}_j$ are nonnegative 
measures such that $\int_{\mathbb R}\tilde{\sigma}_j(dt)/(1+t^2)<\infty$.
Hence, we finally obtain,
for $ z\in\mathbb C^+$ and $j=1,2$,
\begin{equation}\label{8.1}
c+\int\limits_{[b,\infty)}\Big(\frac 1{t-z}-\frac t{1+t^2}\Big)\,\sigma(dt)
=b_j+c_j+\tilde{c}_j+\int\limits_{\mathbb R}
\Big(\frac 1{t-z}-\frac t{1+t^2}\Big)(\nu_j+\tilde{\sigma}_j)(dt).
\end{equation}
Applying the~Stieltjes-Perron inversion formula (\ref{3.5})
to (\ref{8.1}), we see that the~measures $\nu_j,j=1,2$, have
supports which are contained on the~set $[b,\infty)$. 
Hence either $Z_j(z)=z+c_j$ or
\begin{equation}\label{8.2}
Z_j(z)=z+c_j+\int\limits_{[b,\infty)}\Big(\frac 1{t-z}-\frac t{1+t^2}\Big)
\,\nu_j(dt) ,
\end{equation}
where $c_j\in\mathbb R$ and $\int_{[b,\infty)}\nu_j(dt)/(1+t^2)>0$.

Let one of $Z_j(z)$, say $Z_1(z)$, be of the~form $Z_1(z)=z+c_1$, then
$F_{\mu_1}(z)=F_{\mu}(z-c_1)$ and we have 
$\mu_1=\mu\boxplus\delta_{c_1}$ and $\varphi_{\mu_1}(z)=\varphi_{\mu}(z)+
c_1$. From the~relation (\ref{2.4}) $\varphi_{\mu}(z)=\varphi_{\mu_1}(z)+
\varphi_{\mu_2}(z)$, we obtain $\varphi_{\mu_2}(z)=-c_1$ which implies
$\mu_2=\delta_{-c_1}$.

It remains to consider the~case where both $Z_1(z)$ and $Z_2(z)$ have the~form
(\ref{8.2}). In this case $Z_1(z)$ and $Z_2(z)$ are analytic functions
on $\mathbb C\setminus [b,\infty)$. Moreover $Z_1'(x)>0$ and $Z_2'(x)>0$ for 
$x\in (-\infty,b)$, and $Z_1(x)\downarrow-\infty$ and $Z_2(x)\downarrow-\infty$ 
for $x\downarrow -\infty$. 

By the~relation (\ref{2.3}), we obtain
$$
z=Z_1(z)+Z_2(z)-F_{\mu_1}(Z_1(z))=Z_1(z)+Z_2(z)-F_{\mu}(z),
\quad z\in\mathbb C^+.
$$
Recalling that $F_{\mu}(x)\uparrow +\infty$ as $x\uparrow b$ we see
from this formula that one of $Z_j(z)$, say $Z_1(z)$, has the~property
$Z_1(x)\uparrow +\infty$ as $x\uparrow b$. 

Since $Z_1(z)\in\mathcal F$, the~function 
$1/(t-Z_1(z))$ is in $\mathcal N$ and $z/(t-Z_1(z))$ converges to $1$
as $z\to\infty$ for any fixed $t\in\mathbb R$. Moreover the~function 
$1/(t-Z_1(z))$ is analytic on $\mathbb C\setminus [b,\infty)$ with 
the~exception of a~simple pole $\beta=\beta(t)\in(-\infty,b)$.
Note that the~function $\beta(t)$, as a~function on the~variable $t$, 
is a~strictly increasing continuous function such that $\beta(t)\downarrow
-\infty$ as $t\downarrow -\infty$ and $\beta(t)\uparrow b$
as $t\uparrow \infty$.
In view of these properties 
we conclude that the~Nevanlinna integral representation (\ref{3.3}) 
for the~functions $1/(t-Z_1(z)),\,z\in\mathbb C^+$ for every fixed 
$t\in\mathbb R$, has the~form
\begin{equation}\label{8.3}
\frac 1{t-Z_1(z)}=\frac{\alpha(t)}{\beta(t)-z}+W(z,t)
:=\frac{\alpha(t)}{\beta(t)-z}+\int\limits_{[b,\infty)}\frac{\rho_{t}(du)}{u-z},
\end{equation}
where $\alpha(t)>0$ and $\rho_{t}$ is a~nonnegative finite measure.
Moreover $\alpha(t)$ is a~continuous function. We see from (\ref{8.3})
that $W(z,t)\in\mathcal N$ and $W(z,t)$ is bounded by modulus for $z$ from
every compact set in $\mathbb C^+$ and $t\in [-N,N]$ for every $N>0$.
From (\ref{8.0}) with $j=1$ and (\ref{8.3}) we obtain the~relation,
for every $a>0$,
\begin{equation}\label{8.3a}
z+c+\int\limits_{[b,\infty)}\Big(\frac 1{t-z}-\frac t{1+t^2}\Big)\,\sigma(dt) 
=Z_1(z)+\int\limits_{[-N,N]}\frac{\alpha(t)\,\sigma_1(dt)}{\beta(t)-z}
+T_N(z),\quad z\in\mathbb C^+,
\end{equation}
where $T_N(z)\in\mathcal N$.

Since $\alpha(t)>0,\,t\in\mathbb R$, and
$\beta(t)<b,\,t\in\mathbb R$, we note that the~second summand
on the~right hand-side of (\ref{8.3a}) can be written in the form
$$
\int\limits_{[-N,N]}\frac{\alpha(t)\,\sigma_1(dt)}
{\beta(t)-z}=\int\limits_{[t_1,t_2]}
\ffrac{\alpha(\beta^{(-1)}(u))\,\tilde{\sigma}_1(du)}{u-z},
$$ 
where $t_1:=\beta(-N)<b,\,t_2:=\beta(N)<b$, and 
$\tilde{\sigma}_1$ is a~measure such that $\tilde{\sigma}_1(B)
:=\sigma_1(\beta^{(-1)}(B))$ for any Borel set $B$. Note that
$\tilde{\sigma}_1([t_1,t_2])=\sigma_1([-N,N])$. Applying to
both sides of (\ref{8.3a}) the~inversion formula (\ref{3.4}),
we obtain that $\sigma_1([-N,N])=0$ for every $N>0$ and therefore
$\sigma_1(\mathbb R)=0$.
Hence the~relation (\ref{8.3a}) implies that 
$\sigma_1\equiv 0$. Thus, $F_{\mu_1}(z)=z+b_1$ and $\mu_1=\delta_{-b_1}$.
Hence $\mu_2=\mu\boxplus\delta_{b_1}$.
This implies the~statement of the~theorem.
\end{proof}

\vspace{0,5cm}
$\mathbf{8.2}$. We shall now formulate and prove the~result 
for the~semigroup $(\mathcal M_+,\boxtimes)$. 

\begin{theorem}\label{8.2th}
Let $\mu\in\mathcal M_+$ such that $\mu((0,a])>0$, $\mu(\{b\})>0$, and 
$\mu([0,a])+\mu(\{b\})=1$, where $0<a<b$, and $\lim_{x\uparrow 1/a}
\psi_{\mu}(x)\ge -1$. 
Then $\mu$ is an~indecomposable p-measure in $(\mathcal M_+,\boxtimes)$.
\end{theorem}
\begin{proof}
Denote $a_1:=1/a$ and $b_1:=1/b$. We have from the~assumptions of 
the~theorem that $\psi_{\mu}(z)$ is an~analytic function on $\mathbb C\setminus
[a_1,\infty)$ with the~exception of the~simple pole $b_1$. Moreover
$\psi_{\mu}(x)$ is a~strictly monotone function 
on $(-\infty,b_1)$ and on $(b_1,a_1)$,
and $\psi_{\mu}(x)\downarrow -1+\mu(\{0\})$ as $x\downarrow -\infty$,
$\psi_{\mu}(x)\uparrow \infty$ as $x\uparrow b_1$, $\psi_{\mu}(x)\downarrow
-\infty$ as $x\downarrow b_1$, and $\lim_{x\uparrow a_1}\psi_{\mu}(x)\ge -1$. 

Since $K_{\mu}(z):=1-1/(1+\psi_{\mu}(z))$, we conclude that $K_{\mu}(z)$
is an~analytic function on $\mathbb C\setminus [d,\infty)$, where
$b_1<d\le a_1$ is a~point such that $\lim_{x\uparrow d}\psi_{\mu}(x)=-1$.
In addition $K_{\mu}(x)$ is a~strictly monotone function on $(-\infty,d)$,
and $K_{\mu}(0)=0$ and $K_{\mu}(x)\uparrow \infty$ as $x\uparrow d$.
Since $K_{\mu}(z)\in\mathcal K$ it admits the~integral representation
\begin{equation}\label{8.4b}
K_{\mu}(z)=\gamma z+\int\limits_{[d,\infty)}\frac z{u-z}\,\tau(du),
\end{equation}
where $\gamma\ge 0$ and $\tau$ is a~nonnegative measure such that
$\int_{[d,\infty)}\tau(du)/(1+u)<\infty$.

Let $\mu_1,\mu_2$ denote measures in $\mathcal M_+$ and assume that 
$\mu:=\mu_1\boxtimes\mu_2$. We have to prove that either 
$\mu_1=\mu\boxtimes\delta_a$ and $\mu_2
=\delta_{1/a}$ or $\mu_1=\delta_{1/a}$ and $\mu_2=\mu\boxtimes\delta_a$, 
where $a\in (0,+\infty)$.

By Theorem~2.4, there exist functions $Z_j(z)\in\mathcal K,\,j=1,2$, such that
$K_{\mu}(z)=K_{\mu_1}(Z_1(z))=K_{\mu_2}(Z_2(z))$ for $z\in\mathbb C^+$
and
\begin{equation}\label{8.4a}
Z_1(z)Z_2(z)=zK_{\mu}(z),\quad z\in \mathbb C^+,\quad j=1,2.
\end{equation}
Since $K_{\mu_j}(z)\in\mathcal K,\,j=1,2$, rewrite the~first of these relations 
in the~form
\begin{equation}\label{8.4}
K_{\mu}(z)=d_jZ_j(z)+\int\limits_{(0,\infty)}\frac{Z_j(z)}{u-Z_j(z)}
\,\tau_j(du),\quad z\in\mathbb C^+,\,\,j=1,2,
\end{equation}
where $d_j\ge 0$ and $\tau_j$ are nonnegative measures such that
$\int_{(0,\infty)}\tau_j(du)/(1+u)<\infty$. 

Let us show that $Z_j(z),\,j=1,2$, are analytic functions on 
$\Bbb C\setminus [d,\infty)$. 

If $d_j\ne 0$ in (\ref{8.4}), then, applying Stieltjes inversion
formula (\ref{3.5}) to (\ref{8.4}), we note that $Z_j(z)$ is
an~analytic function on $\mathbb C\setminus [d,\infty)$. 

Let $d_j=0$.
Recalling the~definition of the~Krein class
$\mathcal K$, note that the~functions
$$
K_j(z;u)=\frac{Z_j(z)}{u-Z_j(z)},\qquad z\in\mathbb C^+,\quad j=1,2,
$$
belong to $\mathcal K$ for every fixed $u>0$. Therefore they admit 
the~representation
\begin{equation}\label{8.4d}
K_j(z;u)=a_j(u)z+z\int\limits_{(0,\infty)}\frac{\tau_{u,j}(dt)}{t-z},
\quad  z\in\mathbb C^+,\quad u>0,\quad j=1,2,
\end{equation} 
where $a_j(u)\ge 0$ and $\tau_{u,j}$ are nonnegative measures, satisfying
the~condition 
$$
\int\limits_{(0,\infty)}\frac{\tau_{u,j}(dt)}{(1+t)}<\infty.
$$
It is easy to see that $a_j(u)=0,\,u>0,\,j=1,2$. Moreover the~functions
$u\mapsto \tau_{u,j}(B)$, where $B$ is a~Borel set on $(0,\infty)$, 
are measurable. Using (\ref{8.4d}), we easily deduce from (\ref{8.4})
with $d_j=0$ the~relation
\begin{equation}\label{8.4e}
K_{\mu}(z)=\int\limits_{(0,\infty)}K_j(z;u)\,\tau_j(du)=
z\int\limits_{(0,\infty)}\frac{\nu_j(dt)}{t-z},\quad z\in\mathbb C^+,
\end{equation} 
where $\nu_j(B)=\int_{(0,\infty)}\tau_{u,j}(B)\,\tau_j(du)$ for every
Borel set on $\mathbb R_+$ and $\int_{(0,\infty)}\nu_{j}(dt)/(1+t)<\infty$.

By (\ref{8.4b}) and Stieltjes inversion
formula (\ref{3.5}), we deduce from (\ref{8.4e}) that $\nu_j((0,d))=0$.
Hence there exists a~measurable set $B_j$ such that $\int_{B_j}\tau_{j}(dt)
/(1+t)>0$ and  
$\tau_{u,j}((0,d))=0$ for $u\in B_j$. Thus we have from (\ref{8.4d})  
the~formula
$$
K_j(z;u)=z\int\limits_{[d,\infty)}\frac{\tau_{u,j}(dt)}{t-z}, 
\qquad z\in\mathbb C^+,\quad u\in B_j,\quad j=1,2.
$$
Since $Z_j(z)=uK_j(z;u)/(1+K_j(z;u))$ for $u\in B_j$, are analytic functions
on $\mathbb C\setminus [d,\infty)$, we proved the~desired assertion.

Hence $Z_j(z),\,j=1,2$, are
analytic functions on $\mathbb C\setminus [d,\infty)$ and $Z_j(0)=0$. Moreover
$Z_j(x),\,j=1,2$, are strictly monotone functions on $(-\infty,d)$. By
(\ref{8.4a}), we conclude that one of $Z_j(z)$, say $Z_1(z)$, has 
the~following property: $Z_1(x)\uparrow \infty$ as $x\uparrow d$.
Since $Z_1(z)\in\mathcal K$, the~function 
$Z_1(z)/(u-Z_1(z))$ is in $\mathcal K$ for every fixed $u\in(0,\infty)$. 
Moreover the~function 
$Z_1(z)/(u-Z_1(z))$ is analytic on $\mathbb C\setminus [d,\infty)$ with 
the~exception of a~simple pole $\beta=\beta(u)\in(0,d)$.
Note that the~function $\beta(u)$, as a~function of the~variable $u$, 
is strictly increasing continuous function such that $\beta(u)\downarrow
0$ as $u\downarrow 0$ and $\beta(u)\uparrow d$
as $u\uparrow \infty$.
In view of these properties 
we conclude that the~Nevanlinna integral representation (\ref{3.6}) 
for the~functions $Z_1(z)/(z(u-Z_1(z))),\,z\in\mathbb C^+$, for every fixed
$u\in(0,\infty)$, has the~form
\begin{equation}\label{8.4c}
\frac 1z \cdot\frac {Z_1(z)}{u-Z_1(z)}
=\frac{\alpha(u)}{\beta(u)-z}+R(z,u):=\frac{\alpha(u)}{\beta(u)-z}+p
+\int\limits_{[d,\infty)}
\frac{\rho_{u}(ds)}{s-z},
\end{equation}
where $p\ge 0,\,\alpha(u)>0$ and $\rho_{u}(ds)$ is a~nonnegative measure 
such that $\int_{[d,\infty)}\rho_{u}(ds)$ $/(1+s)<\infty$.
Moreover $\alpha(u)$ is a~continuous function. We see from (\ref{8.4c})
that $R(z,u)\in\mathcal N$ and $R(z,u)$ is bounded by modulus for $z$ from
every compact set in $\mathbb C^+$ and $u\in (0,N]$ for every $N>0$.
Recalling (\ref{8.4b}), (\ref{8.4}) with $j=1$ and (\ref{8.4c}), 
we obtain the~relation, for every $N>0$,
\begin{equation}\label{8.5a}
\gamma +\int\limits_{[d,\infty)}\frac {\tau(du)}{u-z} 
=d_1\frac{Z_1(z)}z+\int\limits_{(0,N]}\frac{\alpha(u)\,\tau_1(du)}{\beta(u)-z}
+T_N(z),\quad z\in\mathbb C^+,
\end{equation}
where $T_N(z)\in\mathcal N$.

Since $\alpha(u)>0,\,u\in (0,\infty)$, and
$0<\beta(u)<d,\,u\in (0,\infty)$, we note that the~second summand
on the~right hand-side of (\ref{8.5a}) can be written in the form
$$
\int\limits_{(0,N]}\frac{\alpha(u)\,\sigma_1(du)}
{\beta(u)-z}=\int\limits_{(0,u_1]}
\frac{\alpha(\beta^{(-1)}(u))\,\tilde{\tau}_1(du)}{u-z},
$$ 
where $u_1:=\beta(N)$, and 
$\tilde{\tau}_1$ is a~measure such that $\tilde{\tau}_1(B)
:=\tau_1(\beta^{(-1)}(B))$ for any Borel set $B$. Note that
$\tilde{\tau}_1((0,u_1])=\tau_1((0,N])$. Applying to
both sides of (\ref{8.5a}) the~inversion formula (\ref{3.4}),
we obtain that $\tau_1((0,N])=0$ for every $N>0$ and therefore
$\tau_1((0,\infty))=0$.
Hence relation (\ref{8.5a}) implies that 
$\tau_1\equiv 0$. Thus, $K_{\mu_1}(z)=d_1z$ and $\mu_1=\delta_{1/d_1}$.
Hence $\mu_2=\mu\boxplus\delta_{d_1}$.
This implies the~statement of the~theorem.
\end{proof}

\vspace{0,5cm}

$\mathbf{8.3}$. It remains to describe a~wide class of indecomposable elements
in the~semigroup $(\mathcal M_*,\boxtimes)$. Our result in this case 
has the~form.

Denote $\gamma_{\alpha}:=\{\zeta\in\mathbb T:-\alpha<\arg \zeta<\alpha\},\,
0<\alpha<\pi$, and $\zeta_{\alpha}=e^{i\alpha}$. 

\begin{theorem}\label{8.3th}
Let $\mu\in\mathcal M_*$ such that $\mu(\{1\})>0$, $\mu(\gamma_{\alpha}
\setminus\{1\})=0$, $\mu(\mathbb T\setminus \gamma_{\alpha})>0$, 
and $\lim_{\theta\uparrow \alpha}
\Im\psi_{\mu}(e^{i\theta})=-\infty$, $\lim_{\theta\downarrow -\alpha}
\Im\psi_{\mu}(e^{i\theta})=\infty$. 
Then $\mu$ is an~indecomposable p-measure in $(\mathcal M_*,\boxtimes)$.
\end{theorem}
\begin{proof}
We have from the~assumptions of 
the~theorem that $\psi_{\mu}(z)$ is an~analytic function on $\mathbb C\setminus
(\mathbb T\setminus\gamma_{\alpha})$ with the~exception of 
the~simple pole $1$. Moreover
$\Im\psi_{\mu}(e^{i\theta})$ is a~strictly monotone function 
on $(-\alpha,0)$ and on $(0,\alpha)$,
and $\Im\psi_{\mu}(e^{i\theta})\uparrow \infty$ as $\theta\downarrow -\alpha$,
$\Im\psi_{\mu}(e^{i\theta})\downarrow -\infty$ as $\theta\uparrow 0$, 
$\Im \psi_{\mu}(e^{i\theta})\uparrow\infty$ as $\theta\downarrow 0$, 
and $\Im\psi_{\mu}(e^{i\theta})\downarrow -\infty$
as $\theta\uparrow \alpha$. Introduce the~function 
$H_{\mu}(z)=1+2\psi_{\mu}(z)$ and consider the~function
$$
Q_{\mu}(z)=\frac{H_{\mu}(z)-1}{H_{\mu}(z)+1},\quad z\in\mathbb D.
$$
It is analytic on $\mathbb C\setminus(\mathbb T\setminus\gamma_{\alpha})$
and $|Q_{\mu}(z)|=1$ for $z\in\gamma_{\alpha}$.
It is easy to verify that ${\rm Arg} Q_{\mu}(e^{i\theta}),\,
{\rm Arg} Q_{\mu}(1)=0$, is strictly monotone
function on the~interval $-\alpha<\theta<\alpha$ and 
${\rm Arg} Q_{\mu}(e^{i(\alpha-0)})   
-{\rm Arg} Q_{\mu}(e^{i(-\alpha+0)})=4\pi$.

Let $\mu_1,\mu_2,\in\mathcal M_*$
We have to prove that 
either $\mu_1=\mu\boxtimes\delta_a$ and $\mu_2=\delta_{1/a}$ or
$\mu_1=\delta_{1/a}$ and $\mu_2=\mu\boxtimes\delta_a$, where $a\in \mathbb T$.
By Theorem~2.7, there exists a~function $Z_1(z)\in\mathcal S_*$ such that
$\psi_{\mu}(z)=\psi_{\mu_1}(Z_1(z))=\psi_{\mu_2}(Z_2(z))$ for $z\in\mathbb D$
and
\begin{equation}\label{8.11b}
Z_1(z)Z_2(z)=zQ_{\mu}(z)=zQ_{\mu_j}(Z_j(z))
=z\frac{\psi_{\mu_j}(Z_j(z))}{1+\psi_{\mu_j}(Z_j(z))}=z\frac{H_{\mu_j}(Z_j(z))-1}
{H_{\mu_j}(Z_j(z))+1}
\end{equation}
for $z\in\mathbb D,\quad j=1,2$, where $H_{\mu_j}(z):=2\psi_{\mu_j}(z)+1$.
We obtain from (\ref{8.11b}) the~relation
\begin{equation}\label{8.11}
\int\limits_{\mathbb T}\frac{1+\zeta z}{1-\zeta z}\,\mu(d\zeta)=
\int\limits_{\mathbb T}\frac{1+\zeta Z_j(z)}{1-\zeta Z_j(z)}\,\mu_j(d\zeta),
\quad z\in\mathbb D,\quad j=1,2.
\end{equation}
  
Since the~functions 
\begin{equation}\label{8.11c}
C_j(z;\xi)=\ffrac{1+Z_j(z)\xi}{1-Z_j(z)\xi},\qquad z\in\mathbb D,\quad j=1,2, 
\end{equation}
belong to the~Carath\'eodory class $\mathcal C$ for every fixed 
$\xi\in\mathbb T$, they admit the~representation
$$
C_j(z;\xi)=\int\limits_{\mathbb T}\frac{1+\zeta z}{1-\zeta z}
\,\sigma_{\xi,j}(d\zeta),\qquad z\in\mathbb D,\quad j=1,2,
$$
where $\sigma_{\xi,j}$ are p-measures. 
Moreover for $j=1,2$ the~functions
$\xi\mapsto\sigma_{\xi,j}(B)$, where $B$ is a~Borel set on $\Bbb T$, 
are measurable. 
From (\ref{8.11}) we deduce the~relations
\begin{equation}\label{8.12}
H_{\mu}(z)=\int\limits_{\mathbb T}
C_j(z;\xi)\,\mu_j(d\xi)=\int\limits_{\mathbb T}\ffrac{1+\zeta z}{1-\zeta z}
\,\nu_j(d\zeta),\quad z\in\mathbb D,
\end{equation}
where $\nu_j(B):=\int_{\mathbb T}\sigma_{\xi,j}(B)\mu_j(d\xi)$ for every 
Borel set $B$ on $\mathbb T$.

By (\ref{8.12}) and the~inversion formula, 
the~p-measures $\nu_j$ have the~property: $\nu_j(\gamma_{\alpha}
\setminus\{1\})=0,\,j=1,2$. Hence, there exist measurable sets 
$B_j,\,j=1,2$, such that $\mu_j(B_j)=1$ and 
$\sigma_{\xi,j}(\gamma_{\alpha}\setminus\{1\})=0$ for $\xi\in B_j$. 
Thus we obtain, for such points $\xi$,
\begin{equation}\label{8.14}
\frac{1+\xi Z_j(z)}{1-\xi Z_j(z)}=C_j(z;\xi)=p_{j}(\xi)\frac{1+z}{1-z}+
\int\limits_{\mathbb T\setminus \gamma_{\alpha}}
\frac{1+\zeta z}{1-\zeta z}\,\sigma_{\xi,_j}(d\zeta),\,\,z\in\mathbb D,\,\,j=1,2,
\end{equation}
where $0\le p_{j}(\xi)\le 1$. Moreovere there exist points $\xi_j\in B_j,
\,j=1,2$, such that $p_{j}(\xi_j)>0$. Since
$$
\xi_j Z_j(z)=\frac{C_j(z;\xi_j)-1}{C_j(z;\xi_j)+1},\quad j=1,2,
$$ 
we see from (\ref{8.14}) that the~functions $Z_j(z),\,j=1,2$, are
analytic on $\mathbb C\setminus(\mathbb T\setminus\gamma_{\alpha})$ 
and $|Z_j(z)|=1$ for $z\in\gamma_{\alpha}$. Hence $Z_j(e^{i\theta})
=e^{ig_j(\theta)},\,j=1,2$, for $-\alpha<\theta<\alpha$, where $g_j(\theta)$
are continuous real-valued functions.
In addition, as it is easy to see, $g_j(\theta),\,j=1,2$,
are strictly monotone functions on $(-\alpha,\alpha)$ and $Z_j(1)=1/\xi_j$. 

We note from (\ref{8.11b}) that
\begin{align}\label{8.15}
g_1(\alpha-0)&-g_1(-\alpha+0)+g_2(\alpha-0)
-g_2(-\alpha+0)\notag\\
&=2\alpha+\arg Q_{\mu}(e^{i(\alpha-0)})-\arg Q_{\mu}(e^{i(-\alpha+0)})=2\alpha+4\pi.
\end{align}
We conclude from (\ref{8.15}) that one of $Z_j$, say $Z_1(z)$, has 
the~property: $g_1(\alpha-0)-g_1(\alpha+0)>2\pi$.

Return to the~formula (\ref{8.12}) with $j=1$. We obtain from the~definition
(\ref{8.11c}) of the~functions $C_1(z;e^{i\psi}),\,
z\in\mathbb D,\,-\pi\le\psi<\pi$, 
that they have the~form
\begin{equation}\label{8.16}
C_1(z;e^{i\psi})=\alpha(\psi)\frac{1+e^{i\beta(\psi)}z}
{1-e^{i\beta(\psi)}z}+M(z;\psi),
\end{equation}
where $\alpha(\psi)$ is a~positive continuous function, 
$\beta(\psi)$ is a~strictly monotone continuous function
on $[-\pi,\pi)$ such that $-\alpha<\beta(\psi)<\alpha$, 
$\beta(0)=-\arg \xi_1$, and,
for every fixed $\psi\in[-\pi,\pi)$, $M(z;\psi)$ is an~analytic function
on $\Bbb C\setminus (\mathbb T\setminus\gamma_{\alpha})$.
Comparing (\ref{8.14}) 
and (\ref{8.16}) we arrive 
at a~contradiction.
Theorem~\ref{8.3th} is proved completely.

\end{proof}

\bigskip
          
\address{
Gennadii Chistyakov\\
Fakult\"at f\"ur Mathematik\\
Universit\"at Bielefeld\\
Postfach 100131\\
33501 Bielefeld\\
Germany
}

\email{chistyak@math.uni-bielefeld.de, chistyakov@ilt.kharkov.ua}

\bigskip

\address{
Friedrich G\"otze\\ 
Fakult\"at f\"ur Mathematik\\
Universit\"at Bielefeld\\
Postfach 100131\\
33501 Bielefeld\\
Germany
}

\email{ goetze@math.uni-bielefeld.de}

\end{document}